\documentclass[aos,preprint]{imsart}

\usepackage{amsthm,amsmath,mathrsfs,amssymb,graphicx, epsfig}
\usepackage{graphicx,array,epsfig,fancyheadings,rotating,amsfonts,subfig,latexsym,enumerate,paralist,etoolbox,mathtools,color,dsfont}
\usepackage[dvipsnames]{xcolor}

\RequirePackage[colorlinks,citecolor=blue,urlcolor=blue]{hyperref}

\startlocaldefs

\newcommand{\diag}{\textup{diag}}
\newcommand{\offdiag}{\overline{\textup{diag}}}
\newcommand{\corr}{\textup{Corr}}
\newcommand{\cov}{\textup{Cov}}

\newcommand{\bd}{\boldsymbol}
\newcommand{\mb}{\mathbf}
\newcommand{\be}{\begin{equation}}
 \newcommand{\ee}{\end{equation}}
 \DeclareMathOperator*{\argmin}{arg\,min}
\newtheorem{thm}{Theorem}
\newtheorem{lem}{Lemma}
\newtheorem{cor}{Corollary}

\theoremstyle{definition}

\newtheorem{example}{Example}
\newtheorem{case}{Case}
\newtheorem{remark}{Remark}
\newtheorem{prop}{Proposition}
\endlocaldefs

\usepackage[font={footnotesize}]{caption}

\makeatletter
\newcommand{\labitem}[2]{%
\def\@itemlabel{#1}
\item
\def\@currentlabel{#1}\label{#2}}
\makeatother

\begin{document}

\begin{frontmatter}

\title{ESTIMATION OF LARGE COVARIANCE AND PRECISION MATRICES FROM TEMPORALLY DEPENDENT OBSERVATIONS}
\runtitle{LARGE MATRIX ESTIMATION UNDER DEPENDENCE}


\begin{aug}
  \author{\fnms{Hai}  \snm{Shu}\thanksref{t1}
  \ead[label=e1]{haishu@umich.edu}}
  \and
  \author{\fnms{Bin} \snm{Nan}\ead[label=e2]{bnan@umich.edu}}\thanksref{t1}

\thankstext{t1}{Supported in part by NIH grant R01-AG036802 and NSF grant DMS-1407142. }

 \runauthor{H. Shu and B. Nan}

  \affiliation{University of Michigan}

  \address{
  Department of Biostatistics\\
  University of Michigan\\
  1415 Washington Heights\\
Ann Arbor, Michigan 48109-2029\\
USA
\\
          \printead{e1}\\
          \phantom{E-mail:\ }\printead*{e2}}


\end{aug}

\begin{abstract}
We consider the estimation of large covariance and precision matrices from high-dimensional sub-Gaussian or heavier-tailed observations with slowly decaying temporal dependence.
The temporal dependence is allowed to be long-range so with longer memory than those considered in the current literature. We show that several commonly used methods for independent observations can be applied to the temporally dependent data. In particular, the rates of convergence are obtained for the generalized thresholding estimation of covariance and correlation matrices, and for the constrained $\ell_1$ minimization and the $\ell_1$ penalized likelihood estimation of precision matrix.
Properties of sparsistency and sign-consistency are also established.
A gap-block cross-validation method is proposed for the tuning parameter selection, which performs well in simulations. As a motivating example,
we study the brain functional connectivity using resting-state fMRI time series data with long-range temporal dependence.

\end{abstract}



\end{frontmatter}


\section{Introduction}\label{sec: intro}

Let $\{\bd{X}_1,\dots,\bd{X}_n\}$ be a sample of $p$-dimensional random vectors, each with the same mean $\bd{\mu}_p$, covariance matrix $\mb{\Sigma}$ and precision matrix $\mb{\Omega}=\mb{\Sigma}^{-1}$. It is well known that the sample covariance matrix
is not a consistent estimator of $\mb{\Sigma}$ when
$p$ grows with $n$ \cite{Bai10,Bai93}.
When the sample observations $\bd{X}_1,\dots,\bd{X}_n$ are independent and identically distributed (i.i.d.), several regularization methods have been proposed for the consistent estimation of large $\mb{\Sigma}$, including
thresholding \cite{Bick08a,Cai11,El08,Roth09}, block-thresholding \cite{Cai12},
banding
\cite{Bick08b} and
tapering \cite{Cai10}.
Existing methods also include Cholesky-based method
\cite{Huan06,Roth10},
penalized pseudo-likelihood method
\cite{Lam09}
and
sparse matrix transform
\cite{Cao11}.
Consistent correlation matrix estimation can be obtained similarly from i.i.d. observations
\cite{El08,Jian04}.

The precision matrix $\mb{\Omega}=(\omega_{ij})_{p\times p}$, when it exists,
is closely related to
the partial correlations between the pairs of variables in a vector $\bd{X}$. Specifically, the partial correlation between $X_i$ and $X_j$ given $\{X_k,k\ne i,j\}$ is  equal to $-\omega_{ij}/\sqrt{\omega_{ii}\omega_{jj}}$ \cite{Cram46}.
Zero partial correlation means conditional independence
between Gaussian or nonparanormal random variables \cite{Liu09}.
There is a rich literature on the estimation of large $\mb{\Omega}$ from i.i.d. observations.
Various algorithms for the $\ell_1$ penalized maximum likelihood method ($\ell_1$-MLE) and its variants have been developed \cite{Bane08,Frie08,Hsie14,Yuan07}, and related theoretical properties have been investigated by \cite{Lam09,Ravi11,Roth08}.
Methods of estimating $\mb{\Omega}$ column-by-column thus implementable with parallel computing include nodewise Lasso \cite{Mein06, VGeer14},
graphical Dantzig selector \cite{Yuan10}, constrained $\ell_1$-minimization for inverse matrix estimation (CLIME) 
\cite{CLL11}, and adaptive CLIME \cite{CLZ16}.

Recently, researchers become increasingly interested in estimating the large covariance and precision matrices from temporally dependent observations $\{\bd{X}_t: t =1, \dots, n\}$, here $t$ denotes time.
Such research is particularly useful in analyzing the resting-state functional magnetic resonance imaging (rfMRI) data
to assess the brain functional connectivity \cite{Powe11,Ryal12}. In such imaging studies, the number of brain nodes (voxels or regions of interest) $p$ can be greater than the number of images $n$.
The temporal dependence of time series $\bd{X}_t$ is traditionally dealt with by imposing
the so-called strong mixing conditions \cite{Brad05}. To overcome the difficulties in computing strong mixing coefficients and verifying strong mixing conditions, \cite{Wu05} introduced a new type of dependence measure, the functional dependence measure, and recently applied it to the hard thresholding estimation of large covariance matrix and the $\ell_1$-MLE type methods of large precision matrix
\cite{Chen13}. The functional dependence measure 
may still be difficult
 to understand and to interpret for most data analysts. Practically, it is straightforward to describe the temporal dependence directly by using cross-correlations \cite{Broc91}.
By imposing certain weak dependence conditions on the cross-correlation matrix of samples $\{\bd{X}_t\}_{t=1}^n$,
\cite{Bhat14, Bhat14_2}
extended the banding and tapering regularization methods for estimating covariance matrix.

A univariate stationary time series is said to be long-memory if its  autocorrelation function $\rho(t)$ satisfies $\sum_{t=0}^\infty|\rho(t)|=\infty$, and short-memory otherwise \cite{Palm07}. The rfMRI data have been reported with long-memory in the scientific literature, e.g., \cite{He11,Tagl13}. However, the temporal dependence considered by \cite{Chen13} and that considered by \cite{Bhat14,Bhat14_2} do not cover any long-memory time series. Later we will illustrate that the rfMRI data example does not meet their restrictive temporal dependence conditions. Hence, it is important to show the applicability of the estimating methods for i.i.d. samples to this kind of data. In this article, we propose a new and straightforward temporal dependence measure that solely depends on the Frobenius norm and the spectral norm of the autocorrelation matrices of $\{\bd{X}_t\}_{t=1}^n$. Such a new measure clearly displays the effect of temporal dependence on the convergence rates of our considered matrix estimators, allowing each time series to be long-memory or even to be non-stationary. So the rfMRI data can be well handled by our relaxed assumption (see Figure~\ref{autocorrelation}). To the best of our knowledge, this is the first work that investigates the estimation of large covariance and precision matrices 
from long-memory observations.

Note that the estimation of large correlation matrix was not considered by either \cite{Chen13} or \cite{Bhat14, Bhat14_2}, which is a more interesting problem in, for example,  the study of brain functional connectivity. It was considered in a recent work by \cite{Zhou14} but under the assumption that all $p$ time series have the same temporal decay rate, which is rather restrictive and often violated (see Figure~\ref{autocorrelation} for an example of rfMRI data). Moreover, all four aforementioned articles assumed that $\bd{\mu}_p=(\mu_{pi})_{1\le i\le p}$ is known, which may not be true in practice.  Although the sample mean $\bar{X}_i=\frac{1}{n}\sum_{j=1}^nX_{ij}$ entrywise converges to $\mu_{pi}$ in probability or even almost surely under some dependence conditions \cite{Broc91,Hu08}, extra care will still be needed when true mean is replaced by sample mean in the matrix estimation, especially for long-memory, heavy-tailed, or even non-stationary data. We consider unknown $\bd{\mu}_p$ in this article, and show that the mean estimation indeed affects our derived matrix convergence rates, particularly for data with heavy tail probabilities.

In this article, we study the generalized thresholding estimation \cite{Roth09} for covariance and correlation matrices, and the CLIME approach \cite{CLL11} and an $\ell_1$-MLE type method called sparse permutation invariant covariance
estimation (SPICE; \cite{Roth08})
for precision matrix.
The convergence rates, sparsistency and sign-consistency are provided for temporally dependent data, potentially with long memory, which are generated from 
 a linear spatio-temporal model with all basis random variables coming from sub-Gaussian, or generalized sub-exponential, or distributions with polynomial-type tails.
We also establish the minimax optimal convergence rates of estimating covariance and correlation matrices for a certain class of temporally dependent sub-Gaussian data, including short-memory and some long-memory cases, and show that they can be achieved by the generalized thresholding method. Moreover, if the matrix $\ell_1$ norm of the precision matrix is bounded by a constant for such data, then the CLIME estimator attains the minimax optimal rates for i.i.d. observations shown in \cite{CLZ16}.

The article is organized as follows.
In Section~\ref{s:Temporal dependence},
we introduce the new temporal dependence measure and the considered temporally dependent data generating mechanism.
We provide the theoretical results for 
the estimation of covariance and correlation matrices in Section~\ref{s:cov est} and for the estimation of precision matrix in Section~\ref{s:prec est} for temporally dependent observations with sub-Gaussian tails. We consider extensions to data with generalized sub-exponential tails and polynomial-type tails in Section~\ref{s: heavy tails}.
In Section~\ref{s:numeric}, we introduce a gap-block cross-validation
method for the tuning parameter selection, evaluate the estimating performance via simulations,
and analyze a rfMRI data set for brain functional connectivity.	
The concentration inequalities that the proofs of the theoretical results are based on are given in the Appendix. Detailed proofs together with additional numerical considerations are provided in the Supplementary Material due to the page limitation.

\section{Temporal dependence}\label{s:Temporal dependence}

We start with a brief introduction of useful notation.
For a real matrix $\mb{M}=(M_{ij})$,  
we define: the 
spectral norm $\| \mb{M} \|_2=[\varphi_{\max}{(\mb{M}^\top\mb{M})}]^{1/2}$, where $\varphi_{\max}$ is the largest eigenvalue, also $\varphi_k$ and $\varphi_{\min}$ are the $k$-th and the smallest eigenvalues, respectively; 
the Frobenius norm 
$\| \mb{M} \|_F=(\sum_{i}\sum_{j}M_{ij}^2)^{1/2}$;
the matrix $\ell_1$ norm $\|\mb{M}\|_1=\max_{j}\sum_{i}|M_{ij}|$;
the entrywise $\ell_1$ norm $|\mb{M}|_1=\sum_{i,j}|M_{ij}|$ and its off-diagonal version
$|\mb{M}|_{1,\text{off}}=\sum_{i\ne j}|M_{ij}|$;
and the entrywise $\ell_{\infty}$ norm 
$|\mb{M}|_{\infty}=\max_{i,j}|M_{ij}|$.      

Denote $\textup{vec}(\mb{M})$
$=(\bd{M}_1^\top,\dots,\bd{M}_n^\top)^\top$, where $\bd{M}_j$ is the $j$-th column of $\mb{M}$. Write $\mb{M}\succ 0$ when $\mb{M}$ is positive definite.
Denote the trace and the determinant of a square matrix $\mb{M}$ by $\text{tr}(\mb{M})$ and $\text{det}(\mb{M})$, respectively.
Denote the Kronecker product by $\otimes$.
Write $x_n\asymp y_n$ if $x_n=O(y_n)$ and $y_n=O(x_n)$,
 and denote $x_n\sim y_n$ if $x_n/y_n\to 1$ as $n\to \infty$.
Define $\lceil{x}\rceil$ and $\lfloor{x}\rfloor$ to be the smallest integer $\ge x$ and the largest integer $\le x$, respectively. Let $\mathds{1}(A)$ be the indicator function of event $A$, $x_+=x\mathds{1}(x\ge 0)$ and sign$(x)=\mathds{1}(x\ge 0)-\mathds{1}(x\le 0)$. Let $A\coloneqq B$ denote that $A$ is defined to be $B$.
Denote $X\stackrel{d}{=}Y$ if $X$ and $Y$ have the same distribution.
Denote $\bd{1}_n=(1,1,\dots,1)^\top$ with length $n$ and $\mb{I}_{n\times n}$ to be the $n{\times} n$ identity matrix.
If without further notification, a constant is independent of $n$ and $p$.
Throughout the rest of the article,
we assume $p\to \infty$ as $n\to \infty$ and only use $n\to \infty$ in the
asymptotic arguments.

\subsection{A new temporal dependence measure}

Let $\mb{X}_{p\times n}\coloneqq(\bd{X}_1,...,\bd{X}_n)$, where each column  $\bd{X}_i$ follows a distribution with the same covariance matrix $\mb{\Sigma} = ({\sigma}_{k\ell})_{p\times p}$
and correlation matrix
$\mb{R}=(\rho_{k\ell})_{p\times p}$. 
Let $\bd{X}_{[1]},\dots,\bd{X}_{[p]}$ be the $p$ row vectors of $\mb{X}_{p\times n}$, and $\mb{R}_{[k]}=(\rho_{[k]}^{ij})_{n\times n}$ be the correlation matrix of $\bd{X}_{[k]}$, i.e., the autocorrelation matrix of the $k$-th time series.
For all $k$, we have the following inequalities:
\begin{equation} \label{tempo-ineq}
1\le \frac{1}{n}\| \mb{R}_{[k]}  \|_F^2\le  \| \mb{R}_{[k]}    \|_2   \le\| \mb{R}_{[k]}    \|_1\le n,
\end{equation}
where the second inequality follows from
\begin{align*}
\frac{1}{n}\text{tr}(\mb{R}_{[k]}^2)=\frac{1}{n}\sum_{i=1}^n \varphi_i^2(\mb{R}_{[k]})
\le \frac{1}{n}\varphi_{\max}(\mb{R}_{[k]})\sum_{i=1}^n \varphi_i(\mb{R}_{[k]})
=\frac{1}{n} \| \mb{R}_{[k]} \|_2 \text{tr}(\mb{R}_{[k]}),
\end{align*}
and the third inequality is obtained from Corollary 2.3.2 in \cite{Golu96}.

Define a joint temporal dependence measure that consists of two components bounded by $g_F$ and $g_2$, respectively:
\be\label{g_F,g_2}
\max_{1\le k\le p}\frac{1}{n}\| \mb{R}_{[k]}  \|_F^2\le g_F,
\  
\max_{1\le k\le p} \| \mb{R}_{[k]}    \|_2\le g_2,
\ee
and from (\ref{tempo-ineq}) we set $1\le g_F\le g_2\le n$. The proposed measure naturally describes the overall strength of temporal dependence of the $p$ time series. Particularly,
we can set  $g_2=1$ if all the $p$ time series are white noise processes, and $g_F=n$ if every pair of data points in one of the time series
are perfectly correlated or anti-correlated.

Later we will show that the convergence rates of considered estimators are nicely characterized by the bounds $g_F$ and $g_2$, which is particularly useful in obtaining convergence results for long-memory data.
Note that we do not consider cross-correlations between the multiple time series, neither assume any specific temporal decay model or stationarity for autocorrelations within each individual time series. Here are two special examples of practical interests. 

\begin{case}[High-dimensional short-memory dependence]\label{case1}
Recall that a univariate stationary time series is said to be short-memory if 
its autocorrelation function satisfies $\sum_{t=0}^\infty|\rho(t)|<\infty$.
We extend the ``short-memory" concept to multivariate time series that are allowed to be non-stationary by the property
\be\label{short mem multi}
\max_{1\le k\le p} \| \mb{R}_{[k]}    \|_1 <\infty \   \text{as}  \ n\to\infty.
\ee
Thus from (\ref{tempo-ineq}) we can set the measure bound $g_2< \infty$ as $n\to \infty$.
\end{case}

\begin{case}[Polynomial-dominated decay (PDD) model]\label{case2}
We say $\mb{X}_{p\times n}$ has PDD temporal dependence if
\be\label{PDD}
\max_{1\le k\le p}\big|\rho_{[k]}^{ij}\big|\le C_0|i-j|^{-\alpha} \quad \text{for all $i\ne j$} 
\ee
with some positive constants $C_0$ and $\alpha$.
We can then set the measure bounds
\be\label{g for PDD}
g_F=2C_0^2 H_{\lceil{n/2}\rceil}^{(2\alpha)}+1 \quad \text{and}  \quad 
g_2= 2C_0 H_{\lceil{n/2}\rceil}^{(\alpha)}+1 \ge \max_{1\le k\le p} \| \mb{R}_{[k]}    \|_1 
\ee
with the generalized harmonic number (see (25), (26) in \cite{Chle09}) 
\be\label{harmonic number bound}
H_{n}^{(\alpha)}=\sum_{k=1}^n k^{-\alpha}
<1+
\begin{cases}
\frac{n^{1-\alpha}-1}{1-\alpha},&\alpha\ne 1;\\
\log n,&\alpha=1.	
\end{cases}
\ee
The model is short-memory in the sense of 
\eqref{short mem multi} when $\alpha>1$,
but allows an individual time series to be long-memory when $0<\alpha\le 1$. 
It is worth noting that
the fractional Gaussian noise \cite{Mand68,Taqq03} and the autoregressive fractionally integrated moving average process \cite{Gran80,Hosk81}
are classical examples of stationary univariate time series 
with autocorrelation function $\rho(t)\sim Ct^{-\alpha}$ as $t\to \infty$ with $C\ne 0$ and $\alpha\in(0,1)\cup(1,2)$.
\end{case}

\subsection{Comparisons to existing work} \label{subs:compare}

For banding and tapering estimators of $\mb{\Sigma}$,
\cite{Bhat14} considered 
a weak temporal dependence
$\max_{a_n\le |i-j| \le n}{|\mb{\Lambda}^{ij} |_{\infty}}=$ $ O(n^{-2}a_n),$
 where $\mb{\Lambda}^{ij}=(\varLambda_{k\ell}^{ij})_{p\times p}$
with $\varLambda_{k\ell}^{ij}$ satisfying
 $\textup{Cov}(X_{ki},X_{\ell j})=\varLambda_{k\ell}^{ij}\sigma_{k\ell}$,
$a_n\sqrt{\log p/n}=o(1)$, and $\{a_n\}_{n\ge 1}$ is a non-decreasing sequence of non-negative integers.
That $a_n\sqrt{\log p/n}=o(1)$ implies $a_n=o(\sqrt{n})$.
Thus,
$
\left|\mb{\Lambda}^{ij}{:}|i-j|{=}\lfloor\sqrt{n}\rfloor    \right|_{\infty}\le \max_{a_n\le |i-j| \le n}|\mb{\Lambda}^{ij} |_{\infty}
$
$=O(n^{-2}a_n)=o(n^{-3/2}).
$
Then
$
\sum_{|i-j|=0}^{\infty}|\rho_{[k]}^{ij}|<\infty$
for any given $k$, which means that
the time series cannot be long-memory,
not even $\rho_{[k]}^{ij}\asymp |i-j|^{-\alpha}$ with $\alpha\in(0,3]$.
\cite{Bhat14_2} extended the banding and tapering techniques to the estimation of $\text{Cov}(\bd{X}_{j},\bd{X}_{j+k})$, $k\ge 0$, for the stationary infinite-order moving average model. It is easy to show that their time series cannot be long-memory.

\cite{Chen13} considered the hard thresholding estimation of $\mb{\Sigma}$  and an $\ell_1$-MLE type estimation of $\mb{\Omega}$
using the functional dependence measure of \cite{Wu05}. 
Without loss of generality,	
assume that the first row of $\mb{X}_{p\times n}$, $\{X_{1t}\}$, is a stationary process  with autocovariance function
$\gamma_1(t)$. Following their setup by letting $E(X_{1t})=0$, 
we have $\gamma_1(t)=E(X_{11}X_{1,t+1})$. By the argument in the proof of Theorem 1 in \cite{Wu09} together with Lyapunov's inequality \cite{Bill95} and Theorem~1 of \cite{Wu05}, one can see that their model requires $\sum_{t=0}^{\infty}|\gamma_1(t)|<\infty$, which indicates that
$\{X_{1t}\}$ cannot be long-memory.

\cite{Zhou14} considered estimating a separable covariance $\textup{Cov}(\bd{X}_{pn})=\mb{A}\otimes\mb{B}$,
where $\bd{X}_{pn} \coloneqq \text{vec}(\mb{X}_{p\times n})$.
Her model implies the same autocorrelation coefficients $\{\rho_{[k]}^{ij}\}_{1\le i,j\le n}$ for all $k$, indicating a rather restrictive model with homogeneous decay rate for all $p$ time series.

\begin{figure}[t!]
\centering
\includegraphics[width=4.5in]{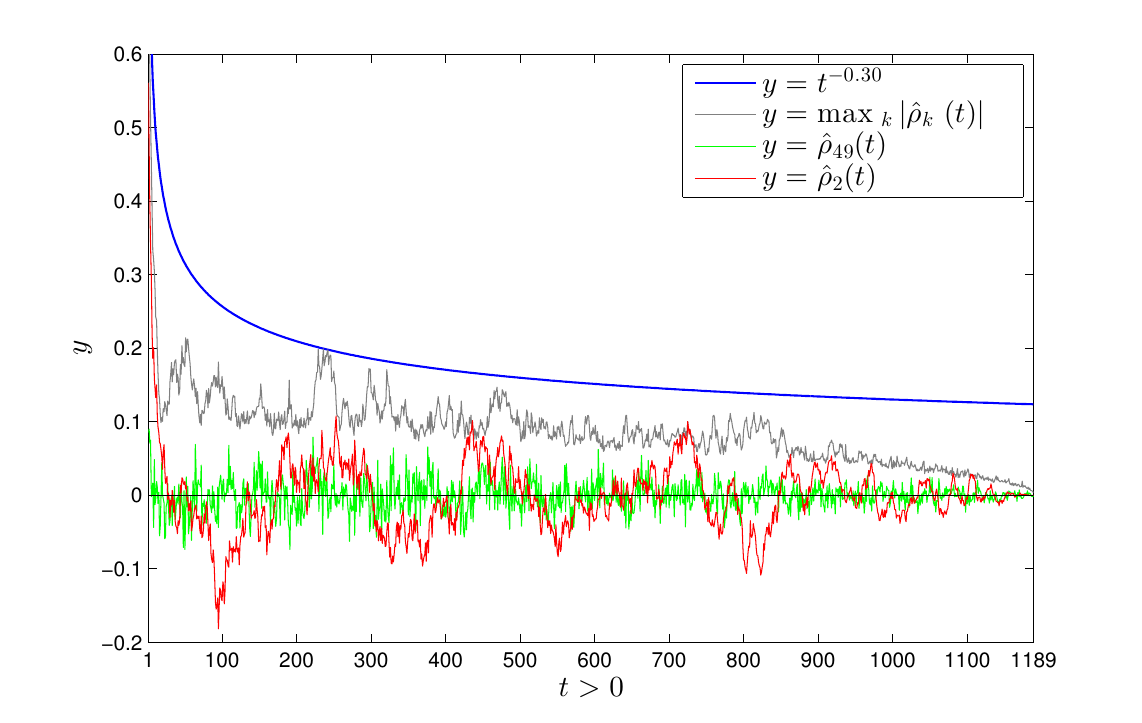}
 \caption{Sample autocorrelations of brain nodes.}
\label{autocorrelation}
\end{figure}

Now take a look at the rfMRI data example of a single subject which will be further analyzed in Subsection~\ref{subs:real data}.
The data set consists of 1190 temporal brain images. We consider 907 functional brain nodes in each image.
All node time series have passed the Priestley-Subba Rao test for stationarity \cite{Prie69}, 
the generalized Jarque-Bera test for Gaussianity \cite{Bai05}, and the Hinich's bispectral test for linearity \cite{Hini82}.
Hence the linear spatio-temporal model that we will define in the next subsection with sub-Gaussian tails seems adequate for the data.	
There are 134 time series detected as long-memory by the GPH test \cite{Gewe83}.  All these tests are conducted
with a significant p-value of 0.05 adjusted by the false discovery rate controlling procedure of \cite{Benj01}.  Hence, the weak temporal dependence models of \cite{Bhat14,Bhat14_2,Chen13} do not apply to these long-memory time series. For node $k$, its
autocorrelation function $\rho_k(t):=\rho_{[k]}^{ij}$, $t=|i-j|$, can be approximated by the sample autocorrelation function $\hat{\rho}_k(t)$.  
Figure~\ref{autocorrelation} shows that the rfMRI data approximately satisfy the PDD model~\eqref{PDD} with
$C_0=1$ and $\alpha=0.30$ since $\max_{1\le k\le p}|\hat{\rho}_k(t)|\le t^{-0.30}$.
 The figure also
illustrates the estimated autocorrelation functions for two randomly selected brain nodes, which  clearly have different patterns, indicating that the assumption of homogeneous decay rates for all time series in \cite{Zhou14} does not hold.

\subsection{Data generating mechanism}
Throughout the article, we assume that the vectorized data are obtained from
the linear spatio-temporal model
\be
\label{linear form}
\bd{X}_{pn}  
	\coloneqq\text{vec}(\mb{X}_{p\times n})
=\mb{H}\bd{e}+\bd{\mu}_{pn},
\ee
where $\mb{H}=(h_{ij})_{pn\times m}$ is a real deterministic matrix,
$\bd{\mu}_{pn}=\bd{1}_n\otimes \bd{\mu}_p$,
and the
random vector $\bd{e}=(e_1,\dots,e_{m})^\top$ consists of $m$ independent 
(not necessarily i.i.d.)	
random variables satisfying $E(e_i)=0$ and $E(e_i^2)=1$ for all $i$.		
We allow $m=\infty$ by requiring that for each $i$, $\sum_{j=1}^{m}h_{ij}e_j$ converges both almost surely and in mean square when $m\to\infty$. A sufficient and necessary condition for both modes of convergence is $\sum_{j=1}^\infty h_{ij}^2<\infty$ for every $i$ (see Theorem~8.3.4 and its proof in \cite{Athr06}).
Under these two modes of convergence, it can be shown that  $E(\mb{H}\bd{e})=\mb{H}E(\bd{e})$ and $\textup{Cov}(\mb{H}\bd{e})=\mb{H}\textup{Cov}(\bd{e})\mb{H}^\top$ (see Proposition~2.7.1 in \cite{Broc91}).
Hence, for either finite or infinite $m$, we have
$E(\bd{X}_{pn})=\bd{\mu}_{pn}$ and $\textup{Cov}(\bd{X}_{pn})=\mb{HH}^\top$ with all $n$ submatrices of dimension $p\times p$  on the diagonal equal to $\mb{\Sigma}$ and temporal correlations determined by the off-diagonal submatrices.
In filtering theory,
matrix $\mb{H}$ is said to be a linear spatio-temporal coloring filter \cite{Fomi99,Mano05}, which generates the output $\bd{X}_{pn}$ by introducing both spatial and temporal dependence in the input independent variables $e_1,\dots,e_{m}$. We will use $\mb{X}_{p\times n}$ and $\bd{X}_{pn}$ exchangeably.      

The following are two examples of \eqref{linear form} which are often seen in the literature.
In particular, two processes used in analyzing fMRI data, i.e., the multivariate fractional Gaussian noise \cite{Ciuc14} and the vector autoregressive model \cite{Harr03}, are special cases of these two examples, respectively.

\begin{example}[Gaussian data]
\label{eg1}
Assume that $\bd{X}_{pn}$ has a multivariate Gaussian distribution $\mathcal{N}(\bd{\mu}_{pn},\mb{\Delta})$.
Then $\mb{\Delta}=\mb{HH}$ with a symmetric real matrix $\mb{H}$.
If $\mb{\Delta}\succ 0$,
then $\bd{X}_{pn}=\mb{H}\bd{e}+\bd{\mu}_{pn}$ with
$\bd{e}=\mb{H}^{-1}(\bd{X}_{pn}-\bd{\mu}_{pn})\sim \mathcal{N}(\bd{0},\mb{I}_{pn\times pn})$. If $\mb{\Delta}$ is singular, then $\bd{X}_{pn}$ has a degenerate multivariate Gaussian distribution, and can be expressed as $\bd{X}_{pn}\stackrel{d}{=}\mb{H}\bd{e}+\bd{\mu}_{pn}$ with any $\bd{e}\sim  \mathcal{N}(\bd{0},\mb{I}_{pn\times pn})$.
In fact, replacing ``$=$" in \eqref{linear form} by ``$\stackrel{d}{=}$" does not affect the theoretical results.
\end{example}

\begin{example}[Moving average processes]
\label{eg2}
Consider the processes
\be
\label{infty sum}
\bd{X}_j=\sum_{\ell=0}^L\mb{B}_\ell\bd{e}_{j-\ell}, \quad \text{with}\quad 0\le L\le \infty,
\ee 
where the case with $L=\infty$ is well-defined in the sense of entrywise almost-sure convergence and mean-square convergence, $\{\mb{B}_\ell\}$ are $p\times p$ real deterministic matrices,
and $\bd{e}_j=(e_{1j},e_{2j}\dots,e_{pj})^\top$ 
is a vector with  independent zero-mean and unit-variance entries $\{e_{ij}, 1\le i \le p,\, -\infty \le j \le n\}$.
Since every $X_{ij}$ is a linear combination of $\{e_{st}\}$, we always can find a matrix $\mb{H}$ such that $\bd{X}_{pn}=\mb{H}\bd{e}$ with $\bd{e}=(\bd{e}_{1-L}^\top,\bd{e}_{2-L}^\top,\dots,\bd{e}_{n}^\top)^\top$.
It is well-known that any causal vector
autoregressive moving average process of the form
$
\bd{X}_j-\mb{A}_1\bd{X}_{j-1}-\dots-\mb{A}_a\bd{X}_{j-a}=\bd{e}_j+\mb{M}_1\bd{e}_{j-1}
+\dots+\mb{M}_b\bd{e}_{j-b}
$
with finite nonnegative integers $a$ and $b$, and real deterministic matrices
$\{\mb{A}_i,\mb{M}_k\}$,
can be written in the form of \eqref{infty sum} with $L=\infty$ \cite[ pp. 418]{Broc91}.
Model~\eqref{infty sum} with $L=\infty$ is widely studied in recent literature of high dimensional time series (see, e.g., \cite{Chen13,Chen16,Bhat14_2,Liu15}).
\end{example}

We will consider the following three types of moment conditions for the basis random variables $e_1,\dots,e_{m}$ in (\ref{linear form}), corresponding to sub-Gaussian, generalized sub-exponential, and polynomial-type tails, respectively.
Let $Z$ be a random variable, and $K$, $\vartheta$, and $\eta_k$ be positive constants.
\begin{enumerate}
\renewcommand{\theenumi}{\textup{(C\arabic{enumi})}}
	\item\label{C1} {\it Sub-Gaussian tails}: For all $k\ge 1$, we have 
	$
	(E|Z|^k)^{1/k}\le K k^{1/2}.
	$
	
	\item\label{C2} {\it Generalized sub-exponential tails}: 
	For some $\vartheta\in (0,2)$ and all $k\ge \vartheta$, we have 
	$
	 (E|Z|^k)^{1/k}\le K (k/\vartheta)^{1/\vartheta}.
	$
	 
	\item\label{C3}  {\it Polynomial-type tails}: For some $k\ge 4$, we have 
	$
	(E|Z|^k)^{1/k}\le \eta_k.
	$
\end{enumerate}

We do not consider condition~\ref{C1} as a special case of condition~\ref{C2} by setting $\vartheta=2$ due to the fact that sharper convergence rates can be obtained under \ref{C1}. We can apply  the Hanson-Wright inequality to \ref{C1} \cite{Rude13}, but not able to extend it to  \ref{C2} because the moment generating function of $Z^2$ is no longer finite in an open interval around zero 
(see Proposition~7.23 and inequality~(7.32) in \cite{Fouc13}, and Lemma~5.5 in \cite{Vers12}).
Conditions \ref{C1} and \ref{C2} can be equivalently written as
$
P(|Z|\ge u)\le 2\exp(-u^{\vartheta}/C)
$
with some constant $C>0$ for all $u\ge 0$, where for the former $\vartheta=2$. 
Condition~\ref{C3} implies $P(|Z|\ge u)\le \eta_k^k/u^k$ for all $u>0$.
Conversely, if $P(|Z|\ge u)=O(u^{-k'})$ with some $k'\in (0,k)$ as $u\to \infty$, then ${(E|Z|^k)^{1/k}<\infty}$.

\section{Estimation of covariance and correlation matrices for sub-Gaussian data}
\label{s:cov est}
Consider the $\ell_q$-ball sparse covariance matrices \cite{Bick08a, Roth09}
\be
\label{sparsity Sigma}
\mathcal{U}(q,c_p,v_0)
=\bigg\{
\mb{\Sigma}:   \max_{1\le i\le p} \sum_{j=1}^p|\sigma_{ij}|^q\le c_p,\max_{1\le i\le p}\sigma_{ii}\le v_0
\bigg\},
\ee
and the corresponding correlation matrices
\be
\label{sparsity R}
\mathcal{R}(q,c_p)
=\bigg\{
\mb{R}:  \max_{1\le i\le p}\sum_{j=1}^p|\rho_{ij}|^q\le c_p
\bigg\},
\ee
where constants $v_0>0$ and $0\le q<1$.
For any thresholding parameter $\tau\ge 0$, define a generalized thresholding function \cite{Roth09}
by $s_{\tau}: \mathbb{R}\to\mathbb{R}$ satisfying the following conditions for all $z\in \mathbb{R}$:
\begin{inparaenum}[(i)]
\item $|s_{\tau}(z)|\le |z|$; \label{(i)}
\item $s_{\tau}(z)=0$ for $|z|\le {\tau}$; \label{(ii)}
\item $|s_{\tau}(z)-z|\le {\tau}$. \label{(iii)}
\end{inparaenum}
Such defined generalized thresholding function  covers many widely used thresholding functions, including
hard thresholding
$
s_{\tau}^H(z)=z\mathds{1}(|z|> {\tau}),
$
soft thresholding
$
s_{\tau}^S(z)=\text{sign}(z)(|z|-{\tau})_+,
$
smoothly clipped absolute deviation
and adaptive lasso
thresholdings. See details about these examples in
\cite{Roth09}.
We define the generalized thresholding estimators of $\mb{\Sigma}$ and $\mb{R}$, respectively, by
$$
S_{\tau_1}(\widehat{\mb{\Sigma}})=(s_{\tau_1}(\hat{\sigma}_{ij}))_{p\times p}
\ \
\mbox{and}
\ \
S_{\tau_2}(\widehat{\mb{R}})=\left(s_{\tau_2}(\hat{\rho}_{ij})\mathds{1}(i\ne j)+\mathds{1}(i=j)\right)_{p\times p},
$$
where $\widehat{\mb{\Sigma}}\coloneqq(\hat{\sigma}_{ij})_{p\times p}$ is the sample covariance matrix given by
\be
\label{sample cov}
\widehat{\mb{\Sigma}}=\frac{1}{n}\sum_{i=1}^n\bd{X}_i\bd{X}_i^\top-\bar{\bd{X}}\bar{\bd{X}}^\top
\ee
with $\bar{\bd{X}}=\frac{1}{n}\sum_{i=1}^n\bd{X}_i$,
and $
\widehat{\mb{R}}\coloneqq\left(\hat{\rho}_{ij}\right)_{p\times p}=\left(\hat{\sigma}_{ij}/\sqrt{\hat{\sigma}_{ii}\hat{\sigma}_{jj}} \right)_{p\times p}
$ is the sample correlation matrix.
Define 	
\be\label{u for subgauss}
u_1=\max\left\{
(\log p)g_2/n, \ [(\log p) g_F /n]^{1/2}
\right\}.
\ee
Then we have the following results.

\begin{thm}
	\label{thm1}
Suppose that $\mb{X}_{p\times n}$ is generated from \eqref{linear form} with all $e_i$ satisfying condition~\ref{C1} with the same $K$. Uniformly on $\mb{\Sigma}\in \mathcal{U}(q,c_p,v_0)$ and
$\{\mb{R}_{[k]}\}_{k=1}^p$ subject to \eqref{g_F,g_2},
for sufficiently large constant $M_1>0$ depending only on $v_0$ and $K$, if
	${\tau_1}=M_1u_1$ and
	$u_1=o(1)$,
	then
	\begin{align}
	|S_{\tau_1}(\widehat{\mb{\Sigma}})-\mb{\Sigma}|_{\infty}&=O_P(u_1),
	\nonumber\\
	\|S_{\tau_1}(\widehat{\mb{\Sigma}})-\mb{\Sigma}\|_2&=O_P\left(c_pu_1^{1-q}\right),\label{prob 2-norm}
	\\
	\frac{1}{p}\|S_{\tau_1}(\widehat{\mb{\Sigma}})-\mb{\Sigma}\|_F^2&=O_P\left(c_pu_1^{2-q}\right).\label{prob Frobenius norm}
	\end{align}
	Moreover, if $p\ge n^{c}$ for some constant $c>0$, then with sufficiently large $M_1$ additionally depending on $c$ and $q$, we have
	\begin{align*}
	E\left(|S_{\tau_1}(\widehat{\mb{\Sigma}})-\mb{\Sigma}|_{\infty}^2\right)&=O(u_1^2),
	\nonumber\\
	E\left(\|S_{\tau_1}(\widehat{\mb{\Sigma}})-\mb{\Sigma}\|_2^2\right)&=O\left(c_p^2u_1^{2-2q}\right),
	\\
	\frac{1}{p}E\left(\|S_{\tau_1}(\widehat{\mb{\Sigma}})-\mb{\Sigma}\|_F^2\right)&=O\left(c_pu_1^{2-q}\right).
	\end{align*}
\end{thm}

\begin{remark}\label{rmk1}
When $(\log p)/n = o(1)$, if $u_1=O(\sqrt{(\log p)/n})$, which is true when $g_2<\infty$ that holds for short-memory multivariate time series 
satisfying \eqref{short mem multi}, then the in-probability convergence rates given in \eqref{prob 2-norm} and \eqref{prob Frobenius norm} are the same rates as those for i.i.d. observations given in \cite{Bick08a} and \cite{Roth09}. 
The same in-probability convergence rates are also obtained by \cite[Proposition 5.1]{Basu15} for certain short-memory stationary Gaussian data using the hard thresholding method.   
\end{remark}

\begin{remark}\label{rmk2}
For the PDD  temporal dependence given in
\eqref{PDD}, by \eqref{g for PDD} and \eqref{harmonic number bound}, together with $u_1=o(1)$, 
we have
\[
u_1\lesssim
\begin{cases}
\sqrt{(\log p)/{n}}, & \alpha>1;\\
\max \big\{[(\log p)(\log n)] /{n},  \sqrt{(\log p)/{n}}  \big\},  & \alpha=1;\\
\max\big\{ (\log p)/{n^\alpha},\sqrt{(\log p)/{n}} 
\big\}, & \alpha\in(1/2,1);\\ 
\max\big\{(\log p)/{n^{1/2}},\sqrt{[(\log p)(\log n)]/{n}} 
\big\}, & \alpha=1/2;\\
(\log p)/{n^\alpha}, & \alpha\in(0,1/2).
\end{cases}
\]
Here  $x_n \lesssim y_n$ denotes $x_n=O(y_n)$. Note that the case with $\alpha>1$ is short-memory in the sense of \eqref{short mem multi}. 
When 
$\alpha=1$ and $(\log n)\sqrt{(\log p)/n}=O(1)$, or when
$\alpha\in (1/2,1)$ and $(\log p)^{1/2}n^{1/2-\alpha}=O(1)$, which allows some individual time series to be long-memory, we also have $u_1=O(\sqrt{(\log p)/n})$, yielding the same 
convergence rates as in the  case with i.i.d. data.
\end{remark}

\begin{thm}[Sparsistency and sign-consistency]
\label{thm2}
Under the conditions for the convergence in probability given in Theorem~\ref{thm1},
we have
${s_{\tau_1}(\hat\sigma_{ij})= 0}$ for all $(i,j)$ where $\sigma_{ij}=0$ with probability tending to 1. If further assume all nonzero entries of $\mb{\Sigma}$ satisfy $|\sigma_{ij}|\ge 2\tau_1$, then we have
$\textup{sign}(s_{\tau_1}(\hat{\sigma}_{ij}))$
$=\textup{sign}(\sigma_{ij})$ for all $(i,j)$ where $\sigma_{ij}\ne 0$ with probability tending to 1.
\end{thm}

\begin{cor}\label{cor1}
Theorems~\ref{thm1} and \ref{thm2} hold with $\widehat{\mb{\Sigma}}, \mb{\Sigma},\hat{\sigma}_{ij},\sigma_{ij}$, $\mathcal{U}(q,c_p,v_0)$, $\tau_1$, and $M_1$  replaced by $\widehat{\mb{R}}, \mb{R},\hat{\rho}_{ij},\rho_{ij}$, $\mathcal{R}(q,c_p)$, $\tau_2$, and $M_2$, respectively, 
where $M_2$ does not depend on $v_0$.
\end{cor}

We now provide the minimax optimal rates for estimating covariance and correlation matrices over certain sets of distributions of $\mb{X}_{p\times n}$, including the short-memory case \eqref{short mem multi} and some long-memory cases \eqref{PDD} with $\alpha\in (1/2,1]$.

\begin{thm}[Minimax rates]
\label{minimax}
Let $K_G=\sup_{k\ge 1} \sqrt{2/k}\left[ \Gamma(\frac{1+k}{2})/\sqrt{\pi} \right]^{1/k}$, where $\Gamma(x)$ is the gamma function.
Let $\mathcal{P}_1(q,c_p,v_0,g_F,g_2,K,\kappa)$ be the set of distributions of $\mb{X}_{p\times n}$  generated from \eqref{linear form} with all $e_i$ satisfying \ref{C1} with constant $K\ge K_G$,
$\mb{\Sigma}\in \mathcal{U}(q,c_p,v_0)$, and $\{\mb{R}_{[k]}\}_{k=1}^p$ subject to \eqref{g_F,g_2},
where the constant ${\kappa\ge 1}$ is used in setting $u_1\le \kappa\sqrt{(\log p)/n}$. Let $\mathcal{P}_2(q,c_p,g_F,g_2,K,\kappa)$ be the corresponding set to $\mathcal{P}_1$ with $\mb{\Sigma}\in \mathcal{U}(q,c_p,v_0)$ replaced by
$\mb{R}\in \mathcal{R}(q,c_p)$.
Let $\mathfrak{D}$ denote the distribution of $\mb{X}_{p\times n}$.
If $\sqrt{(\log p)/n}=o(1)$, 
\be\label{Cai cond}
p\ge n^{c_1} \quad \text{and}\quad c_p\le c_2n^{(1-q)/2}(\log p)^{-(3-q)/2}
\ee
with some constants $c_1>1$ and $c_2>0$, then for any estimator $\widetilde{\mb{\Sigma}}$ we have 
\[
\inf_{\widetilde{\mb{\Sigma}}}\sup_{\mathfrak{D}\in\mathcal{P}_1}E_{\mb{X}_{p\times n}|\mathfrak{D}}
\left(\| \widetilde{\mb{\Sigma}}-\mb{\Sigma}  \|_2^2 \right)
\asymp  c_p^2\left(\frac{\log p}{n}\right)^{1-q},
\]
\[
\inf_{\widetilde{\mb{\Sigma}}}\sup_{\mathfrak{D}\in\mathcal{P}_1}\frac{1}{p}E_{\mb{X}_{p\times n}|\mathfrak{D}}
\left(\| \widetilde{\mb{\Sigma}}-\mb{\Sigma}  \|_F^2 \right)
\asymp  c_p\left(\frac{\log p}{n}\right)^{1-q/2}.
\]
Additionally if $c_p>1$, then for any estimator $\widetilde{\mb{R}}$ we have
\[
\inf_{\widetilde{\mb{R}}}\sup_{\mathfrak{D}\in\mathcal{P}_2}E_{\mb{X}_{p\times n}|\mathfrak{D}}
\left(\| \widetilde{\mb{R}}-\mb{R}  \|_2^2 \right)
\asymp  c_p^2\left(\frac{\log p}{n}\right)^{1-q},
\]
\[
\inf_{\widetilde{\mb{R}}}\sup_{\mathfrak{D}\in\mathcal{P}_2}\frac{1}{p}E_{\mb{X}_{p\times n}|\mathfrak{D}}
\left(\| \widetilde{\mb{R}}-\mb{R}  \|_F^2 \right)
\asymp  c_p\left(\frac{\log p}{n}\right)^{1-q/2}.
\]
For sufficiently large positive constants $M_1$ and $M_2$, with $\tau_1=M_1u_1$ and $\tau_2=M_2u_1$, the generalized thresholding estimators $S_{\tau_1}(\widehat{\mb{\Sigma}})$ and $S_{\tau_2}(\widehat{\mb{R}})$
attain the above minimax optimal rates, respectively.
\end{thm}

The assumption~\eqref{Cai cond} follows \cite{CZ12} who studied the minimax optimal rates of the covariance matrix estimation for i.i.d. data. From Remarks~\ref{rmk1} and \ref{rmk2} we see that suitable $(g_F,g_2,\kappa)$ in $\mathcal{P}_1$ and $\mathcal{P}_2$ allow $\mb{X}_{p\times n}$ to have short memory \eqref{short mem multi}, or to 
follow the PDD model
\eqref{PDD} 
with $\alpha\in(1/2,1]$ (thus with some time series to be long-memory) under some additional conditions for $n$ and $p$ discussed in Remark~\ref{rmk2}.

\section{Estimation of precision matrix for sub-Gaussian data}\label{s:prec est}
We consider both the CLIME and the SPICE methods for the estimation of $\mb{\Omega}$, which were originally developed for i.i.d. observations.

\subsection{CLIME estimation}
Following \cite{CLL11},
we consider the following set of precision matrices
\[
\mathcal{G}_1(q,c_p,M_p,v_0)
{=}\Big\{
\bd\Omega\succ 0:
\max_{1\le i\le p}\sum_{j=1}^p|\omega_{ij}|^q\le c_p,
  \| \bd\Omega \|_1\le M_p,
\max_{1\le i\le p}\sigma_{ii} \le v_0
\Big\},
\]
where constant $0\le q<1$, and $(c_p,M_p)$ are allowed to depend on $p$. 
Though the condition 	
	$\max_i\sigma_{ii} \le v_0$
is not explicitly provided by \cite{CLL11} in their original set definition, it is implied by
their moment conditions (see their (C1) and (C2)).
Note that the above ${\cal G}_1$ contains $\ell_q$-ball sparse matrices such as those with exponentially decaying entries from the diagonal, for example, AR(1) matrices. For an invertible band matrix $\mb{\Sigma}$, its inverse matrix $\mb{\Omega}$ generally has exponentially decaying entries from the diagonal \cite{Demk84}.

Let  $\widehat{\mb{\Theta}}_{\varepsilon,\lambda_1} \coloneqq(\hat{\theta}_{ij}^{(\varepsilon,\lambda_1)})_{p\times p}$
be a solution of the following optimization:
\be
\label{optim1}
\min | \mb{\Theta}|_1 \quad \mbox{subject to}\quad
|\widetilde{\mb{\Sigma}}_{\varepsilon}\mb{\Theta}-\mb{I}_{p\times p}|_{\infty}\le \lambda_1,\quad \mb{\Theta}\in \mathbb{R}^{p\times p},
\ee   
where
$\widetilde{\bd\Sigma}_{\varepsilon}=\widehat{\bd\Sigma}+\varepsilon \mb{I}_{p\times p}$, $\widehat{\bd\Sigma}$ is given in \eqref{sample cov}, $\varepsilon\ge 0$ is a perturbation parameter introduced for the same reasons given in \cite{CLL11} and can
be set to be $n^{-1/2}$ in practice (see Remark~\ref{perturb} below),
and $\lambda_1$ is a tuning parameter.
The CLIME estimator 
$\widehat{\mb{\Omega}}_{\varepsilon,\lambda_1}\coloneqq(\hat{\omega}_{ij}^{(\varepsilon,\lambda_1)})_{p\times p}$ 
is then obtained by symmetrizing
$\widehat{\mb{\Theta}}_{\varepsilon,\lambda_1}$
with
\[
\hat{\omega}_{ij}^{(\varepsilon,\lambda_1)}=\hat{\omega}_{ji}^{(\varepsilon,\lambda_1)}
=\hat{\theta}_{ij}^{(\varepsilon,\lambda_1)}\mathds{1}(|\hat{\theta}_{ij}^{(\varepsilon,\lambda_1)}|\le |\hat{\theta}_{ji}^{(\varepsilon,\lambda_1)}|)
+\hat{\theta}_{ji}^{(\varepsilon,\lambda_1)}\mathds{1}(|\hat{\theta}_{ij}^{(\varepsilon,\lambda_1)}|> |\hat{\theta}_{ji}^{(\varepsilon,\lambda_1)}|).
\]

For $1\le i\le p$, let 
$\hat{\bd{\beta}}_{i}^{(\varepsilon,\lambda_1)}$
be a solution of the following convex optimization problem:
\be
\label{optim2}
\min | \bd{\beta}_{i} |_1 \quad \mbox{subject to}\quad
|\widetilde{\mb{\Sigma}}_{\varepsilon}\bd{\beta}_{i}-\bd{e}_i|_{\infty}\le \lambda_1,
\ee
where $\bd{\beta}_{i}$
is a real vector and $\bd{e}_i$ is the vector with 1 in the $i$-th coordinate and 0 in all other coordinates.
\cite{CLL11} showed that solving the optimization problem \eqref{optim1} is equivalent to solving the $p$ optimization problems given in
\eqref{optim2}, i.e., 
$\widehat{\mb{\Theta}}_{\varepsilon,\lambda_1}=(\hat{\bd{\beta}}_{1}^{(\varepsilon,\lambda_1)},\dots,\hat{\bd{\beta}}_{p}^{(\varepsilon,\lambda_1)})$.
This equivalence is useful for both numerical implementation and theoretical analysis.
The following theorem gives the convergence results of CLIME under temporal dependence.

\begin{thm}
\label{clime_thm1}
Suppose that $\mb{X}_{p\times n}$ is generated from \eqref{linear form} with all $e_i$ satisfying condition~\ref{C1} with the same $K$. Uniformly on $\mb{\Omega}\in \mathcal{G}_1(q,c_p,M_p,v_0)$ and
$\{\mb{R}_{[k]}\}_{k=1}^p$ subject to \eqref{g_F,g_2}, 
for sufficiently large constant $M>0$ depending only on $v_0$ and $K$, if
	${\lambda_1}=MM_pu_1$, $0\le \varepsilon\le u_1$ and $u_1=o(1)$
with $u_1$ defined in \eqref{u for subgauss},
	then
\begin{align*}
|\widehat{\mb{\Omega}}_{\varepsilon,\lambda_1}-\bd\Omega|_{\infty}&= O_P(M_p^2u_1),\\
\|\widehat{\mb{\Omega}}_{\varepsilon,\lambda_1}-\bd\Omega\|_2&=O_P\left(c_p(M_p^2u_1)^{1-q}\right),
\\
\frac{1}{p}\|\widehat{\mb{\Omega}}_{\varepsilon,\lambda_1}-\bd\Omega\|_{F}^2&=O_P\left(c_p(M_p^2u_1)^{2-q}\right).
\end{align*}
Moreover, if $p\ge n^c$ 
and $\min\{p^{-C},u_1\}\le \varepsilon\le u_1$ 
for some positive constants $c$ and $C$, then with sufficiently large $M$ additionally depending on 
$c,C$ and $q$, we have
\begin{align}
E\left(|\widehat{\mb{\Omega}}_{\varepsilon,\lambda_1}-\bd\Omega  |_{\infty}^2\right)&=O\left((M_p^2u_1)^2\right),
\nonumber\\
E\left(\|\widehat{\mb{\Omega}}_{\varepsilon,\lambda_1}-\bd\Omega\|_2^2\right)&=O\left(c_p^2(M_p^2u_1)^{2-2q}\right),\label{clime 2 in msq}
\\
\frac{1}{p}E\left(\|\widehat{\mb{\Omega}}_{\varepsilon,\lambda_1}-\bd\Omega\|_F^2\right)&=O\left(c_p(M_p^2u_1)^{2-q}\right).\label{clime F in msq}
\end{align}
\end{thm}

\begin{remark}
If $(\log p)/n=o(1)$ and $u_1=O(\sqrt{(\log p)/n})$, then the above convergence rates are the same as those for i.i.d. data given in  \cite{CLL11}. Additionally, if $M_p$ is a constant, then the mean-square convergence rates of CLIME in \eqref{clime 2 in msq} and \eqref{clime F in msq} attain
the minimax optimal convergence rates for i.i.d. data shown in \cite{CLZ16} under slightly different assumptions.
From Remarks~\ref{rmk1} and \ref{rmk2} we see that $u_1=O(\sqrt{(\log p)/n})$ when $u_1=o(1)$ can be achieved for the short-memory case \eqref{short mem multi} and also for the long-memory cases satisfying \eqref{PDD} with $\alpha\in (1/2,1]$ under some additional conditions for $n$ and $p$. 
\end{remark}

\begin{remark}\label{perturb}
As discussed in \cite{CLL11}, the perturbation parameter $\varepsilon>0$ is used
for a proper initialization of $\{\bd{\beta}_{i}\}$
in the numerical implementation of \eqref{optim2}, and also ensures the existence of $E(\|\widehat{\mb{\Omega}}_{\varepsilon,\lambda_1}-\bd\Omega\|_{2}^2)$. Since  $g_F \ge 1$, from (\ref{u for subgauss}) we have that $u_1\ge \sqrt{(\log p)/n}\ge  n^{-1/2}$.
Hence when $p\ge n^c$,
let $C=1/(2c)$, we have $p^{-C}\le n^{-1/2}\le u_1$.
Thus,
we can simply let $\varepsilon=n^{-1/2}$ in practice, which is also the default setting of the R package {\tt flare} \cite{Li15} that
implements the CLIME algorithm. The same choice of $\varepsilon$ is given in (10) of \cite{CLL11} for i.i.d. observations.
\end{remark}

To better recover the sparsity structure of $\mb{\Omega}$, \cite{CLL11} introduced additional thresholding on $\widehat{\mb{\Omega}}_{\varepsilon,\lambda_1}$. Similarly, we may
define a hard-thresholded CLIME estimator 
$\widetilde{\mb{\Omega}}_{\varepsilon,\lambda_1,\xi}=(\tilde{\omega}_{ij}^{(\varepsilon,\lambda_1,\xi)})_{p\times p}$ 
by
$
\tilde{\omega}_{ij}^{(\varepsilon,\lambda_1,\xi)}=\hat{\omega}_{ij}^{(\varepsilon,\lambda_1)}\mathds{1}(|\hat{\omega}_{ij}^{(\varepsilon,\lambda_1)}|> \xi)
$
with a tuning parameter $\xi\ge 4M_p\lambda_1$. Although such an estimator enjoys nice theoretical properties given below, how to practically select $\xi$ remains unknown.

\begin{thm}[Sparsistency and sign-consistency]
\label{clime_thm2}
Under the conditions for the convergence in probability given in Theorem~\ref{clime_thm1}, we have
$
\tilde{\omega}_{ij}^{(\varepsilon,\lambda_1,\xi)}
=0$ for all $(i,j)$ where $\omega_{ij}=0$ with probability tending to 1. If further assume all nonzero entries of $\bd\Omega$ satisfy $|\omega_{ij}|>\xi+4M_p\lambda_1$, then we have
$
\textup{sign}(
\tilde{\omega}_{ij}^{(\varepsilon,\lambda_1,\xi)}  
)=\textup{sign}(\omega_{ij})
$ for all $(i,j)$ where $\omega_{ij}\ne 0$ with probability tending to 1.
\end{thm}

\subsection{SPICE estimation}
For i.i.d. observations, \cite{Roth08} proposed the SPICE method for estimating the following precision matrix $\mb{\Omega}$
\[
\mathcal{G}_2(s_p,v_0){=}
\Big\{
\mb{\Omega}: \sum_{1\le i\ne j\le p}\mathds{1}(\omega_{ij}\ne 0)\le s_p,
 0{<}v_0^{-1}{\le} \varphi_{\min}(\mb{\Omega}){\le}\varphi_{\max}(\mb{\Omega}){\le} v_0
\Big\},
\]
where $s_p$ determines the sparsity of $\mb{\Omega}$ and can depend on $p$, and $v_0$ is a constant.
Two types of SPICE estimators were proposed:
\be
\label{spice1}
\widetilde{\mb{\Omega}}_{\lambda_2}=\argmin_{\substack{{\mb{\Theta}}\succ 0, {\mb{\Theta}}={\mb{\Theta}}^\top}}
\left\{
\text{tr}(\mb{\Theta}\widehat{\mb{\Sigma}})-\log\det(\mb{\Theta})+\lambda_2 |\mb{\Theta}|_{1,\text{off}}
\right\},
\ee
and
\begin{align}
\label{spice2}
\widehat{\mb{\Omega}}_{\lambda_2}&
  \coloneqq (\hat{\omega}_{ij}^{(\lambda_2)})_{p\times p} 
=\widehat{\mb{W}}^{-1} \widehat{\mb{K}}_{\lambda_2}     \widehat{\mb{W}}^{-1}
\quad \text{with}\\
\widehat{\mb{K}}_{\lambda_2}&=
\argmin_{\substack{{\mb{\Theta}}\succ 0, {\mb{\Theta}}={\mb{\Theta}}^\top}}
\left\{
\text{tr}(\mb{\Theta}\widehat{\mb{R}})-\log\det(\mb{\Theta})+\lambda_2 |\mb{\Theta}|_{1,\text{off}}
\right\},\nonumber
\end{align}
where
$\lambda_2>0$ is a tuning parameter, and
$\widehat{\mb{W}} = {\rm diag}\{\sqrt{\hat\sigma_{11}}, \dots, \sqrt{\hat\sigma_{pp}}\}$. We can see that $\widehat{\mb{K}}_{\lambda_2}$ is the SPICE estimator of $\mb{K}\coloneqq\mb{R}^{-1}$.
The SPICE estimator \eqref{spice1} is a slight modification of the graphical Lasso (GLasso) estimator of \cite{Frie08}.
GLasso uses $|\mb{\Omega}|_1$ rather than $|\mb{\Omega}|_{1,\text{off}}$ in the penalty, but the SPICE estimators
\eqref{spice1} and \eqref{spice2}
are more amenable to theoretical analysis \cite{Lam09,Ravi11,Roth08}, and numerically they give similar results for i.i.d. data \cite{Roth08}.
It is worth noting that for i.i.d. data,  \eqref{spice1} requires $\sqrt{(p+s_p)(\log p)/n}=o(1)$ but \eqref{spice2} relaxes it to  $\sqrt{(1+s_p)(\log p)/n}=o(1)$. Similar requirements also hold for temporally dependent observations. Hence in this article, we only consider the SPICE estimator given in \eqref{spice2}.

\begin{thm}
\label{spice_thm2}
Suppose that $\mb{X}_{p\times n}$ is generated from \eqref{linear form} with all $e_i$ satisfying condition~\ref{C1} with the same $K$.
Uniformly on $\mb{\Omega}\in \mathcal{G}_2(s_p,v_0)$ and
$\{\mb{R}_{[k]}\}_{k=1}^p$ subject to \eqref{g_F,g_2},
for sufficiently large constant $M>0$ depending only on $v_0$ and $K$, if
	${\lambda_2}=Mu_1$ and $u_1=o(1/\sqrt{1+s_p})$
with $u_1$ defined in \eqref{u for subgauss},
then we have
\begin{align*}
\|\widehat{\mb{K}}_{\lambda_2}-\mb{K}\|_F&=O_P(u_1\sqrt{s_p}),
\label{diff spice corr PDD}
\\
\| \widehat{\mb{\Omega}}_{\lambda_2}-\mb{\Omega}  \|_2&=O_P(u_1\sqrt{1+s_p}),
\nonumber\\
\frac{1}{\sqrt{p}}\| \widehat{\mb{\Omega}}_{\lambda_2}-\mb{\Omega}  \|_F&=O_P\left(u_1\sqrt{1+s_p/p}\right). \nonumber 
\end{align*}
\end{thm}

Again by Remarks \ref{rmk1} and \ref{rmk2}, $u_1=O(\sqrt{(\log p)/n})$ is achievable for the short-memory case \eqref{short mem multi} and also for some long-memory cases, thus for such temporally dependent data Theorem~\ref{spice_thm2} gives the same convergence rates as those given in \cite{Roth08} for i.i.d. observations.
The condition $u_1=o(1/\sqrt{1+s_p})$ implies $s_p=o(u_1^{-2})=o(n/\log p)$, meaning that $\mb{\Omega}$ needs to be very sparse. Such a condition easily fails for many simple band matrices when $p\ge n$.

Under the irrepresentability condition, however, the sparsity requirement can be relaxed \cite{Ravi11}. In particular, define
$\mb{\Gamma}=\mb{R}\otimes \mb{R}$. By $(i,j)$-th row of $\mb{\Gamma}$ we refer to its $[i+(j-1)p]$-th row, and by $(k,\ell)$-th column to its $[k+(\ell-1)p]$-th column.
For any two subsets $T$ and $T'$ of $\{1,...,p\}{\times}\{1,...,p\}$,
denote $\mb{\Gamma}_{TT'}$ be the $\text{card}(T){\times} \text{card}(T')$ matrix with rows and columns of $\mb{\Gamma}$ indexed by $T$ and $T'$ respectively, where $\text{card}(T)$ denotes the cardinality of set $T$.
Let $S$ be the set of nonzero entries of $\mb{\Omega}$ and $S^c$ be the complement of $S$ in $\{1,...,p\}{\times}\{1,...,p\}$.
Define
$\kappa_{\mb{R}}=\| \mb{R} \|_1$
and
$\kappa_{\mb{\Gamma}}=\| \mb{\Gamma}_{SS}^{-1}    \|_1$.
Assume the following irrepresentability condition of \cite{Ravi11}:
\begin{equation}
\label{irr cond}
\max_{e\in S^c}\left | \mb{\Gamma}_{eS}\mb{\Gamma}_{SS}^{-1}  \right|_1 \le1-\beta
\end{equation}
for some $\beta\in (0,1]$.
Define $d$ to be the maximum number of nonzeros per row in $\mb{\Omega}$.
Then we have the following result.

\begin{thm}
\label{spice_thm3}
Let $r=\left(
0.5+2.5(1+8/\beta)\kappa_{\mb{\Gamma}}
\right)Mu_1v_0$,  where $u_1$ is defined in \eqref{u for subgauss}. Suppose that $\mb{X}_{p\times n}$ is generated from \eqref{linear form} with all $e_i$ satisfying condition~\ref{C1} with the same $K$.
Uniformly on $\mb{\Omega}\in \mathcal{G}_2(s_p,v_0)$ and
$\{\mb{R}_{[k]}\}_{k=1}^p$ subject to \eqref{g_F,g_2},
for sufficiently large constant $M>0$ depending on $v_0$ and $K$,
if $\lambda_2=8Mu_1/\beta
\le[6(1+\beta/8)d\max\{\kappa_{\mb{R}} \kappa_{\mb{\Gamma}}, \kappa_{\mb{R}}^3  \kappa_{\mb{\Gamma}}^2   \}]^{-1}$
and $u_1=o(\min\{1,[(1+8/\beta)\kappa_{\mb{\Gamma}}]^{-1}\})$, then
with probability tending to 1 we have
\begin{align*}
|\widehat{\mb{\Omega}}_{\lambda_2}-\mb{\Omega}|_{\infty} &\le r,
\\
\| \widehat{\mb{\Omega}}_{\lambda_2}-\mb{\Omega} \|_2 &\le  r \min\left\{d,\sqrt{p+s_p}\right\},
\\
\frac{1}{\sqrt{p}}\|\widehat{\mb{\Omega}}_{\lambda_2}-\mb{\Omega}  \|_F
&\le r\sqrt{1+s_p/p},
\end{align*}
and
$\hat\omega_{ij}^{(\lambda_2)}=0$ for all $(i,j)$ with $\omega_{ij}=0$.
If we further assume all nonzero entries of $\bd\Omega$ satisfy $|\omega_{ij}|>r$,
then with probability tending to 1,
$\textup{sign}(\hat\omega_{ij}^{(\lambda_2)})=\textup{sign}(\omega_{ij})$
for all $(i,j)$ where $\omega_{ij}\ne 0$.
\end{thm}

Consider the case
when $\beta$ remains constant and $\max\{\kappa_{\mb{R}},
\kappa_{\mb{\Gamma}}\}$ has a constant upper bound.
Then the conditions in Theorem~\ref{spice_thm3} about $\lambda_2$ and $u_1$ reduce to $\lambda_2=M' u_1$ and $u_1=o(1)$ with a constant $M'=8M/\beta$, and meanwhile we have
$
\| \widehat{\mb{\Omega}}_{\lambda_2}-\mb{\Omega} \|_2=O_P(u_1d).
$
Then the desired result of $
\| \widehat{\mb{\Omega}}_{\lambda_2}-\mb{\Omega} \|_2=o_P
(1)$ is achieved under a relaxed sparsity condition $d=o(u_1^{-1})$.
If $d^2 > 1+s_p$, then  $s_p=o(u_1^{-2})$ and the condition of Theorem~\ref{spice_thm2} satisfies. Hence $
\| \widehat{\mb{\Omega}}_{\lambda_2}-\mb{\Omega} \|_2=O_P(u_1\sqrt{\min\{d^2,1+s_p\}})=o_P(1),
$
which is the better rate between those from Theorems~\ref{spice_thm2} and \ref{spice_thm3}.

\section{Extension to heavy tail data}\label{s: heavy tails}
In this section, we generalize the theoretical results for sub-Gaussian data to the cases
when all the basis random variables $\{e_i\}$ have the generalized sub-exponential tails under 
condition \ref{C2}
or the polynomial-type tails under condition \ref{C3}.
Define
\be\label{u for subexp}
u_2=\max\left\{
(\log p)^{1+2/\vartheta}g_2/n,(\log p)^{1+2/\vartheta}(g_F/n)^{1/2}
\right\},
\ee
and
\be\label{u for poly}
u_3=\max\left\{
p^{(2+2C)/k} g_2/n, p^{(4+2C)/k}(g_F/n)^{1/2}
\right\}
\ee
with an arbitrary constant $C>0$.
The quantities $u_2$ and $u_3$ will substitute $u_1$ in characterizing the matrix estimation convergence rates under the tail conditions ~\ref{C2} and \ref{C3}, respectively. The first term in either $u_2$ or $u_3$ can be dropped if $\bd{\mu}_p$ is known thus no need to be estimated.

\begin{thm}[Generalized sub-exponential tails]\label{subexp thm}
Theorems~\ref{thm1}, \ref{thm2}, \ref{clime_thm1}$-$\ref{spice_thm3}, and Corollary~\ref{cor1} hold
with condition~\ref{C1}, parameter $K$, and $u_1$ replaced by 
condition \ref{C2}, parameters $\{K,\vartheta\}$, and $u_2$, respectively.
\end{thm}

\begin{thm}[Polynomial-type tails]\label{poly thm}
Theorems~\ref{thm1}, \ref{thm2}, \ref{clime_thm1}$-$\ref{spice_thm3}, and Corollary~\ref{cor1},
except those mean-square convergence results therein,
hold with condition~\ref{C1}, parameter $K$, and $u_1$ replaced by condition \ref{C3}, parameters $\{k,\eta_k\}$, and $u_3$, respectively.
\end{thm}

For data with polynomial-type tails, the mean-square convergence results may require higher order moment conditions, thus are not pursued here.

A referee pointed out a potential connection to the recent work of \cite{Chen16}, where the estimation of $\mb{\Omega}\bd{b}$ with $\bd{b}\in \mathbb{R}^p$ was considered for  high-dimensional mean-zero stationary processes given in \eqref{infty sum} with a temporal dependence satisfying PDD given in $\eqref{PDD}$. Our exploration shows that their concentration inequalities (see page 3 of their Supplementary Material) can be used to obtain the same convergence rates for our concerned estimators to their sub-Gaussian time series. Together with the concentration inequalities in \cite{Wu16}, their inequalities also can be applied to their time series with the generalized sub-exponential tails, but result in slower convergence rates.  If applied to their time series data with polynomial-type tails, however, their inequalities seem to yield faster convergence rates than ours. We leave the details to interested readers. Note that it is not clear if the concentration inequalities in \cite{Chen16} and \cite{Wu16} are applicable to the more general temporal dependence that we consider here in this article.

\section{Numerical Results}\label{s:numeric}
\subsection{Gap-block cross-validation}

For tuning parameter selection, we propose a gap-block cross-validation method that includes the following steps:
\begin{enumerate}[1.]
	\item  Split the data $\mb{X}_{p\times n}$ into $H_1\ge 4$ approximately equal-sized non-overlapping blocks  $\mb{X}^*_{i}$, $i=1,\dots,H_1$, such that $\mb{X}_{p\times n}=(\mb{X}_{1}^*,\mb{X}_{2}^*,\dots,\mb{X}_{H_1}^*)$. For each $i$, set aside block $\mb{X}_i^*$ that will be used as the validation data, and use the remaining data after further dropping the neighboring blocks at both sides of $\mb{X}_i^*$ as the training data that are denoted by $\mb{X}_i^{**}$.
	\label{stp1}

	\item Randomly sample $H_2$ blocks $\mb{X}_{H_1+1}^*,\dots,\mb{X}_{H_1+H_2}^*$ from $\mb{X}_{p\times n}$, where $\mb{X}_{H_1+j}^*$ consists of $\lceil{n/H_1}\rceil$ consecutive columns of $\mb{X}_{p\times n}$ for each $j=1,\dots,H_2$. Note that these sampled blocks can overlap. For each $i =H_1+1,\dots,H_1+H_2$, set aside block $\mb{X}_i^*$  as the validation data, and use the remaining data by further excluding the $\lceil{n/H_1}\rceil$ columns at both sides of $\mb{X}_i^*$ from $\mb{X}_{p\times n}$ as the training data that are denoted by $\mb{X}_i^{**}$.
	\label{stp2}
	
	\item Let $H=H_1+H_2$.
Select the optimal values of tuning parameters $\tau_1,\tau_2,\lambda_1$ and $\lambda_2$ among their corresponding prespecified candidate values $\{\tau_{1j}\}_{j=1}^J$, $\{\tau_{2j}\}_{j=1}^J$, $\{\lambda_{1j}\}_{j=1}^J$ and $\{\lambda_{2j}\}_{j=1}^J$, and denote them respectively by
	\begin{align*}
	\hat{\tau}_1&=\mathop{\arg\min}_{1\le j\le J}\frac{1}{H}\sum_{i=1}^H\| S_{{\tau}_{1j}}(\widehat{\mb{\Sigma}}_i^{**})-\widehat{\mb{\Sigma}}_i^*   \|_F^2,\\
		\hat{\tau}_2&=\mathop{\arg\min}_{1\le j\le J}\frac{1}{H}\sum_{i=1}^H\| S_{{\tau}_{2j}}(\widehat{\mb{R}}_i^{**})-\widehat{\mb{R}}_i^*   \|_F^2,\\
\hat{\lambda}_1&=	\mathop{\arg\min}_{1\le j\le J}\frac{1}{H}\sum_{i=1}^H	
\left[\text{tr}(\widehat{\mb{\Omega}}_{\varepsilon,\lambda_{1j},i}^{**}\widehat{\mb{\Sigma}}_i^*)-\log\det(\widehat{\mb{\Omega}}_{\varepsilon,\lambda_{1j},i}^{**})\right],\\
\hat{\lambda}_2&=	\mathop{\arg\min}_{1\le j\le J}\frac{1}{H}\sum_{i=1}^H	
\left[\text{tr}(\widehat{\mb{\Omega}}_{\lambda_{2j},i}^{**}\widehat{\mb{\Sigma}}_i^*)-\log\det(\widehat{\mb{\Omega}}_{\lambda_{2j},i}^{**})\right],
	\end{align*}
where $\widehat{\mb{\Sigma}}_i^{*}$ 	and $\widehat{\mb{R}}_i^*$ are obtained from
$\mb{X}_i^*$, 
$\widehat{\mb{\Sigma}}_i^{**}$ and $\widehat{\mb{R}}_i^{**}$ are obtained from $\mb{X}_i^{**}$,
and $\widehat{\mb{\Omega}}_{\varepsilon,\lambda_{1j},i}^{**}$ and $\widehat{\mb{\Omega}}_{\lambda_{2j},i}^{**}$
are the CLIME and SPICE estimators, respectively, obtained from $\mb{X}_i^{**}$.

	\label{stp3}
\end{enumerate}

In the above cross-validation (CV), due to lack of independent observations,
we use gap blocks, each of size $\approx \lceil{n/H_1}\rceil$,  to separate training and validation datasets 
so that they are nearly uncorrelated.
The idea of using gap blocks has been employed by the $hv$-block CV of \cite{Raci00} for linear models with dependent data.
Similar to the $k$-fold CV for i.i.d. data, Step~\ref{stp1} guarantees all observations are used for both training and validation, but is limited due to
the constrain of keeping the temporal ordering of the observations.
Step~\ref{stp2} allows more data splits. This is particularly useful when Step~\ref{stp1} only allows a small number of data splits due to large-size of the gap block and/or limited sample size $n$.
Step~\ref{stp2} is inspired by the commonly used repeated random subsampling CV for i.i.d. observations \cite{Syed12}.
The above loss functions for selecting tuning parameters are widely used in the literature \cite{Bick08a,CLL11,CLZ16, Roth09}.
The theoretical justification for the gap-block CV remains open. In our numerical examples, we simply set $H_1=H_2=10$, mimicing the $10$-fold CV recommended by \cite{Fang16,Hast09}. 
Our simulation studies show that the method performs well for temporally dependent data.

\subsection{Simulation studies}\label{Simulations}
We evaluate the numerical performance of the hard and soft thresholding estimators for large correlation matrix and the CLIME and SPICE estimators for large precision matrix.
We generate Gaussian data with zero mean and
covariance matrix $\mb{\Sigma}$ or precision matrix $\mb{\Omega}$ from one of the following four models:
\begin{enumerate}
\renewcommand{\theenumi}{\textup{Model \arabic{enumi}}}
\item \label{Model 1}$\sigma_{ij}=0.6^{|i-j|}$;
\item \label{Model 2}$\sigma_{ii}=1$, $\sigma_{i,i+1}=\sigma_{i+1,i}=0.6$, $\sigma_{i,i+2}=\sigma_{i+2,i}=0.3$, and $\sigma_{ij}=0$ for $|i-j|\ge 3$;
\item \label{Model 3}$\omega_{ij}=0.6^{|i-j|}$;
\item \label{Model 4}$\omega_{ii}=1$, $\omega_{i,i+1}=\omega_{i+1,i}=0.6$, $\omega_{i,i+2}=\omega_{i+2,i}=0.3$, and $\omega_{ij}=0$ for $|i-j|\ge 3$.
\end{enumerate}
Similar models have been considered in \cite{Bick08a,CLL11,CLZ16,Roth08,Roth09}.
For the temporal dependence,
we set $\textup{Corr}(X_{ki},X_{\ell j})=\varLambda_{k\ell}^{ij}\rho_{k\ell}$
with
 \be
\label{simu rate matrix}
\varLambda_{k\ell}^{ij}=
(|i-j|+1)^{-\alpha}, \quad 1\le i, j \le n,
\ee
so that
$\rho_{[k]}^{ij}\sim |i-j|^{-\alpha}$.
It is computationally expensive to simulate data $\bd{X}_{pn}$ directly from a multivariate Gaussian random number generator because of the large dimension of its covariance matrix $\textup{Cov}(\bd{X}_{pn})$. Instead, we simulate data using the method of \cite{Boch07}, which approximately satisfy
\eqref{simu rate matrix} (see details in Supplementary Material).

Simulations are conducted
with sample size $n=200$, variable dimension $p$ ranging from 100 to 400, and 100 replications under each setting, for which $\alpha$ varies from 0.1 to 2. The i.i.d. case is also considered, for which an ordinary 10-fold CV is implemented.
For each simulated data set, we choose the optimal tuning parameter from a set of 50 specified values (see Section S.4.1 in Supplementary Material).  
The CLIME and SPICE are computed by the R packages {\tt flare} \cite{Li15} and {\tt QUIC} \cite{Hsie14}, respectively.
For CLIME, we use the default perturbation of {\tt flare} with $\varepsilon=n^{-1/2}$.

The estimation performance is measured by both the spectral norm and the Frobenius norm. True-positive rate (TPR) and false-positive rate (FPR) are used for evaluating sparsity recovering for the correlation matrix:
\[
\text{TPR}=\frac{\#\{(i,j):s_{\tau}(\hat{\rho}_{ij})\ne 0 \ \text{and}\ \rho_{ij}\ne 0,i\ne j\}}{\#\{(i,j):\rho_{ij}\ne 0,i\ne j\}},
\]
\[
\text{FPR}=\frac{\#\{(i,j):s_{\tau}(\hat{\rho}_{ij})\ne 0 \ \text{and} \ \rho_{ij}=0,i\ne j\}}{\#\{(i,j):\rho_{ij}= 0,i\ne j\}}.
\]
Similar quantities are also reported for the precision matrix. The TPR and FPR are not provided for Models 1 and 3.

Simulation results are summarized in Tables~\ref{model 1_2}-\ref{TPR/FPR}. In all setups,
the sample correlation matrix and the inverse of sample covariance matrix (whenever possible) perform the worst.
It is not surprising that the performance of all the regularized estimators generally is better for weaker temporal dependence or smaller $p$. The soft thresholding method performs slightly better than the hard thresholding method in terms of matrix losses for small $\alpha$ and slightly worse for large $\alpha$, and always has higher TPRs but bigger FPRs.
The CLIME estimator performs similarly as the SPICE estimator in matrix norms,
but generally yields lower FPRs.

We notice that the SPICE algorithm in the R package {\tt QUIC} is much faster than the CLIME algorithm in the R package {\tt flare}
by
using a single computer core. However, the column-by-column estimating nature of CLIME can speed up using parallel computing on multiple cores.

\begin{table}[ht!]
\caption{Comparison of average (SD) matrix losses for correlation matrix estimation}
\centering
\scalebox{0.88}{
\begin{tabular}{c l r r r r r r r r}
\hline\hline\\[-2ex]
$p$  & $\ \alpha$
&& $\widehat{\mb{R}}\quad\ $ & Hard & Soft
&& $\widehat{\mb{R}}\quad\ $ & Hard & Soft
\\
\hline\\[-2ex]
& & &\multicolumn{3}{c}{Spectral norm}
&& \multicolumn{3}{c}{Frobenius norm}\\
\cline{4-6}\cline{8-10}\\[-2ex]
& & &\multicolumn{7}{c}{Model 1}\\
100
&0.1
&
&13.7(1.68) &2.8(0.09) &2.6(0.07)
&
&22.6(1.08) &9.9(0.28) &8.7(0.24)
\\

&0.25
&
&10.5(1.59) &2.4(0.15) &2.4(0.08)
&
&17.4(0.95) &8.1(0.42) &7.5(0.26)
\\

&0.5
&
&7.8(1.14) &2.0(0.15) &2.2(0.08)
&
&14.3(0.69) &6.8(0.33) &6.6(0.23)
\\

&1
&
&4.2(0.45)& 1.5(0.10) &1.7(0.08)
&
&9.9(0.29)& 5.2(0.23)& 5.1(0.20)
\\

&2
&
&2.6(0.24) & 1.1(0.09) &1.4(0.08)
&
&7.5(0.17) & 3.9(0.15) & 4.0(0.19)
\\

& i.i.d.
&
&2.4(0.18) &1.0(0.08) &1.3(0.08)
&
&7.0(0.15) &3.5(0.13) &3.7(0.15)
\\

200
&0.1
&
&27.2(2.69) &2.9(0.05) &2.8(0.04)
&
&45.6(1.54) &14.5(0.25) &13.1(0.22)
\\

&0.25
&
& 20.6(2.54)&2.5(0.14) &2.5(0.06)
&
&35.0(1.39) &12.2(0.56) &11.4(0.29)
\\

&0.5
&
&15.2(1.77) &2.2(0.12) &2.3(0.06)
&
&28.7(0.99) &10.2(0.40) &10.1(0.25)
\\

&1
&
&7.8(0.64)& 1.6(0.08)& 1.9(0.06)
&
&20.1(0.35)&  7.9(0.24) & 7.9(0.21)
\\

&2
&
&4.3(0.24)& 1.3(0.08)& 1.6(0.06)
&
&15.1(0.15)&  5.9(0.19)&  6.3(0.18)
\\

& i.i.d.
&
&3.8(0.22) &1.1(0.07) &1.5(0.06)
&
& 14.1(0.15)&5.3(0.14) &5.8(0.17)
\\

300
&0.1
&
&40.6(3.39) &3.0(0.03) &2.8(0.03)
&
&68.5(1.88) & 18.0(0.21)&16.5(0.24)
\\

&0.25
&
&30.9(3.23) &2.6(0.11) &2.6(0.04)
&
&52.6(1.75) &15.4(0.63) &14.5(0.30)
\\

&0.5
&
&22.5(2.16) & 2.3(0.12)& 2.4(0.04)
&
&43.2(1.16) &12.8(0.47)& 12.9(0.27)
\\

&1
&
&11.2(0.79)&  1.7(0.05)&  2.0(0.05)
&
&30.2(0.42)&  9.9(0.21)& 10.1(0.25)
\\

&2
&
&5.8(0.27)& 1.3(0.08)& 1.7(0.05)
&
&22.8(0.16)&  7.5(0.25)&  8.2(0.19)
\\

& i.i.d.
&
&5.0(0.20) &1.2(0.08) &1.6(0.05)
&
&21.2(0.15) &6.7(0.12)&7.5(0.17)
\\

400
&0.1
&
&54.2(4.01) &3.0(0.02) &2.9(0.02)
&
&91.7(2.17) &20.9(0.17) &19.4(0.22)
\\

&0.25
&
&41.0(3.88) &2.7(0.09) &2.7(0.04)
&
&70.1(2.09) &18.4(0.61) &17.1(0.29)
\\

&0.5
&
&29.8(2.62) & 2.3(0.12)&  2.5(0.04)
&
&57.7(1.38)& 15.2(0.59)& 15.3(0.30)
\\

&1
&
&14.6(0.91)&  1.7(0.05) & 2.1(0.04)
&
&40.3(0.48)& 11.6(0.22) &12.1(0.20)
\\

&2
&
&7.2(0.26)& 1.4(0.07)& 1.8(0.04)
&
&30.4(0.16)&  9.0(0.27)&  9.8(0.23)
\\

& i.i.d.
&
&6.0(0.21) &1.2(0.08) &1.6(0.05)
&
&28.2(0.15) &7.9(0.14) &8.9(0.17)
\\
\\
& & &\multicolumn{7}{c}{Model 2}\\
100
&0.1
&
&13.8(1.71) &1.8(0.04) &1.6(0.04)
&
&22.6(1.05) &8.7(0.29) &7.7(0.22)
\\

&0.25
&
&10.5(1.61) &1.5(0.18) &1.4(0.09)
&
&17.5(0.94) &6.7(0.48) &6.5(0.24)
\\

&0.5
&
&7.8(1.10) & 1.2(0.17)&1.3(0.07)
&
&14.3(0.66) &5.2(0.34) &5.6(0.21)
\\

&1
&
&4.2(0.40) &0.7(0.09) &1.0(0.05)
&
& 10.0(0.27)&4.0(0.17) &4.1(0.16)
\\

&2
&
&2.5(0.18) &0.6(0.05) &0.8(0.04)
&
&7.5(0.14) &2.6(0.25) &3.2(0.13)
\\

& i.i.d.
&
&2.3(0.15) &0.5(0.07) &0.7(0.04)
&
&7.0(0.13) &2.0(0.23) &2.8(0.12)
\\

200
&0.1
&
&27.2(2.62) &1.8(0.02) &1.7(0.03)
&
&45.6(1.51) &12.9(0.28) &11.6(0.21)
\\

&0.25
&
&20.6(2.29) &1.6(0.15) &1.5(0.07)
&
&35.0(1.29) &10.3(0.56) &9.9(0.27)
\\

&0.5
&
&15.1(1.58) &1.3(0.14) &1.4(0.05)
&
&28.8(0.88) &7.9(0.43) &8.6(0.21)
\\

&1
&
& 7.7(0.57)&0.8(0.10) &1.1(0.04)
&
&20.1(0.34) &5.8(0.15) &6.5(0.20)
\\

&2
&
&4.2(0.18) & 0.6(0.05)&0.9(0.04)
&
&15.2(0.14) &4.2(0.30) &5.0(0.13)
\\

& i.i.d.
&
&3.6(0.16) &0.6(0.06) &0.8(0.04)
&
&14.1(0.14) & 3.2(0.23)&4.4(0.12)
\\

300
&0.1
&
&40.8(3.54) &1.8(0.05) &1.7(0.02)
&
&68.7(1.84) &16.0(0.27) &14.6(0.24)
\\

&0.25
&
&30.8(2.95) &1.7(0.17) &1.6(0.13)
&
&52.6(1.62) & 13.2(0.69)&12.5(0.28)
\\

&0.5
&
&22.4(2.04) & 1.4(0.12)&1.4(0.09)
&
&43.3(1.10) &10.1(0.57) &10.9(0.25)
\\

&1
&
&11.1(0.73) &0.9(0.08) &1.1(0.03)
&
&30.3(0.41) &7.3(0.16) &8.3(0.20)
\\

&2
&
&5.6(0.22) &0.6(0.04) &0.9(0.04)
&
&22.8(0.14) & 5.5(0.29)&6.5(0.18)
\\

& i.i.d.
&
&4.7(0.15) &0.6(0.05) &0.8(0.03)
&
&21.2(0.13) &4.1(0.21) &5.7(0.12)
\\

400
&0.1
&
&54.0(3.61) &1.8(0.04) &1.7(0.02)
&
&91.7(1.97) &18.6(0.16) &17.2(0.15)
\\

&0.25
&
&41.1(3.58) &1.7(0.09) &1.7(0.12)
&
&70.2(1.89) & 15.8(0.63)&14.9(0.33)
\\

&0.5
&
&29.7(2.53) &1.5(0.17) &1.5(0.08)
&
&57.7(1.29) &12.1(0.62) &13.0(0.24)
\\

&1
&
&14.5(0.86) & 0.9(0.08)&1.1(0.03)
&
&40.4(0.46) &8.6(0.16) &10.0(0.23)
\\

&2
&
&7.0(0.26) &0.7(0.04) &0.9(0.03)
&
&30.4(0.14) &6.6(0.26) &7.7(0.15)
\\

& i.i.d.
&
&5.7(0.18) &0.6(0.05) &0.9(0.03)
&
&28.3(0.13) &4.9(0.21) &6.8(0.12)
\\
\hline
\end{tabular}
}
\label{model 1_2}
\end{table}

\begin{table}[ht!]
\caption{Comparison of average (SD) matrix losses for precision matrix estimation}
\centering
\scalebox{0.88}{
\begin{tabular}{c l r r r r r r r rr}
\hline\hline\\[-2ex]
$p$  & $\ \alpha$
&& $\widehat{\mb\Sigma}^{-1}\quad\ $ & CLIME & SPICE
&& $\widehat{\mb\Sigma}^{-1}\quad\ $& CLIME & SPICE
\\
\hline\\[-2ex]
& & &\multicolumn{3}{c}{Spectral norm}
&& \multicolumn{3}{c}{Frobenius norm}
\\
\cline{4-6} \cline{8-10}\\[-2ex]
& & &\multicolumn{7}{c}{Model 3}
\\
100
&0.1
&
&381.7(40.07) &4.9(0.26)&5.7(0.53)
&
&850.5(38.22) &28.8(1.54) &27.1(1.46)

\\

&0.25
&
&97.6(9.23) &1.8(0.09)&2.2(0.08)
&
&214.6(9.38) &9.5(0.34)&9.3(0.20)

\\

&0.5
&
&43.3(4.60) &2.4(0.09)&2.7(0.06)
&
&93.9(4.36) &7.7(0.15)&8.6(0.15)

\\

&1
&
&21.8(2.74) &2.6(0.06)&2.9(0.04)
&
&45.4(2.73) &8.0(0.19)&9.2(0.15)
\\

&2
&
&14.1(1.80) &2.7(0.05)&2.9(0.04)
&
&28.9(1.86) &8.0(0.20)&9.1(0.14)

\\

& i.i.d.
&
&12.6(1.66) &2.5(0.06)&2.8(0.04)
&
&25.5(1.56) &7.4(0.20)&8.6(0.15)

\\

200
&0.1
&
&N/A &6.2(0.38)&5.8(0.48)
&
&N/A &49.6(2.46)&38.4(1.48)

\\

&0.25
&
&N/A&2.1(0.12) &2.4(0.06)
&
&N/A&14.8(0.52)&13.7(0.18)

\\

&0.5
&
&N/A&2.6(0.07)&2.8(0.04)
&
&N/A &11.9(0.18)&12.8(0.12)
\\

&1
&
&N/A&2.9(0.05)&3.1(0.03)
&
&N/A &12.4(0.23)&13.7(0.14)

\\

&2
&
&N/A&2.9(0.04)&3.1(0.02)
&
&N/A &12.6(0.21)&13.8(0.09)

\\

& i.i.d.
&
&N/A&2.7(0.04)&3.0(0.02)
&
&N/A &11.6(0.24)&13.3(0.14)

\\

300
&0.1
&
&N/A &5.3(0.36)&5.9(0.45)
&
&N/A &51.2(2.85)&47.1(1.48)

\\

&0.25
&
&N/A&2.4(0.11) &2.4(0.05)
&
&N/A &18.0(0.36)&17.1(0.18)

\\

&0.5
&
&N/A&2.8(0.07)&2.9(0.03)
&
&N/A&15.7(0.27) &15.9(0.13)
\\

&1
&
&N/A &3.0(0.04)&3.1(0.02)
&
&N/A&15.9(0.28)&17.1(0.12)

\\

&2
&
&N/A&3.0(0.03)&3.1(0.01)
&
&N/A &16.1(0.22)&17.3(0.09)
\\

& i.i.d.
&
&N/A &2.8(0.04)&3.1(0.02)
&
&N/A&15.0(0.26)&16.8(0.11)

\\

400
&0.1
&
&N/A&5.8(0.44)&6.0(0.37)
&
&N/A &63.9(4.29)&54.7(1.60)

\\

&0.25
&
&N/A&2.6(0.08)&2.5(0.05)
&
&N/A &20.8(0.22)&20.0(0.19)
\\

&0.5
&
&N/A &2.9(0.06)&2.9(0.03)
&
&N/A &19.0(0.31)&18.6(0.12)

\\

&1
&
&N/A&3.0(0.04)&3.1(0.02)
&
&N/A &19.0(0.32)&19.9(0.13)
\\

&2
&
&N/A &3.1(0.03)&3.2(0.01)
&
&N/A &19.0(0.24)&20.2(0.10)

\\

& i.i.d.
&
&N/A&2.9(0.04) &3.1(0.01)
&
&N/A &17.9(0.31)&19.7(0.10)

\\
\\
&&& \multicolumn{7}{c}{Model 4}\\
100
&0.1
&

&355.4(37.62) &4.9(0.40)&5.9(0.72)
&
&829.5(35.78) &28.0(2.05)&26.5(1.68)
\\

&0.25
&
&91.1(8.42)&1.9(0.31)&1.7(0.19)
&
&209.0(8.63) &8.2(1.03) &7.3(0.30)
\\

&0.5
&
&40.7(4.29) &1.1(0.10)&1.4(0.07)
&
&91.6(3.96) &4.7(0.17)&5.8(0.19)
\\

&1
&
&20.5(2.44)&1.3(0.07)&1.6(0.06)
&
&44.4(2.44) &5.1(0.26)&6.2(0.21)
\\

&2
&
&13.3(1.62) &1.4(0.07)&1.6(0.05)
&
&28.3(1.70) &5.3(0.25)&6.3(0.17)
\\

& i.i.d.
&
&11.8(1.44)&1.2(0.06)&1.4(0.05)
&
&25.0(1.37) &4.6(0.24)&5.7(0.18)
\\

200
&0.1
&
&N/A&5.4(0.50)&5.6(0.61)
&
&N/A &41.4(2.89)&33.9(1.61)
\\

&0.25
&
&N/A  &1.8(0.19)&1.6(0.14)
&
&N/A &11.5(0.59)&10.5(0.18)
\\

&0.5
&
&N/A &1.4(0.11)&1.7(0.04)
&
&N/A &8.5(0.32)&9.6(0.17)
\\

&1
&
&N/A&1.6(0.06)&1.8(0.03)
&
&N/A&9.1(0.38)&10.5(0.21)
\\

&2
&
&N/A &1.6(0.05)&1.8(0.03)
&
&N/A &9.2(0.32)&10.8(0.17)
\\

& i.i.d.
&
&N/A &1.4(0.06)&1.7(0.03)
&
&N/A &7.8(0.34)&9.9(0.17)
\\

300
&0.1
&
&N/A &6.0(0.54)&5.6(0.67)
&
&N/A &54.7(4.26)&39.8(1.58)
\\

&0.25
&
&N/A&1.6(0.12)&1.6(0.14)
&
&N/A &14.0(0.30)&13.2(0.13)
\\

&0.5
&
&N/A&1.8(0.07)&1.8(0.04)
&
&N/A&13.1(0.51)&12.5(0.20)
\\

&1
&
&N/A&1.9(0.06) &1.9(0.03)
&
&N/A &13.1(0.53)&13.8(0.20)
\\

&2
&
&N/A &1.8(0.05)&2.0(0.03)
&
&N/A &12.6(0.39)&14.2(0.20)
\\

& i.i.d.
&
&N/A&1.5(0.05)&1.8(0.02)
&
&N/A &10.5(0.38)&13.2(0.19)
\\

400
&0.1
&
&N/A&5.1(0.46)&5.4(0.62)
&
&N/A&54.4(4.12)&44.6(1.43)
\\

&0.25
&
&N/A&1.8(0.09)&1.7(0.14)
&
&N/A&17.5(0.28)&15.5(0.11)
\\

&0.5
&
&N/A&2.0(0.06)&1.9(0.03)
&
&N/A&17.3(0.55)&14.9(0.19)
\\

&1
&
&N/A&2.0(0.06)&2.0(0.02)
&
&N/A&16.7(0.59)&16.5(0.20)
\\

&2
&
&N/A&1.9(0.05)&2.0(0.02)
&
&N/A&15.9(0.50)&17.1(0.20)
\\

& i.i.d.
&
&N/A &1.7(0.06)&1.9(0.02)
&
&N/A &13.5(0.48)&16.0(0.20)
\\
\hline
\end{tabular}
}
\label{model 3_4}
\end{table}

\begin{table}[ht!]
\caption{Comparison of average (SD) TPR(\%)/FPR(\%) for Models 2 and 4}
\centering
\scalebox{0.88}{
\begin{tabular}{c l l l l  l l}
\hline\hline\\[-2ex]
$p$  & $\ \alpha$
&& Model 2, Hard && Model 2, Soft
\\
\cline{1-2}\cline{4-4}\cline {6-6}\\[-2ex]
100
&0.1
&
&10.86(4.35)/0.02(0.03)  &&  54.19(4.41)/4.98(1.26)
 \\

&0.25
&
&35.16(5.43)/0.07(0.06)  && 70.72(3.96)/6.10(1.16)
\\

&0.5
&
&48.43(3.76)/0.06(0.06)  && 80.43(3.19)/6.75(1.19)
  \\

&1
&
&60.92(4.25)/0.02(0.03)  && 94.34(2.12)/7.23(1.39)
 \\

&2
&
&83.93(4.08)/0.04(0.05)  && 99.33(0.73)/7.47(1.57)
  \\

& i.i.d.
&
&93.42(2.63)/0.13(0.09)  && 99.91(0.21)/11.42(1.82)
 \\

200
&0.1
&
&\hspace{0.5ex}  5.57(2.93)/0.00(0.00)  && 45.91(3.86)/2.40(0.55)
 \\

&0.25
&
&28.31(4.75)/0.02(0.02)  && 64.71(3.23)/3.20(0.69)
 \\

&0.5
&
&44.48(3.02)/0.02(0.02)  && 74.38(2.42)/3.40(0.59)
  \\

&1
&
&57.45(2.14)/0.01(0.01)  && 91.40(2.11)/3.84(0.81)
  \\

&2
&
& 79.04(3.66)/0.02(0.01) && 98.71(0.67)/3.73(0.58)
 \\

& i.i.d.
&
&90.74(2.68)/0.07(0.05)  && 99.68(0.31)/6.64(0.65)
  \\

300
&0.1
&
&\hspace{0.5ex} 4.15(2.50)/0.00(0.00)  && 40.61(3.94)/1.50(0.43)
 \\

&0.25
&
&24.28(4.85)/0.01(0.01)  && 61.27(2.70)/2.13(0.42)
 \\

&0.5
&
&41.75(3.51)/0.01(0.01)  && 71.65(2.51)/2.43(0.47)
 \\

&1
&
&55.42(2.10)/0.00(0.00)  && 89.41(1.80)/2.61(0.44)
  \\

&2
&
&74.39(3.23)/0.01(0.01)  && 98.11(0.69)/2.49(0.57)
 \\

& i.i.d.
&
&88.97(2.29)/0.04(0.02)  &&99.57(0.34)/4.77(0.84)
\\

400
&0.1
&
&\hspace{0.5ex} 2.65(1.29)/0.00(0.00)  && 36.80(2.27)/1.02(0.23)
\\

&0.25
&
&20.81(3.74)/0.01(0.00)  && 58.30(2.86)/1.54(0.35)
 \\

&0.5
&
&40.14(3.58)/0.01(0.01)  && 68.74(2.06)/1.68(0.35)
  \\

&1
&
&53.82(1.65)/0.00(0.00)  && 87.51(1.87)/1.80(0.40)
 \\

&2
&
&72.19(2.58)/0.00(0.00)  && 97.79(0.66)/1.97(0.22)
  \\

& i.i.d.
&
&87.51(1.65)/0.03(0.01)  && 99.38(0.30)/3.93(0.40)
  \\
  \\[-2ex]
  \cline{4-6}\\[-2ex]
&&
& Model 4, CLIME && Model 4, SPICE
\\
\cline{4-4}\cline {6-6}\\[-2ex]
 100
&0.1
&
&91.28(2.76)/25.49(2.37)&&82.99(2.76)/28.97(1.04)  \\

&0.25
&
&92.65(2.35)/17.82(1.84)&&90.93(2.19)/29.68(1.31) \\

&0.5
&
&95.30(1.73)/17.80(1.47)&&96.00(1.54)/31.58(1.49)  \\

&1
&
 &98.47(0.90)/14.37(1.21)&&99.24(0.66)/30.65(1.49)  \\

&2
&
&99.71(0.36)/11.99(1.27) &&99.94(0.17)/27.77(1.34)  \\

& i.i.d.
&
&99.91(0.20)/16.21(1.63) &&99.99(0.07)/31.40(1.28) \\

200
&0.1
&
&82.24(2.70)/12.72(0.64)&&76.07(1.95)/17.78(0.56)  \\

&0.25
&
&84.83(2.28)/15.70(2.62)&&84.75(1.90)/18.87(0.59) \\

&0.5
&
&89.55(2.39)/13.21(3.00)&&91.65(1.45)/20.07(0.64)  \\

&1
&
&93.81(1.52)/7.27(0.58)&&97.12(0.97)/19.07(0.85)  \\

&2
&
&97.77(0.97)/4.86(0.55)&&99.31(0.42)/16.25(0.81)  \\

& i.i.d.
&
&99.56(0.36)/7.24(0.79)&&99.88(0.18)/19.42(0.81)  \\

300
&0.1
&
&82.60(3.59)/12.71(2.55)&&71.66(1.71)/13.05(0.34)  \\

&0.25
&
&77.62(2.62)/14.39(2.62)&&81.09(1.71)/14.06(0.39)  \\

&0.5
&
&82.23(2.48)/14.33(3.57)&&88.71(1.44)/14.98(0.42)  \\

&1
&
&86.84(2.58)/4.71(0.67)&&94.87(1.02)/14.20(0.54)  \\

&2
&
&94.88(1.38)/2.84(0.41)&&98.27(0.68)/11.59(0.65)  \\

& i.i.d.
&
&98.83(0.49)/4.89(0.58) &&99.56(0.29)/14.32(0.70) \\

400
&0.1
&
&83.04(2.84)/14.91(2.84)&&68.51(1.49)/10.36(0.24) \\

&0.25
&
&76.76(3.46)/15.11(3.40)&&78.50(1.41)/11.41(0.32)  \\

&0.5
&
&78.58(2.35)/15.67(3.64)&&86.19(1.44)/12.20(0.35)  \\

&1
&
&79.44(3.05)/4.40(0.77)&&92.85(1.09)/11.55(0.41) \\

&2
&
&90.47(2.32)/1.92(0.35)&&96.68(0.85)/8.97(0.55)  \\

& i.i.d.
&
&97.63(0.82)/3.50(0.52)&&99.09(0.39)/11.34(0.60)  \\

\hline
\end{tabular}
}
\label{TPR/FPR}
\end{table}

\subsection{rfMRI data analysis}\label{subs:real data}

Here we analyze a rfMRI data set for the estimation of brain functional connectivity. The preprocessed rfMRI data of a healthy young woman are provided by the WU-Minn Human Connectome Project (www.humanconnectome.org).
The original data consist of 1,200 temporal brain images and each image contains 229,404 brain voxels with size $2{\times} 2{\times}2\mbox{ mm}^3$. We discard the first 10 images due to concerns of early nonsteady magnetization. For the ease of implementation, we use a grid-based method \cite{Srip14} to reduce the image dimension to 907 functional brain nodes that are placed in a regular three-dimensional grid spaced at 12-mm intervals throughout the brain. Each node consists of a 3-mm voxel-center-to-voxel-center radius pseudosphere, which encompasses 19 voxels, and the time series at the node is a spatially averaged time series of these 19 voxels. 
The temporal dependence of the 907 time series
is approximated by 
the PDD model~\eqref{PDD}
with $C_0=1$ and $\alpha=0.30$
(see Figure~\ref{autocorrelation}).

The functional connectivity between two brain nodes can be evaluated by
either correlation or partial correlation, here we follow the convention by simply
calling them the marginal connectivity and the direct connectivity, respectively. For the marginal connectivity,
we only apply the hard thresholding method for estimating the correlation matrix, which usually yields less number of false discoveries than the soft thresholding. We find that 1.47\% of all the pairs of nodes are connected with a threshold value of 0.12 to the sample correlations.
For the direct connectivity, we calculate the estimated partial correlations $\{-\hat{\omega}_{ij}/\sqrt{\hat{\omega}_{ii}\hat{\omega}_{jj}},i\ne j\}$
from the precision matrix estimator $\widehat{\mb{\Omega}}\coloneqq(\hat{\omega}_{ij})_{p\times p}$. Both CLIME and SPICE yield similar result, hence we only report the result of CLIME. We find that 2.71\% of all the pairs of nodes are connected conditional on all other nodes. Most of the nonzero estimated partial correlations have small absolute values, with the medium at 0.01 and the maximum at 0.45.
About 0.62\% of all the pairs of nodes are connected both marginally and directly.

Define the degree of a node to be the number of its connected nodes, and a hub to be a high-degree node. 
The marginal connectivity node degrees range from 0 to 164 with the medium at 2, and the direct connectivity node degrees range from 5 to 85 with the medium at 22.
The top 10 hubs found by either method are provided in Supplementary Material with 
six overlapping hubs. Seven of the top 10 hubs of marginal connectivity are spatially close to those in \cite{Buck09} and \cite{Cole10} obtained from multiple subjects. Note that they arbitrarily used 0.25 as the threshold value for the sample correlations, whereas our threshold value of 0.12 is selected from cross-validation. Some additional results are provided in Supplemental Material.

\section*{Acknowledgements}
The rfMRI data were provided by the Human Connectome Project, WU-Minn Consortium (Principal Investigators: David Van Essen and Kamil Ugurbil; 1U54MH091657), which is funded by the 16 NIH Institutes and Centers that support the NIH Blueprint for Neuroscience Research, and by the McDonnell Center for Systems Neuroscience at Washington University. We are grateful to the Editor Professor Edward George, the Associate Editor, and three anonymous referees for their constructive comments that greatly reshaped this article.

\renewcommand\theequation{A.\arabic{equation}}
\setcounter{equation}{0}

\renewcommand\thesection{A.\arabic{section}}
\setcounter{section}{0}

\renewcommand\thelem{A.\arabic{lem}}
\setcounter{lem}{0}

\section*{Appendix: Technical lemmas}\label{Appendix}
The keys to the proofs of Theorems \ref{thm1}-\ref{poly thm} are the proper concentration inequalities for $|\widehat{\mb{\Sigma}}-\mb{\Sigma} |_{\infty}$ and ${|\widehat{\mb{R}}-\mb{R} |_{\infty}}$ under temporal dependence. Once these inequalities are established, the rest of the proofs are straightforward extensions of those in \cite{Bick08a,CLL11,Ravi11,Roth08,Roth09}. We provide these inequalities in the following lemmas, where Part (i) in Lemma \ref{lem: ineq infty} is an extension of the Hoeffding-type inequality \cite[Theorem 7.27]{Fouc13}  and 
the Hanson-Wright inequality \cite[Theorem 1.1]{Rude13} 
from finite-dimensional to infinite-dimensional sub-Gaussian random vectors.
These lemmas can also be applied to the estimation of large band matrix \cite{Bick08b} and other high-dimensional time series problems such as  linear regression \cite{Wu16} and linear functionals \cite{Chen16}.

\begin{lem}\label{lem: ineq infty}	
Let $\bd{e}=(e_1,e_2,\dots)^\top$ be an infinite-dimensional random vector with 
each entry $e_i$ satisfying $E(e_i)=0$ and $E(e_i^2)=1$. Let $\bd{X}=\mb{A}\bd{e}$ and $\bd{Y}=\mb{B}\bd{e}$
be two well-defined random vectors with length $n$ in the sense of 
entrywise almost-sure convergence and mean-square convergence, 
where $\mb{A}$ and $\mb{B}$ are two deterministic matrices. 
For any $n$-dimensional deterministic vector $\bd{b}$ and all $u>0$, 

\smallskip

\textup{(i)}	if all $e_i$ satisfy condition~\ref{C1} with the same $K$, then
\be
\label{P_bX for subgauss}
P[|\bd{b}^\top\bd{X}|\ge u]\le 
	2\exp\left\{
	-\frac{Cu^2}{K^2\|\bd{b}\|_F^2\| \mb{A}\mb{A}^\top \|_2}
	\right\}
\ee
and
\begin{align}\label{P_XY for subgauss}
	P\left[
	\left|
	\bd{X}^\top\bd{Y}-E(\bd{X}^\top\bd{Y})
	\right|\ge u
	\right]\le
	2\exp\Bigg\{ &
	-C\min\bigg(
	\frac{u^2}{K^4\|\mb{A}\mb{A}^\top\|_F\|\mb{B}\mb{B}^\top\|_F},\nonumber\\
	&\frac{u}{K^2\sqrt{\|\mb{A}\mb{A}^\top\|_2\|\mb{B}\mb{B}^\top\|_2}}
	\bigg)
	\Bigg\}
\end{align}	
with an absolute constant $C>0$;

\smallskip

\textup{(ii)} if all $e_i$ satisfy condition~\ref{C2} with the same $K$ and $\vartheta$, then
\be\label{P_bX for subexp}
P[|\bd{b}^\top\bd{X}|\ge u]\le 
2\exp\left\{
-\frac{(u\|\bd{b} \|_F^{-1}\|\mb{AA}^\top  \|_2^{-1/2})^{\frac{1}{1/2+1/\vartheta}}}{CK(2/\vartheta)^{1/\vartheta}}
\right\}
\ee
and
\begin{align}\label{P_XY for subexp}
P\left[
|\bd{X}^\top\bd{Y}-E(\bd{X}^\top\bd{Y})|\ge u
\right] &\le
2\exp\left\{
-
\frac{
	(u\| \mb{AA}^\top\|_F^{-1/2}\| \mb{BB}^\top\|_F^{-1/2})^{\frac{1}{1/2+2/\vartheta}}
}{
C  K^2 (4/\vartheta)^{4/\vartheta}
}
\right\}
\nonumber\\
&\quad+
2\exp\left\{
-\frac{(u \| \mb{AA}^\top\|_F^{-1/2}\| \mb{BB}^\top\|_F^{-1/2}   )^{\frac{1}{1+2/\vartheta}}}{C K^2(2/\vartheta)^{2/\vartheta}}
\right\}
\end{align}
with an absolute constant $C>0$;

\smallskip

\textup{(iii)} if all $e_i$ satisfy condition ~\ref{C3} with the same $k$ and $\eta_k$, then
\be\label{P_bX for poly}
P[|\bd{b}^\top\bd{X}|\ge u]\le 
(C\eta_{k}k^{1/2}/u)^{k} \|\mb{A}\mb{A}^\top  \|_2^{k/2}  \| \bd{b}\|_F^{k}
\ee
and
\begin{align}\label{P_XY for subexp}
P\left[
|\bd{X}^\top\bd{Y}-E(\bd{X}^\top\bd{Y})|\ge u
\right]
&
\le (C \eta_{k}^2 k^{1/2}/u  )^{k/2}  \|  \mb{AA}^\top \|_F^{k/4}  \| \mb{BB}^\top \|_F^{k/4}, \nonumber\\
&\quad \ +
(C\eta_{k}^2 k/u)^{k}  \|  \mb{AA}^\top\|_F^{k/2}  \|\mb{BB}^\top  \|_F^{k/2}
\end{align}
with an absolute constant $C>0$.
\end{lem}

\begin{lem}\label{lem: max diff in sample sigma}
Let $v_0>0$ be an absolute constant. Suppose that $\mb{X}_{p\times n}$ is generated from \eqref{linear form}. Uniformly on $\mb{\Sigma}$ satisfying $|\mb{\Sigma}|_{\infty}\le v_0$ and
$\{\mb{R}_{[k]}\}_{k=1}^p$ subject to \eqref{g_F,g_2},
for any absolute constant $C>0$, if any of the following three conditions holds: 
\begin{enumerate}
\labitem{\textup{(i)}}{part(I)}
all $e_i$ satisfy condition~\ref{C1} with the same $K$, $u^*=C_1u_1$ with $C_1>0$ being a sufficiently large constant depending only on $(v_0, K, C)$, and 
$u_1$ is given in~\eqref{u for subgauss};

\labitem{\textup{(ii)}}{part(II)}
all ${e_i}$ satisfy condition~\ref{C2} with the same $K$ and $\vartheta$, $u^*=C_2u_2$ with $C_2>0$
being a sufficiently large constant depending only on $(v_0, K, \vartheta,C)$, and $u_2$ is given in~\eqref{u for subexp};

\labitem{\textup{(iii)}}{part(III)}
all $e_i$ satisfy condition~\ref{C3} with the same $k$ and $\eta_k$, $u^*=C_3u_3$ with $C_3>0$ being a sufficiently large constant depending only on $(v_0, k,\eta_k)$, 
and $u_3$ is given in~\eqref{u for poly} with the same constant $C$ given above;

\end{enumerate}
then we have
\[
P\left[
		|\widehat{\mb{\Sigma}}-\mb{\Sigma}|_{\infty}\ge u^*
		\right]=O(p^{-C}).
\]

\end{lem}

\begin{lem}\label{lem: max diff in sample R}
Suppose that $\mb{X}_{p\times n}$ is generated from \eqref{linear form}. 
Uniformly on $\{\mb{R}_{[k]}\}_{k=1}^p$ subject to \eqref{g_F,g_2},
for any absolute constant $C>0$, if any of the three conditions given in Lemma~\ref{lem: max diff in sample sigma} holds and the corresponding $u_j=o(1)$, $j\in\{1,2,3\}$, then we have
\[
P\left[
		|\widehat{\mb{R}}-\mb{R}|_{\infty}\ge u^*
		\right]=O(p^{-C}).
\]

\end{lem}

\renewcommand\theequation{S.\arabic{equation}}
\setcounter{equation}{0}

\renewcommand\thesection{S.\arabic{section}}
\setcounter{section}{0}
 \renewcommand{\thetable}{S\arabic{table}}
 \setcounter{table}{0}
 \renewcommand{\thefigure}{S\arabic{figure}}
 \setcounter{figure}{0}

\vspace{0.25cm}
\section*{Supplementary material}

\section{Technical Preparations}

\begin{prop}\label{subexp prop}
	Let $Z$ be a random variable, and $\vartheta>0$. The following three properties are equivalent with positive parameters $\{K_1,K_2,K_3\}$ satisfying $K_i\le C K_j$ for an absolute constant $C>0$:
	\begin{enumerate}[1.]
		\item\label{property1} $P(|Z|\ge u)\le 2\exp(-u^\vartheta/K_1)$ for all $u\ge0$;
		
		\item\label{property2} $(E|Z|^k)^{1/k}\le K_2 (k/\vartheta)^{1/\vartheta}$ for all $k\ge \vartheta$;
		
		\item\label{property3} $E(\exp\{|Z|^\vartheta/K_3\})\le 2$.	
	\end{enumerate}	
\end{prop}

\begin{proof}
	Let $X=|Z|^\vartheta$. Then applying Proposition~1.3.1 of \cite{Vers16} we obtain the results.
\end{proof}

In the following, 
we denote 
the diagonal and the off-diagonal of a square matrix $\mb{A}$ by
$\diag(\mb{A})$
and $\offdiag(\mb{A})$,
respectively. 
The proposition below is from Lemmas~B.2 and B.4 of \cite{Erdo13}. 
\begin{prop}\label{Erdo13 Lemma}
	Assume $\bd{Z}=(Z_1,\dots, Z_m)^\top$, $m<\infty$, where  $\{Z_i\}_{i=1}^m$ are independent mean-zero random variables.
	Let $(E|Z_i|^k)^{1/k}\le \eta_k$ for all $i=1,\dots,m$ and all $k>0$, where $\eta_k$ can be $\infty$. Then for 
	any deterministic vector $\bd{b}=(b_1,\dots,b_m)^\top$
	and all $k\ge 2$,
	\be\label{|bZ|_k}
	(E|\bd{b}^\top\bd{Z}|^k)^{1/k}\le C_1k^{1/2}\eta_k \|\bd{b} \|_F
	\ee
	with an absolute constant $C_1>0$.
	If additionally assume $E(Z_i^2)=1$, $i=1,\dots,m$,
	then for any deterministic matrix $\mb{A}=(a_{ij})_{m\times m}$ and all $k\ge 2$,
	\be\label{|aZi^2-1|_p}
	\left[E\left|\sum_{i=1}^m  a_{ii}(Z_i^2-1)
	\right|^k\right]^{1/k}\le C_2 k^{1/2} \eta_{2k}^2
	\|\diag(\mb{A}) \|_F
	\ee
	and
	\be\label{|aZiZj|_p}
	\left[E\left|
	\sum_{1\le i\ne j\le m} a_{ij}Z_iZ_j
	\right|^k\right]^{1/k}
	\le C_2 k \eta_k^2 \| \offdiag(\mb{A}) \|_F
	\ee
	with an absolute constant $C_2>0$.
\end{prop}
\begin{proof}
	When  $E(Z_i^2)=1$ for all $i$, applying Lemma~B.2 of \cite{Erdo13} gives \eqref{|bZ|_k}, but we found their proof applies without this additional restriction.
	
	Now assume $E(Z_i^2-1)=0$ for all $i$.
	For all $k\ge 1$, by Lyapounov's inequality,
	see (5.37) in \cite{Bill95},
	we have $\eta_{2k}\ge\eta_2\ge E(Z_i^2)=1$. Then by Minkowski's inequality, see (19.4) in \cite{Bill95}, we obtain
	\be\label{2mu_2k^2}
	(E|Z_i^2-1|^k)^{1/k}
	\le (E|Z_i^2|^k)^{1/k}+1
	=\eta_{2k}^2+1\le 2\eta_{2k}^2.
	\ee
	Thus, we obtain \eqref{|aZi^2-1|_p} from \eqref{|bZ|_k}. Furthermore, applying  Lemma~B.4 of \cite{Erdo13} gives \eqref{|aZiZj|_p}.
\end{proof}

\begin{prop}\label{lem: P_bZ}
	Let $\bd{Z}=(Z_1,\dots, Z_m)^\top$, $m<\infty$, where $\{Z_i\}_{i=1}^m$ are independent mean-zero random variables with the generalized sub-exponential tails defined in \ref{C2} in the main text with the same constants $K$ and $\vartheta$. Then for 
	any deterministic vector $\bd{b}=(b_1,\dots,b_m)^\top$
	and all $u>0$, we have
	\be\label{P_bZ}
	P\left[
	\left|
	\bd{b}^\top \bd{Z}
	\right|\ge u
	\right]
	\le 
	2\exp\left\{
	-\frac{(u/\|\bd{b} \|_F)^{\frac{1}{1/2+1/\vartheta}}}{C_1K(2/\vartheta)^{1/\vartheta}}
	\right\}
	\ee
	with an absolute constant $C_1>0$.
	If additionally assume $E(Z_i^2)=1$, $i=1,\dots,m$,
	then for any deterministic matrix $\mb{A}=(a_{ij})_{m\times m}$ and all $u> 0$, 
	\begin{align}\label{P_ZAZ}
		P\left[
		|\bd{Z}^\top\mb{A}\bd{Z}-E(\bd{Z}^\top\mb{A}\bd{Z})|\ge u
		\right] 
		&\le 
		2\exp\left\{
		-
		\frac{
			(u/ \| \diag(\mb{A}) \|_F)^{\frac{1}{1/2+2/\vartheta}}
		}{
		C_2  K^2 (4/\vartheta)^{4/\vartheta}
	}
	\right\}
	\nonumber\\
	&\qquad+
	2\exp\left\{
	-\frac{\left(u/\left\|\offdiag(\mb{A}) \right\|_F\right)^{\frac{1}{1+2/\vartheta}}}{C_2K^2(2/\vartheta)^{2/\vartheta}}
	\right\}
\end{align}
with an absolute constant $C_2>0$.
\end{prop}
\begin{proof}
	We only consider the non-trivial case for
	\eqref{P_bZ} with $\bd{b}\ne \bd{0}$.
	By \eqref{|bZ|_k} in Proposition~\ref{Erdo13 Lemma} and condition~\ref{C2} in the main text, for all $k\ge 2$, we have
	\[
	(E|\bd{b}^\top\bd{Z}|^k)^{1/k}
	\le c_1k^{1/2} K(k/\vartheta)^{1/\vartheta}\|\bd{b} \|_F
	=c_1K(1/\vartheta)^{1/\vartheta}k^{1/2+1/\vartheta}\| \bd{b} \|_F
	\]
	with an absolute constant $c_1>0$.
	Rewriting the above inequality yields
	\be\label{E_bZ}
	\left(
	E\left|
	\bd{b}^\top\bd{Z}/\|\bd{b}  \|_F
	\right|^k
	\right)^{1/k}
	\le
	c_1K(1/\vartheta)^{1/\vartheta}k^{1/2+1/\vartheta}.
	\ee
	By Lyapounov's inequality, see (5.37) in \cite{Bill95}, for all $k\in [1/(1/2+1/\vartheta),2]$,
	\begin{align*}
		\left(
		E\left|
		\bd{b}^\top\bd{Z}/\|\bd{b}  \|_F
		\right|^k
		\right)^{1/k}
		&\le
		\left(
		E\left|
		\bd{b}^\top\bd{Z}/\|\bd{b}  \|_F
		\right|^2
		\right)^{1/2}
		\le c_1K(1/\vartheta)^{1/\vartheta}2^{1/2+1/\vartheta}\\
		&\le
		c_1K(1/\vartheta)^{1/\vartheta}2^{1/2+1/\vartheta}
		[(1/2+1/\vartheta)k]^{1/2+1/\vartheta}.
	\end{align*}
	Together with~\eqref{E_bZ}, for all $k\ge 1/(1/2+1/\vartheta)$, we have
	\be\label{E_bZ/b}
	\left(
	E\left|
	\bd{b}^\top\bd{Z}/\|\bd{b}  \|_F
	\right|^k
	\right)^{1/k}
	\le
	c_12^{1/2}K(2/\vartheta)^{1/\vartheta}[(1/2+1/\vartheta)k]^{1/2+1/\vartheta}.
	\ee
	Then from 
	the equivalence between properties \ref{property1} and \ref{property2} in
	Proposition~\ref{subexp prop}, we have that for all $u>0$,
	\[
	P\left[
	\left|\bd{b}^\top\bd{Z}/\|\bd{b} \|_F\right|\ge u
	\right]
	\le 2\exp\left\{
	-\frac{u^{\frac{1}{1/2+1/\vartheta}}}{C_1 K(2/\vartheta)^{1/\vartheta}}
	\right\}
	\]
	with an absolute constant $C_1>0$, which yields~\eqref{P_bZ}.
	
	Furthermore, assume $E(Z_i^2)=1$ for all $i=1,\dots,m$.
	Note that for all $u>0$, we have
	\begin{align}\label{P_ZAZ terms}
		&P\left[
		|\bd{Z}^\top\mb{A}\bd{Z}-E(\bd{Z}^\top\mb{A}\bd{Z})|\ge u
		\right]\\
		&\ =
		P\left[
		\left|
		\sum_{1\le i,j\le m} a_{ij}Z_iZ_j
		-\sum_{1\le i,j\le m} a_{ij}E(Z_iZ_j)
		\right|\ge u
		\right]\nonumber\\
		&\ \le 
		P\left[
		\left|\sum_{i=1}^m a_{ii}(Z_i^2-1)
		\right|\ge u/2
		\right]
		+
		P\left[
		\left|
		\sum_{1\le i\ne j\le m} a_{ij}Z_iZ_j
		\right|\ge u/2
		\right]. \nonumber
	\end{align}
	Now consider the first term on the right hand side of~\eqref{P_ZAZ terms}. 
	Since $E(Z_i^2-1)=0$ for all $i$,
	by \eqref{2mu_2k^2} and condition~\ref{C2}, for all $k\ge 1$, we obtain
	\[
	(E|Z_i^2-1|^k)^{1/k}
	\le 
	2K^2(2k/\vartheta)^{2/\vartheta}.
	\]
	Then from Lyapounov's inequality, for all $k\in [\vartheta/2,1]$, we have
	\[
	(E|Z_i^2-1|^k)^{1/k}
	\le E|Z_i^2-1|
	\le 2K^2(2/\vartheta)^{2/\vartheta}
	\le 2K^2(2/\vartheta)^{2/\vartheta}(2k/\vartheta)^{2/\vartheta}.
	\]
	Hence, for all $k\ge \vartheta/2$,
	\[
	(E|Z_i^2-1|^k)^{1/k}\le 
	2K^2(2/\vartheta)^{2/\vartheta}(2k/\vartheta)^{2/\vartheta}.
	\]
	Thus, by inequality~\eqref{P_bZ},
	we obtain
	\begin{align}\label{P_ZAZ term 2}
		P\left[
		\left|\sum_{i=1}^m a_{ii}(Z_i^2-1)
		\right|\ge u/2
		\right]
		&\le
		2\exp\left\{
		-
		\frac{
			\left[u/\left(2 \| \diag(\mb{A}) \|_F\right)\right]^{\frac{1}{1/2+2/\vartheta}}
		}{
		c_2 [ K^2(2/\vartheta)^{2/\vartheta}]  (4/\vartheta)^{2/\vartheta}
	}
	\right\}\nonumber\\
	&\le 
	2\exp\left\{
	-
	\frac{
		\left[u/ \left(2\| \diag(\mb{A}) \|_F\right)\right]^{\frac{1}{1/2+2/\vartheta}}
	}{
	c_2  K^2 (4/\vartheta)^{4/\vartheta}
}
\right\}\nonumber\\
&\le 2\exp\left\{
-
\frac{
	(u/ \| \diag(\mb{A}) \|_F)^{\frac{1}{1/2+2/\vartheta}}
}{
c_3  K^2 (4/\vartheta)^{4/\vartheta}
}
\right\}
\end{align}
for all $u>0$, where 
$c_2$ is an absolute constant and
$c_3=2^{2/3}c_2\ge 2^{\frac{1}{1/2+2/\vartheta}}c_2$.

For the second term
on the right hand side of~\eqref{P_ZAZ terms}, we only consider the non-trivial case that 
$\offdiag(\mb{A})\ne \bd{0}$.
By inequality~\eqref{|aZiZj|_p} and condition~\ref{C2},
for all $k\ge 2$, we have
\begin{align*}
	\left[E\Big|
	\sum_{1\le i\ne j\le m} a_{ij}Z_iZ_j
	\Big|^k\right]^{1/k}
	&\le c_4 k K^2 (k/\vartheta)^{2/\vartheta}\|\offdiag(\mb{A}) \|_F\\
	&=c_4  K^2 (1/\vartheta)^{2/\vartheta}
	k^{1+2/\vartheta}
	\|\offdiag(\mb{A}) \|_F
\end{align*}
with an absolute constant $c_4>0$.
Then by Lyapounov's inequality and the same approach used to obtain~\eqref{E_bZ/b}, for all $k\ge 1/(1+2/\vartheta)$,
\begin{align*}
	\left[E\Big|
	\sum_{1\le i\ne j\le m} \frac{a_{ij}Z_iZ_j}{\left\|\offdiag(\mb{A}) \right\|_F}
	\Big|^k\right]^{1/k}
	&\le 
	c_4  K^2 (1/\vartheta)^{2/\vartheta}
	2^{1+2/\vartheta}
	[(1+2/\vartheta)k]^{1+2/\vartheta}\\
	&=2c_4K^2(2/\vartheta)^{2/\vartheta}
	[(1+2/\vartheta)k]^{1+2/\vartheta}.
\end{align*}
Hence, from the equivalence between properties \ref{property1} and \ref{property2} in
Proposition~\ref{subexp prop}, we have that for all $u>0$,
\[
P\left[
\left|
\sum_{1\le i\ne j\le m} \frac{a_{ij}Z_iZ_j}{\left\|\offdiag(\mb{A}) \right\|_F}
\right|\ge u
\right]
\le 2\exp\left\{
-\frac{u^{\frac{1}{1+2/\vartheta}}}{c_5K^2(2/\vartheta)^{2/\vartheta}}
\right\}
\]
with an absolute constant $c_5>0$. Replacing $u$ by
$u/(2\left\|\offdiag(\mb{A}) \right\|_F)$
in the above inequality yields 
\begin{align}\label{P_ZAZ term 1}
	P\left[
	\left|
	\sum_{1\le i\ne j\le m} a_{ij}Z_iZ_j
	\right|\ge u/2
	\right]
	&\le
	2\exp\left\{
	-\frac{\left[u/(2\left\|\offdiag(\mb{A}) \right\|_F)\right]^{\frac{1}{1+2/\vartheta}}}{c_5K^2(2/\vartheta)^{2/\vartheta}}
	\right\}\nonumber\\
	&\le 2\exp\left\{
	-\frac{\left(u/\left\|\offdiag(\mb{A}) \right\|_F\right)^{\frac{1}{1+2/\vartheta}}}{c_6K^2(2/\vartheta)^{2/\vartheta}}
	\right\}
\end{align}
for all $u>0$, where $c_6=c_5\sqrt{2}\ge  2^{1/(1+2/\vartheta)}c_5 $.
Plugging~\eqref{P_ZAZ term 1} and~\eqref{P_ZAZ term 2} into~\eqref{P_ZAZ terms} yields~\eqref{P_ZAZ} with $C_2=\max(c_3,c_6)$.
\end{proof}

Now denote $\mb{P}_i$ to be the $n\times pn$ 
matrix with $\sigma_{ii}^{-1/2}$ in the $\big(j,i+(j-1)p\big)$ entries and $0$ in all other entries, $j=1, \dots, n$, and	
define $\widetilde{\mb{R}}=(\tilde{\rho}_{ij})_{p\times p}=(\hat{\sigma}_{ij}/\sqrt{\sigma_{ii}\sigma_{jj}})_{p\times p}$.
Note that $\widetilde{\mb{R}}$ differs from the sample correlation matrix
$\widehat{\mb{R}}=(\hat{\sigma}_{ij}/\sqrt{\hat{\sigma}_{ii}\hat{\sigma}_{jj}})_{p\times p}$. 

\begin{prop}\label{lem: max diff}
	Suppose $\mb{X}_{p\times n}$ is generated from \eqref{linear form}. For all $u>0$, we have
	\begin{eqnarray}	
		\label{max difference in sigmas}	
		&& \hspace{-0.2in} P\left[
		|\widehat{\mb{\Sigma}}-\mb{\Sigma}|_{\infty}\ge 2u
		\right] \nonumber\\	
		&& \quad \le
		\sum_{i=1}^p
		P\left[|\bd{1}_{n}^\top\mb{P}_i\mb{H}\bd{e}|
		\ge n\sqrt{\frac{u}{\sigma_{ii}}}
		\right]
		\\
		&& \qquad   +
		\sum_{1\le i,j\le p}P\left[
		\left|
		\bd{e}^\top\mb{H}^\top\mb{P}_i^\top
		\mb{P}_j\mb{H}\bd{e}
		-
		E(\bd{e}^\top\mb{H}^\top\mb{P}_i^\top
		\mb{P}_j\mb{H}\bd{e})
		\right|\ge \frac{nu}{\sqrt{\sigma_{ii}\sigma_{jj}}}
		\right] \nonumber
	\end{eqnarray}	
	and 	
	\begin{eqnarray}
		\label{max difference in pseudo-rhos}	
		&&	\hspace{-0.2in} P\left[
		|\widetilde{\mb{R}}-\mb{R}|_{\infty}\ge 2u
		\right]\nonumber\\	
		&& \quad \le
		\sum_{i=1}^p
		P\left[|\bd{1}_{n}^\top\mb{P}_i\mb{H}\bd{e}|
		\ge nu^{1/2}
		\right]
		\\
		&& \qquad+
		\sum_{1\le i,j\le p}P\left[
		\left|
		\bd{e}^\top\mb{H}^\top\mb{P}_i^\top
		\mb{P}_j\mb{H}\bd{e}
		-
		E(\bd{e}^\top\mb{H}^\top\mb{P}_i^\top
		\mb{P}_j\mb{H}\bd{e})
		\right|\ge nu
		\right].\nonumber
	\end{eqnarray}	
	Furthermore, if
	\be\label{general scale of max diff pseudo-rhos}
	P\left[
	|\widetilde{\mb{R}}-\mb{R}|_{\infty}\ge 2u
	\right]=O(p^{-C})
	\ee
	with an absolute constant $C>0$ and $0<u=o(1)$, then we have
	\be\label{general scale of max diff rhos}
	P\left[
	|\widehat{\mb{R}}-\mb{R}|_{\infty}\ge 8u
	\right]=O(p^{-C}).
	\ee	
\end{prop}

\begin{proof}
	Without loss of generality, we assume $\bd{\mu}_p=\bd{0}$. We only consider $\mb{X}_{p\times n}$ generated from \eqref{linear form} with $m=\infty$ because any case with finite $m$ can be constructed by adding infinite number of zero columns in $\mb{H}$.
	Define
	$\bd{Z}_{[i]}=(Z_{i1},\dots,Z_{in})^\top$ with
	$Z_{ij}=X_{ij}/\sqrt{\sigma_{ii}}$,
	then by
	\eqref{linear form} we have $\bd{Z}_{[i]}=\mb{P}_i\mb{H}\bd{e}$.
	Since
	\begin{align*}
		\max_{1\le i,j\le p}|\hat{\sigma}_{ij}-\sigma_{ij}|
		&\le
		\max_{1\le i,j\le p}\left| \bar{X}_i\bar{X}_j\right|
		+
		\max_{1\le i,j\le p}\left|\frac{1}{n}\sum_{k=1}^nX_{ik}X_{jk}-\sigma_{ij}\right|
		\nonumber\\
		&=
		\max_{1\le i\le p}\left| \bar{X}_i  \right|^2
		+
		\max_{1\le i,j\le p}\left|\frac{1}{n}\sum_{k=1}^nX_{ik}X_{jk}-\sigma_{ij}\right|,
	\end{align*}
	where $\bar{X}_i=\frac{1}{n}\sum_{j=1}^n X_{ij}$,
	for any $u>0$, we obtain \eqref{max difference in sigmas} in the following:
	\begin{eqnarray*}
		&& P\left[
		\max_{1\le i,j\le p}|\hat{\sigma}_{ij}-\sigma_{ij}|\ge 2u
		\right]\nonumber\\
		&& \qquad \le	
		P\left[
		\max_{1\le i\le p}\left| \bar{X}_i  \right|^2\ge u
		\right]
		+
		P\left[
		\max_{1\le i,j\le p}\left|\frac{1}{n}\sum_{k=1}^nX_{ik}X_{jk}-\sigma_{ij}\right|\ge u
		\right]
		\nonumber\\
		&& \qquad \le
		\sum_{i=1}^pP\left[
		\left| \frac{1}{ n} \sum_{j=1}^{n} X_{ij}\right|
		\ge u^{1/2}
		\right] \\
		&& \qquad \qquad +
		\sum_{1\le i,j\le p}P\left[
		\left|\frac{1}{ n}\sum_{k=1}^{n}X_{ik}X_{jk}-\sigma_{ij}\right|\ge u
		\right]
		\nonumber\\
		&& \qquad =	
		\sum_{i=1}^pP\left[
		\left|\frac{1}{ n}\sum_{j=1}^{n} Z_{ij}\right|\ge  \sqrt{\frac{u}{\sigma_{ii}}}
		\right] \\
		&& \qquad \qquad	+
		\sum_{1\le i,j\le p}P\left[
		\left|
		\frac{1}{ n}\sum_{k=1}^{n}Z_{ik}Z_{jk}-\rho_{ij}
		\right|\ge \frac{u}{\sqrt{\sigma_{ii}\sigma_{jj}}}
		\right]
		\nonumber\\
		&& \qquad =
		\sum_{i=1}^pP\left[
		\left|\bd{1}_{n}^\top\bd{Z}_{[i]}\right|\ge  n\sqrt{\frac{u}{\sigma_{ii}}}
		\right] \\
		&& \qquad \qquad 	+
		\sum_{1\le i,j\le p}P\left[
		\left|
		\bd{Z}_{[i]}^\top\bd{Z}_{[j]}-E(\bd{Z}_{[i]}^\top\bd{Z}_{[j]})
		\right|\ge \frac{nu}{\sqrt{\sigma_{ii}\sigma_{jj}}}
		\right]
		\nonumber\\
		&& \qquad =
		\sum_{i=1}^p
		P\left[|\bd{1}_{n}^\top\mb{P}_i\mb{H}\bd{e}|
		\ge n\sqrt{\frac{u}{\sigma_{ii}}}
		\right]
		\nonumber\\
		&& \qquad \qquad +
		\sum_{1\le i,j\le p}P\left[
		\left|
		\bd{e}^\top\mb{H}^\top\mb{P}_i^\top
		\mb{P}_j\mb{H}\bd{e}
		-
		E(\bd{e}^\top\mb{H}^\top\mb{P}_i^\top
		\mb{P}_j\mb{H}\bd{e})
		\right|\ge \frac{nu}{\sqrt{\sigma_{ii}\sigma_{jj}}}
		\right].
	\end{eqnarray*}
	
	Note that
	\begin{align*}
		\max_{1\le i,j\le p} |\tilde{\rho}_{ij}-\rho_{ij}|&=\max_{1\le i,j\le p}\left|\frac{\hat{\sigma}_{ij}-\sigma_{ij}}{\sqrt{\sigma_{ii}\sigma_{jj}}}\right|
		\\
		&\le
		\max_{1\le i,j\le p}\left|\frac{\bar{X}_i\bar{X}_j}{\sqrt{\sigma_{ii}\sigma_{jj}}}\right|
		+
		\max_{1\le i,j\le p}\left|\frac{1}{n\sqrt{\sigma_{ii}\sigma_{jj}}}\sum_{k=1}^nX_{ik}X_{jk}-\rho_{ij}\right|\nonumber\\
		&= \max_{1\le i\le p}\left|\frac{\bar{X}_i}{\sqrt{\sigma_{ii}}}   \right|^2
		+
		\max_{1\le i,j\le p}\left|\frac{1}{n\sqrt{\sigma_{ii}\sigma_{jj}}}\sum_{k=1}^nX_{ik}X_{jk}-\rho_{ij}\right|\nonumber\\
		&=\max_{1\le i\le p}\left|\bar{Z}_i\right|^2
		+
		\max_{1\le i,j\le p}\left|  \frac{1}{n}\sum_{k=1}^nZ_{ik}Z_{jk}-\rho_{ij} \right|.
	\end{align*}
	Then following the lines
	of the proof of \eqref{max difference in sigmas},
	we obtain~\eqref{max difference in pseudo-rhos}.
	
	Now assume that $u=o(1)$ and \eqref{general scale of max diff pseudo-rhos} holds.
	We have
	\begin{align*}
		O(p^{-C})
		&=P\left[
		\max_{1\le i,j\le p}\left|\frac{\hat{\sigma}_{ij}-\sigma_{ij}}{\sqrt{\sigma_{ii}\sigma_{jj}}}\right|\ge 2u
		\right]\\
		&\ge
		P\left[
		\max_{1\le i\le p}\left|\frac{\hat{\sigma}_{ii}-\sigma_{ii}}{\sigma_{ii}}\right|\ge 2u
		\right]\nonumber\\
		&=P\left[
		\max_{1\le i\le p}\left(
		\left|\sqrt{\frac{\hat{\sigma}_{ii}}{\sigma_{ii}}}-1\right|\left|\sqrt{\frac{\hat{\sigma}_{ii}}{\sigma_{ii}}}+1\right|
		\right)
		\ge 2u
		\right]\nonumber\\
		&\ge P\left[
		\max_{1\le i\le p}\left|\sqrt{\frac{\hat{\sigma}_{ii}}{\sigma_{ii}}}-1\right|\ge 2u
		\right],
	\end{align*}
	and
	\begin{align*}
		O(p^{-C})
		&=
		P\left[
		\max_{1\le i\le p}\left|\frac{\hat{\sigma}_{ii}-\sigma_{ii}}{\sigma_{ii}}\right|\ge 2u
		\right]\\
		&=P\left[ \max_{1\le i\le p}\left| \frac{\hat{\sigma}_{ii}}{\sigma_{ii}}-1  \right|^2 \ge 4u^2\right]\nonumber\\
		&=
		P\left[
		\max_{1\le i,j\le p}
		\left(\left| \frac{\hat{\sigma}_{ii}}{\sigma_{ii}}-1  \right| \left| \frac{\hat{\sigma}_{jj}}{\sigma_{jj}}-1  \right|\right)
		\ge 4u^2
		\right]\nonumber\\
		&=
		P\left[
		\max_{1\le i,j\le p}
		\left(
		\left| \sqrt{\frac{\hat{\sigma}_{ii}}{\sigma_{ii}}}-1  \right| \left| \sqrt{\frac{\hat{\sigma}_{jj}}{\sigma_{jj}}}-1  \right|
		\left| \sqrt{\frac{\hat{\sigma}_{ii}}{\sigma_{ii}}}+1  \right| \left| \sqrt{\frac{\hat{\sigma}_{jj}}{\sigma_{jj}}}+1  \right| 
		\right)
		\ge 4u^2
		\right]\nonumber \\
		&\ge
		P\left[
		\max_{1\le i,j\le p}
		\left| \sqrt{\frac{\hat{\sigma}_{ii}\hat{\sigma}_{jj}}{\sigma_{ii}\sigma_{jj}}}-\sqrt{\frac{\hat{\sigma}_{ii}}{\sigma_{ii}}}-\sqrt{\frac{\hat{\sigma}_{jj}}{\sigma_{jj}}}+1  \right|
		\ge 4u^2
		\right]\nonumber \\
	\end{align*}
	\begin{align*}
		&\ge
		P\left[
		\max_{1\le i,j\le p}\left(
		\left| \sqrt{\frac{\hat{\sigma}_{ii}\hat{\sigma}_{jj}}{\sigma_{ii}\sigma_{jj}}}-1\right|
		-
		\left|\sqrt{\frac{\hat{\sigma}_{ii}}{\sigma_{ii}}}-1 \right|
		-
		\left|\sqrt{\frac{\hat{\sigma}_{jj}}{\sigma_{jj}}}-1 \right|
		\right)
		\ge 4u^2
		\right]\nonumber \\
		&\ge
		P\left[
		\max_{1\le i,j\le p}\left| \sqrt{\frac{\hat{\sigma}_{ii}\hat{\sigma}_{jj}}{\sigma_{ii}\sigma_{jj}}}-1\right|
		\ge 4u^2
		+2\max_{1\le i\le p}\left|\sqrt{\frac{\hat{\sigma}_{ii}}{\sigma_{ii}}}-1 \right|
		\right].
	\end{align*}
	Then by $u=o(1)$, we have
	\begin{align*}
		\lefteqn{P\left[
			\max_{1\le i,j\le p}	\left| \sqrt{\frac{\hat{\sigma}_{ii}\hat{\sigma}_{jj}}{\sigma_{ii}\sigma_{jj}}}-1\right|
			\ge 6u
			\right]}\nonumber\\
		&\qquad \le
		P\left[
		\max_{1\le i,j\le p}
		\left| \sqrt{\frac{\hat{\sigma}_{ii}\hat{\sigma}_{jj}}{\sigma_{ii}\sigma_{jj}}}-1\right|
		\ge 6u, 
		\max_{1\le i \le p}\left|\sqrt{\frac{\hat{\sigma}_{ii}}{\sigma_{ii}}}-1 \right|\le 2u
		\right]\nonumber\\
		&\qquad\qquad +P\left[
		\max_{1\le i \le p}\left|\sqrt{\frac{\hat{\sigma}_{ii}}{\sigma_{ii}}}-1 \right|\ge 2u \right]\nonumber\\
		&\qquad \le
		P\left[
		\max_{1\le i,j\le p} \left| \sqrt{\frac{\hat{\sigma}_{ii}\hat{\sigma}_{jj}}{\sigma_{ii}\sigma_{jj}}}-1\right|
		\ge 2u+
		2\max_{1\le i\le p}\left|\sqrt{\frac{\hat{\sigma}_{ii}}{\sigma_{ii}}}-1 \right|\right]+O(p^{-C})\nonumber\\
		&\qquad \le
		P\left[
		\max_{1\le i,j\le p} \left| \sqrt{\frac{\hat{\sigma}_{ii}\hat{\sigma}_{jj}}{\sigma_{ii}\sigma_{jj}}}-1\right|
		\ge 4u^2+
		2\max_{1\le i\le p}\left|\sqrt{\frac{\hat{\sigma}_{ii}}{\sigma_{ii}}}-1 \right|\right]+O(p^{-C})\nonumber\\
		&\qquad =O(p^{-C}).
	\end{align*}
	Thus,
	\begin{align*}
		&P\left[
		\max_{1\le i,j\le p}\left|\hat{\rho}_{ij}-\rho_{ij}\right|\ge 8u
		\right]\nonumber\\
		&\qquad \le
		P\left[
		\max_{1\le i,j\le p}\left|
		\frac{\hat{\sigma}_{ij}}{\sqrt{\hat{\sigma}_{ii}\hat{\sigma}_{jj}}}-\frac{\hat{\sigma}_{ij}}{\sqrt{{\sigma}_{ii}{\sigma}_{jj}}}
		\right|+
		\max_{1\le i,j\le p}\left|\frac{\hat{\sigma}_{ij}-\sigma_{ij}}{\sqrt{\sigma_{ii}\sigma_{jj}}}\right|\ge 8u
		\right]\nonumber\\
		&\qquad \le
		P\left[
		\max_{1\le i,j\le p}\left(\left|
		\frac{\hat{\sigma}_{ij}}{\sqrt{\hat{\sigma}_{ii}\hat{\sigma}_{jj}}}
		\right|
		\left| \sqrt{\frac{\hat{\sigma}_{ii}\hat{\sigma}_{jj}}{\sigma_{ii}\sigma_{jj}}}-1\right|\right)
		\ge 6u\right] \\
		& \qquad \qquad +
		P\left[
		\max_{1\le i,j\le p}\left|\frac{\hat{\sigma}_{ij}-\sigma_{ij}}{\sqrt{\sigma_{ii}\sigma_{jj}}}\right|\ge 2u
		\right]\nonumber\\
		&\qquad \le
		P\left[
		\max_{1\le i,j\le p}\left| \sqrt{\frac{\hat{\sigma}_{ii}\hat{\sigma}_{jj}}{\sigma_{ii}\sigma_{jj}}}-1\right|\ge 6u\right]
		+
		P\left[
		\max_{1\le i,j\le p}\left|\frac{\hat{\sigma}_{ij}-\sigma_{ij}}{\sqrt{\sigma_{ii}\sigma_{jj}}}\right|\ge 2u
		\right]\nonumber\\
		&\qquad = O(p^{-C}).
	\end{align*}
\end{proof}

\section{Proofs of Technical Lemmas}	
\begin{proof}[Proof of Lemma~\ref{lem: ineq infty}]	
	First consider part (i). Let $\mb{A}=(a_{ij})_{n\times \infty}$ and 
	$\mb{B}=(b_{ij})_{n\times \infty}$. Let 
	$\mb{A}_{m}=(a_{ij})_{n\times m}$ 
	and $\mb{B}_{m}=(b_{ij})_{n\times m}$
	be the first $m$ columns of $\mb{A}$
	and $\mb{B}$, respectively. Let 
	$\bd{e}_{m}=(e_1,e_2,...,e_{m})^\top$ be the first $m$ elements of $\bd{e}$,
	$\bd{X}_{m}=(X_1^{m},...,X_n^{m})^\top=\mb{A}_{m}\bd{e}_{m}$, and
	$\bd{Y}_{m}=(Y_1^{m},...,Y_n^{m})^\top=\mb{B}_{m}\bd{e}_{m}$.
	By the entrywise almost-sure convergence and mean-square convergence,
	for each $i$, when $m\to \infty$, we have 
	$X_i^{m}=\sum_{j=1}^{m}a_{ij}e_j\stackrel{P}{\to} X_i=\sum_{j=1}^{\infty}a_{ij}e_j$,
	$Y_i^{m}=\sum_{j=1}^{m}b_{ij}e_j\stackrel{P}{\to} Y_i=\sum_{j=1}^{\infty}b_{ij}e_j$,
	$\sum_{j=1}^{\infty} a_{ij}^2<\infty$ and $\sum_{j=1}^{\infty} b_{ij}^2<\infty$.
	Thus, for any positive $n$, $\varepsilon_1$, $\varepsilon_2$, and $\delta$, there exists a number $N$ such
	that for any $m>N$, we have
	\be\label{XY-XY}
	P\left[|\bd{X}^\top\bd{Y}-\bd{X}_{m}^\top\bd{Y}_{m}|\ge \varepsilon_1\right]\le \delta,
	\ee
	\be\label{bX-bX}
	P\left[|\bd{b}^\top\bd{X}-\bd{b}^\top\bd{X}_{m}|\ge \varepsilon_1\right]\le \delta,
	\ee
	and for each $1\le i,j\le n$,
	\be\label{ab-ab}
	\left|\sum_{k=1}^ma_{ik}b_{ik}-\sum_{k=1}^{\infty}a_{ik}b_{ik}\right|\le  \varepsilon_2/n,
	\ee
	\be\label{aa-aa}
	\left| \sum_{k=1}^{m}a_{ik}a_{jk}- \sum_{k=1}^{\infty}a_{ik}a_{jk}\right|\le \delta/n,
	\ee
	\be\label{bb-bb}
	\left| \sum_{k=1}^{m}b_{ik}b_{jk}- \sum_{k=1}^{\infty}b_{ik}b_{jk}\right|\le \delta/n.
	\ee
	The convergence of  $\sum_{k=1}^{m}a_{ik}b_{ik}$ given in
	\eqref{ab-ab} holds because 
	$$\sum_{k=1}^\infty|a_{ik}b_{ik}|\le\sqrt{\sum_{k=1}^\infty a_{ik}^2\sum_{k=1}^\infty b_{ik}^2}<\infty,$$ 
	and similarly we have
	\eqref{aa-aa} and \eqref{bb-bb}.
	Thus,
	\begin{align}\label{AA_F}
		\|\mb{A}_m\mb{A}_m^\top\|_F
		&\le \|\mb{A}\mb{A}^\top\|_F+
		\|\mb{A}_m\mb{A}_m^\top-\mb{A}\mb{A}^\top\|_F\\
		&\le  \|\mb{A}\mb{A}^\top\|_F+\sqrt{\sum_{1\le i,j\le n}\left|\sum_{k=1}^ma_{ik}a_{jk}-
			\sum_{k=1}^{\infty}a_{ik}a_{jk}
			\right|^2}\nonumber\\
		&\le \|\mb{A}\mb{A}^\top\|_F+\sqrt{n^2(\delta/n)^2} \nonumber \\
		& =\|\mb{A}\mb{A}^\top\|_F+\delta,\nonumber
	\end{align}
	and
	\begin{align}\label{AA_2}
		\|\mb{A}_m\mb{A}_m^\top\|_2
		&\le \|\mb{A}\mb{A}^\top\|_2+
		\|\mb{A}_m\mb{A}_m^\top-\mb{A}\mb{A}^\top\|_2\\
		&\le
		\|\mb{A}\mb{A}^\top\|_2+
		\|\mb{A}_m\mb{A}_m^\top-\mb{A}\mb{A}^\top\|_F \nonumber \\
		&\le \|\mb{A}\mb{A}^\top\|_2+\delta.\nonumber
	\end{align}
	Similarly,	
	\be\label{BB_F2}
	\|\mb{B}_m\mb{B}_m^\top\|_F\le \|\mb{B}\mb{B}^\top\|_F+\delta
	\qquad\text{and}\qquad
	\|\mb{B}_m\mb{B}_m^\top\|_2\le \|\mb{B}\mb{B}^\top\|_2+\delta.
	\ee
	By Theorem 1.1 in \cite{Rude13},
	for all $u>0$, we have
	\begin{align*}
		\lefteqn{P\left[
			|\bd{X}_m^\top\bd{Y}_m-E(\bd{X}_m^\top\bd{Y}_m)|
			\ge u/4
			\right]}\nonumber\\
		& \qquad \le 
		2\exp\left\{
		-c_1\min\left(
		\frac{u^2}{K^4\|\mb{A}_m^\top\mb{B}_m\|_F^2},\frac{u}{K^2\|\mb{A}_m^\top\mb{B}_m\|_2}
		\right)
		\right\}
	\end{align*}
	with an absolute constant $c_1>0$.
	By Lemma 5.5 of \cite{Vers12} and Theorem 7.27 of \cite{Fouc13},	
	for all $u>0$, we have
	\begin{align*}
		P\left[
		|\bd{b}^\top\bd{X}_m|\ge u/2
		\right]\le 2\exp\left\{
		-\frac{c_2u^2}{K^2\|\bd{b}^\top\mb{A}_m\|_F^2}
		\right\}
	\end{align*}
	with an absolute constant $c_2>0$.	
	Since 
	\begin{align}\label{|AB|_F}
		\|\mb{A}_m^\top\mb{B}_m\|_F^2 &=\text{tr}(\mb{A}_m^\top\mb{B}_m\mb{B}_m^\top\mb{A}_m) \\
		&=\text{tr}(\mb{A}_m\mb{A}_m^\top\mb{B}_m\mb{B}_m^\top)\nonumber\\
		&\le \sqrt{\text{tr}(\mb{A}_m\mb{A}_m^\top\mb{A}_m\mb{A}_m^\top)\text{tr}(\mb{B}_m\mb{B}_m^\top\mb{B}_m\mb{B}_m^\top)}\nonumber\\
		&=\|\mb{A}_m\mb{A}_m^\top\|_F\|\mb{B}_m\mb{B}_m^\top\|_F, \nonumber \\
		\|\mb{A}_m^\top\mb{B}_m\|_2 & \le \|\mb{A}_m^\top\|_2\|\mb{B}_m\|_2 \nonumber \\
		&=\sqrt{\varphi_{\max}(\mb{A}_m\mb{A}_m^\top)\varphi_{\max}(\mb{B}_m^\top\mb{B}_m)}\nonumber \\
		&=\sqrt{\varphi_{\max}(\mb{A}_m\mb{A}_m^\top)\varphi_{\max}(\mb{B}_m\mb{B}_m^\top)}\nonumber \\
		&=\sqrt{\|\mb{A}_m\mb{A}_m^\top\|_2\|\mb{B}_m\mb{B}_m^\top\|_2}, \nonumber 
	\end{align}
	and
	\begin{align}\label{|bA|_F}
		\|\bd{b}^\top\mb{A}_m\|_F^2 &=\|\mb{A}_m^\top\bd{b}\|_F^2
		\le \|\mb{A}_m^\top\|_2^2\|\bd{b}\|_F^2\\
		&	= \varphi_{\max}(\mb{A}_m\mb{A}_m^\top)\|\bd{b}\|_F^2
		=\| \mb{A}_m\mb{A}_m^\top \|_2\|\bd{b}\|_F^2,\nonumber
	\end{align}
	where \eqref{|bA|_F} is obtained from Lemma 1 in \cite{Lam09}, we have
	\begin{align}\label{cctieq1}
		\lefteqn{P\left[
			\left|
			\bd{X}_m^\top\bd{Y}_m-E(\bd{X}_m^\top\bd{Y}_m)
			\right|\ge u/4
			\right]}\\
		&\qquad \le 2\exp\Bigg\{
		-c_1\min\bigg(
		\frac{u^2}{K^4\|\mb{A}_m\mb{A}_m^\top\|_F\|\mb{B}_m\mb{B}_m^\top\|_F},\nonumber\\
		&\qquad\qquad\qquad\qquad\quad
		\frac{u}{K^2\sqrt{\|\mb{A}_m\mb{A}_m^\top\|_2\|\mb{B}_m\mb{B}_m^\top\|_2}}
		\bigg)
		\Bigg\}\nonumber
	\end{align}
	and
	\begin{align}\label{cctieq2}
		P\left[
		|\bd{b}^\top\bd{X}_m|\ge u/2
		\right]\le 2\exp\left\{
		-\frac{c_2u^2}{K^2\|\bd{b}\|_F^2\| \mb{A}_m\mb{A}_m^\top \|_2}
		\right\}.
	\end{align}
	Let $\varepsilon_1=u/2$ and $\varepsilon_2=u/4$, then by 
	\eqref{XY-XY}, \eqref{ab-ab}, \eqref{cctieq1}, \eqref{AA_F}, \eqref{AA_2}, and \eqref{BB_F2} 
	we have
	\begin{align}\label{long1}
		\lefteqn{
			P\left[
			\left|
			\bd{X}^\top\bd{Y}-E(\bd{X}^\top\bd{Y})
			\right|\ge u
			\right]}\nonumber\\
		&\qquad\le
		P\left[
		|E(\bd{X}_m^\top\bd{Y}_m)-E(\bd{X}^\top\bd{Y})|
		+
		|\bd{X}_m^\top\bd{Y}_m-E(\bd{X}_m^\top\bd{Y}_m)|\ge u/2
		\right]\nonumber\\ 
		&\qquad\qquad+P\left[
		\left|
		\bd{X}^\top\bd{Y}-\bd{X}_m^\top\bd{Y}_m
		\right|\ge u/2\right]\nonumber\\
		&\qquad\le
		P\left[
		\left|\sum_{i=1}^n\sum_{k=1}^ma_{ik}b_{ik}-\sum_{i=1}^n\sum_{k=1}^{\infty}a_{ik}b_{ik}\right|
		+
		|\bd{X}_m^\top\bd{Y}_m-E(\bd{X}_m^\top\bd{Y}_m)|\ge u/2
		\right]\nonumber\\ 
		&\qquad\qquad+P\left[
		\left|
		\bd{X}^\top\bd{Y}-\bd{X}_m^\top\bd{Y}_m
		\right|\ge \varepsilon_1\right]\nonumber\\
		&\qquad\le
		P\left[
		\sum_{i=1}^n\left|\sum_{k=1}^ma_{ik}b_{ik}-\sum_{k=1}^{\infty}a_{ik}b_{ik}\right|
		+
		|\bd{X}_m^\top\bd{Y}_m-E(\bd{X}_m^\top\bd{Y}_m)|\ge u/2
		\right]+\delta\nonumber\\ 
		&\qquad\le 
		P\left[
		|\bd{X}_m^\top\bd{Y}_m-E(\bd{X}_m^\top\bd{Y}_m)|\ge u/2-\varepsilon_2
		\right]+\delta\nonumber\\
		&\qquad\le2\exp\Bigg\{
		-c_1\min\Bigg(
		\frac{u^2}{K^4\|\mb{A}_m\mb{A}_m^\top\|_F\|\mb{B}_m\mb{B}_m^\top\|_F},\nonumber\\
		&\qquad\qquad\qquad\qquad\qquad
		\frac{u}{K^2\sqrt{\|\mb{A}_m\mb{A}_m^\top\|_2\|\mb{B}_m\mb{B}_m^\top\|_2}}
		\Bigg)
		\Bigg\}+\delta\nonumber\\
		&\qquad \le2\exp\Bigg\{
		-c_1\min\Bigg(
		\frac{u^2}{K^4(\|\mb{A}\mb{A}^\top\|_F+\delta)(\|\mb{B}\mb{B}^\top\|_F+\delta)},\nonumber\\
		&\qquad\qquad\qquad\qquad\qquad
		\frac{u}{K^2\sqrt{(\|\mb{A}\mb{A}^\top\|_2+\delta)(\|\mb{B}\mb{B}^\top\|_2+\delta)}}
		\Bigg)
		\Bigg\}+\delta,
	\end{align}
	and by  \eqref{cctieq2}, \eqref{bX-bX}, and \eqref{AA_2} we obtain
	\begin{align}\label{long2}
		P\left[
		|\bd{b}^\top\bd{X}|
		\ge u
		\right]
		&\le
		P\left[
		|\bd{b}^\top\bd{X}_m|
		\ge u/2
		\right]
		+
		P\left[
		|\bd{b}^\top\bd{X}-\bd{b}^\top\bd{X}_m|
		\ge u/2
		\right]\nonumber\\
		&\le
		2\exp\left\{
		-\frac{c_2u^2}{K^2\|\bd{b}\|_F^2\| \mb{A}_m\mb{A}_m^\top \|_2}
		\right\}+\delta\nonumber\\
		&\le
		2\exp\left\{
		-\frac{c_2u^2}{K^2\|\bd{b}\|_F^2(\| \mb{A}\mb{A}^\top \|_2+\delta)}
		\right\}+\delta.
	\end{align}
	Letting $\delta\to 0$ in both \eqref{long1} and \eqref{long2}, we obtain \eqref{P_bX for subgauss} and \eqref{P_XY for subgauss} with $C=\max(c_1,c_2)$.
	
	Part (ii) can be easily shown following the same techniques for the proof of part (i) using  	
	Proposition~\ref{lem: P_bZ} with  both $\| \diag(\cdot) \|_F$ and $\left\|\offdiag(\cdot) \right\|_F$ 
	in \eqref{P_ZAZ} replaced
	by $\| \cdot \|_F$, and also using \eqref{|AB|_F} and \eqref{|bA|_F}. Details are omitted.
	
	Now consider part (iii).
	Let $\mb{C}_m\coloneqq (c_{ij})_{m\times m}=\mb{A}_m^\top\mb{B}_m$.
	By Markov's inequality and Proposition~\ref{Erdo13 Lemma}, for all $u>0$, we have
	\begin{align*}
		P\left[
		|\bd{b}^\top\bd{X}_m|
		\ge u/2
		\right]
		&\le \frac{E|\bd{b}^\top\bd{X}_m|^{k}}{(u/2)^{k}} \\
		&\le \frac{C_1^{k} k^{k/2}\eta_{k}^{k}\| \bd{b}^\top\mb{A}_m \|_F^{k}}{(u/2)^{k}}\\
		&=(2C_1\eta_{k}k^{1/2}/u)^{k} \| \bd{b}^\top\mb{A}_m \|_F^{k},
	\end{align*}	
	and
	\begin{align*}
		\lefteqn{P\left[
			\left|
			\bd{X}_m^\top\bd{Y}_m-E(\bd{X}_m^\top\bd{Y}_m)
			\right|\ge u/4
			\right]}\\
		& \qquad =
		P\left[
		\left|
		\bd{e}_m^\top\mb{C}_m\bd{e}_m-E(\bd{e}_m^\top\mb{C}_m\bd{e}_m)
		\right|\ge u/4
		\right]
		\\
		& \qquad \le
		P\bigg[
		\Big|\sum_{i=1}^m c_{ii}(e_i^2-1)\Big|\ge u/8
		\bigg]
		+
		P\bigg[
		\Big |\sum_{1\le i\ne j\le m} c_{ij}e_ie_j\Big |\ge u/8
		\bigg]\\
		& \qquad \le 
		\frac{E\Big|\sum_{i=1}^m c_{ii}(e_i^2-1)\Big|^{k/2}}{(u/8)^{k/2}}
		+\frac{E\Big |\sum_{1\le i\ne j\le m} c_{ij}e_ie_j\Big |^{k}}{(u/8)^{k}}\\
		& \qquad \le 
		\frac{C_2^{k/2}(k/2)^{k/4}\eta_{k}^{k}\| \diag(\mb{C}_m) \|_F^{k/2}}{(u/8)^{k/2}}
		+\frac{C_2^{k}  k^{k}  \eta_{k}^{2k} \| \offdiag(\mb{C}_m) \|_F^{k} }{(u/8)^{k}}\\
		& \qquad \le (4\sqrt{2}C_2 \eta_{k}^2 k^{1/2}/u  )^{k/2}  \|  \mb{A}_m^\top\mb{B}_m  \|_F^{k/2}
		+
		(8C_2\eta_{k}^2 k/u)^{k}  \|  \mb{A}_m^\top\mb{B}_m  \|_F^{k}.
	\end{align*}	
	Then by \eqref{|AB|_F} and \eqref{|bA|_F}, and using the same proof technique for part (i), we can obtain part (iii) with $C=\max(2C_1, 8C_2)$.
	Details are omitted.
\end{proof}

\begin{proof}[Proof of Lemma~\ref{lem: max diff in sample sigma}]
	From Proposition 2.7.1 in \cite{Broc91}, it is easily seen that
	$\corr(\bd{X}_{[i]})=\cov(\bd{Z}_{[i]})=\cov(\mb{P}_i\mb{H}\bd{e})=\mb{P}_i\mb{H}\cov(\bd{e}){\mb{H}}^\top{\mb{P}_i}^\top=\mb{P}_i\mb{H}{\mb{H}}^\top{\mb{P}_i}^\top$.	
	
	We first consider part (i). By \eqref{max difference in sigmas} in Proposition~\ref{lem: max diff} and part (i) of Lemma~\ref{lem: ineq infty}, for all $u>0$, we have
	\begin{align*}
		\lefteqn{P\left[
			|\widehat{\mb{\Sigma}}-\mb{\Sigma}|_{\infty}\ge 2u
			\right]}\nonumber\\	
		&\qquad  \le  2\sum_{i=1}^p\exp\left\{
		-\frac{c_1n^2u/\sigma_{ii}}{K^2\| \bd{1}_n^\top \|_F^2\|  \mb{P}_i\mb{H} \mb{H}^\top \mb{P}_i^\top\|_2}\right\}\\
		&\qquad \quad+ 2\sum_{1\le i,j\le p}\exp\bigg\{
		-c_1\min\Big(
		\frac{n^2u^2/(\sigma_{ii}\sigma_{jj})}{K^4\| \mb{P}_i\mb{H} \mb{H}^\top \mb{P}_i^\top  \|_F\| \mb{P}_j\mb{H} \mb{H}^\top \mb{P}_j^\top  \|_F},\\
		&\qquad\qquad\qquad\qquad\qquad\qquad\quad
		\frac{nu/\sqrt{\sigma_{ii}\sigma_{jj}}}{K^2\sqrt{ \| \mb{P}_i\mb{H} \mb{H}^\top \mb{P}_i^\top \|_2\|  \mb{P}_j\mb{H} \mb{H}^\top \mb{P}_j^\top \|_2   }}
		\Big)
		\bigg\}\\
		&\qquad \le 2p\exp\left\{
		-\frac{c_1nu/v_0}{K^2g_2}
		\right\}
		+2p^2\exp\left\{
		-c_1\min\left(
		\frac{nu^2/v_0^2}{K^4g_F},
		\frac{nu/v_0}{K^2 g_2}
		\right)
		\right\},	
	\end{align*}
	where $c_1>0$ is an absolute constant. Let
	\[
	u=\tilde{u}_1:=\max\Big\{
	v_0c_2K^2(\log p)g_2/n,v_0[c_2K^4(\log p )g_F/n]^{1/2}
	\Big\},
	\]
	with $c_2=(C+2)/c_1$,
	then
	\[
	P\left[
	|\widehat{\mb{\Sigma}}-\mb{\Sigma}|_{\infty}\ge 2\tilde{u}_1
	\right]
	\le 2p^{-(c_1c_2-1)}+2p^{-(c_1c_2-2)}=O(p^{-C}).
	\]
	Letting $u^*=C_1u_1$ with a constant $C_1\ge 2\max\{v_0c_2K^2, v_0 c_2^{1/2}K^2  \}$
	yields
	\[
	P\left[
	|\widehat{\mb{\Sigma}}-\mb{\Sigma}|_{\infty}\ge u^*
	\right]
	\le 
	P\left[
	|\widehat{\mb{\Sigma}}-\mb{\Sigma}|_{\infty}\ge 2\tilde{u}_1
	\right]=O(p^{-C}).
	\]

	Then consider part (ii).	
	By \eqref{max difference in sigmas} and part (ii) of Lemma~\ref{lem: ineq infty},
	for all~$u>0$, we have
	\begin{align*}
		\lefteqn{P\left[
			|\widehat{\mb{\Sigma}}-\mb{\Sigma}|_{\infty}\ge 2u
			\right]}\nonumber\\	
		& \le 
		2\sum_{i=1}^p \exp\left\{
		-\frac{(n^{1/2}u^{1/2}\sigma_{ii}^{-1/2}\| \mb{P}_i\mb{HH}^\top\mb{P}_i^\top  \|_2^{-1/2})^\frac{1}{1/2+1/\vartheta}}{c_3K(2/\vartheta)^{1/\vartheta}}
		\right\}	\\	
		&\quad+2\sum_{1\le i,j\le p} \exp\left\{
		-\frac{(\frac{nu}{\sqrt{\sigma_{ii}\sigma_{jj}}}\|  \mb{P}_i\mb{HH}^\top\mb{P}_i^\top \|_F^{-1/2}\| \mb{P}_j\mb{HH}^\top\mb{P}_j^\top \|_F^{-1/2})^{\frac{1}{1/2+2/\vartheta}}}{c_3K^2(4/\vartheta)^{4/\vartheta}}
		\right\}\\
	\end{align*}
	\begin{align*}
		&\quad+	
		2\sum_{1\le i,j\le p} \exp\left\{
		-\frac{(\frac{nu}{\sqrt{\sigma_{ii}\sigma_{jj}}}\|  \mb{P}_i\mb{HH}^\top\mb{P}_i^\top   \|_F^{-1/2}\| \mb{P}_j\mb{HH}^\top\mb{P}_j^\top \|_F^{-1/2})^{\frac{1}{1+2/\vartheta}}}{c_3K^2(2/\vartheta)^{2/\vartheta}}
		\right\}\\
		&\le 	2p\exp\left\{
		-\frac{(n^{1/2}u^{1/2}v_0^{-1/2}g_2^{-1/2})^{\frac{1}{1/2+1/\vartheta}}}{c_3K(2/\vartheta)^{1/\vartheta}}
		\right\}\\
		&\quad+2p^2\exp\left\{
		-\frac{(n^{1/2}uv_0^{-1}g_F^{-1/2})^{\frac{1}{1/2+2/\vartheta}}}{c_3K^2(4/\vartheta)^{4/\vartheta}}
		\right\}
		\\
		&\quad+2p^2\exp\left\{
		-\frac{(n^{1/2}uv_0^{-1}g_F^{-1/2})^{\frac{1}{1+2/\vartheta}}}{c_3K^2(2/\vartheta)^{2/\vartheta}}
		\right\},
	\end{align*}
	where $c_3>0$ is an absolute constant. Let 
	\begin{align*}
		u=\max\Big\{
		&v_0 [c_4K(2/\vartheta)^{1/\vartheta}\log p]^{1+2/\vartheta} g_2/n, \\
		& v_0[c_4K^2(4/\vartheta)^{4/\vartheta}\log p]^{1+2/\vartheta}(g_F/n)^{1/2}
		\Big\}
	\end{align*}
	with $c_4=\max\{K^{-2},(C+2)c_3\}$,
	we have
	\[
	P\left[
	|\widehat{\mb{\Sigma}}-\mb{\Sigma}|_{\infty}\ge 2u
	\right]\le 2p^{-(c_4/c_3-1)}+4p^{-(c_4/c_3-2)}
	=O(p^{-C}).
	\]
	Then part (ii) is established by letting $u^*=C_2u_2$ with a constant 
	$$C_2\ge 2\max\{  v_0 [c_4K(2/\vartheta)^{1/\vartheta}]^{1+2/\vartheta},v_0[c_4K^2(4/\vartheta)^{4/\vartheta}]^{1+2/\vartheta}   \}.$$

	Now consider part (iii).
	By \eqref{max difference in sigmas} and part (iii) of Lemma~\ref{lem: ineq infty}, for all $u>0$,
	there exists a constant $c_5>0$ such that
	\begin{align*}
		\lefteqn{P\left[|\widehat{\mb{\Sigma}}-\mb{\Sigma}|_{\infty}\ge 2u\right]}\\
		&\le \sum_{i=1}^p[c_5\eta_{k}k^{1/2}\sigma_{ii}^{1/2}/(nu^{1/2})]^{k}\| \mb{P}_i\mb{HH}^\top  \mb{P}_i^\top\|_2^{k/2} n^{k/2}\\
		&\quad+
		\sum_{1\le i,j\le p}[c_5\eta_{k}^2 k^{1/2}\sigma_{ii}^{1/2}\sigma_{jj}^{1/2}/(nu)]^{k/2}
		\| \mb{P}_i\mb{HH}^\top  \mb{P}_i^\top\|_F^{k/4}
		\| \mb{P}_j\mb{HH}^\top  \mb{P}_j^\top\|_F^{k/4}\\
		&\quad+
		\sum_{1\le i,j\le p}[c_5\eta_{k}^2k\sigma_{ii}^{1/2}\sigma_{jj}^{1/2}/(nu)]^{k}
		\| \mb{P}_i\mb{HH}^\top  \mb{P}_i^\top\|_F^{k/2}
		\| \mb{P}_j\mb{HH}^\top  \mb{P}_j^\top\|_F^{k/2}
		\\
		&\le p [c_5^{2}\eta_{k}^{2}kv_0g_2/(nu)]^{k/2}
		+p^2 \{
		[c_5^2\eta_k^4kv_0^2g_F/(nu^2)]^{k/4}
		{+}
		[c_5^2\eta_k^4k^2v_0^2g_F/(nu^2)]^{k/2}\},
	\end{align*}
	which is $O(p^{-C})$
	when 
	\[
	u=\max\left\{
	c_5^{2}\eta_{k}^{2}kv_0p^{(2+2C)/k}g_2/n,
	(c_5^2\eta_k^4k^2v_0^2)^{1/2}p^{(4+2C)/k}(g_F/n)^{1/2}
	\right\}.
	\]
	Then part (iii) is obtained by letting $u^*=C_3u_3$ with a constant 
	\[
	C_3\ge 2\max\left\{
	c_5^{2}\eta_{k}^{2}kv_0,
	(c_5^2\eta_k^4k^2v_0^2)^{1/2}
	\right\}.
	\]
\end{proof}

\begin{proof}[Proof of Lemma~\ref{lem: max diff in sample R}]
	Using inequality \eqref{max difference in pseudo-rhos} and following the proof of Lemma~\ref{lem: max diff in sample sigma},
	we can obtain 
	\be\label{P for max difference in pseudo-rhos}
	P\left[
	|\widetilde{\mb{R}}-\mb{R}|_{\infty}\ge u^*/4
	\right]=O(p^{-C})
	\ee
	with sufficiently large positive constants $C_1,C_2$ and $C_3$ in $u^*$ for the tail conditions  \ref{C1}, \ref{C2}, and \ref{C3}, respectively. 
	Then by Proposition~\ref{lem: max diff},
	we have
	\[
	P\left[
	|\widehat{\mb{R}}-\mb{R}|_{\infty}\ge u^*
	\right]=O(p^{-C}).
	\]
\end{proof}

\section{Proofs of Main Results}
\begin{proof}[Proof of Theorem~\ref{thm1}]
	Without loss of generality, we assume $\bd{\mu}_p=\bd{0}$ and $\mb{X}_{p\times n}$ is generated from \eqref{linear form} with $m=\infty$. 
	
	By part (i) of Lemma~\ref{lem: max diff in sample sigma},
	for any constant $c_1>0$, 
	there exists a constant $c_2>0$ depending only on $v_0,K,c_1$ such that
	when $M_1\ge c_2$, we have
	\be
	\label{max dif sigma}
	P\left[\max_{1\le i,j\le p}|\hat{\sigma}_{ij}-\sigma_{ij}|\ge \tau_1\right]=O(p^{-c_1})
	\quad \text{with}\ \tau_1=M_1u_1.
	\ee
	Then following the similar lines of the proof of Theorem 1 after  equation (12) in
	\cite{Bick08a}
	and the proof of Theorem 1 in
	\cite{Roth09},
	we obtain that for any constant $c_1>0$, there exists a constant $c_3\ge  c_2$ such that
	\begin{align}
		\label{S_t(Sigma)-Sigma matrix 2 norm}
		P\left[\|S_{\tau_1}(\widehat{\bd{\Sigma}})-\bd{\Sigma}\|_2\ge C_1c_pu_1^{1-q}\right]
		&\le P\left[\|S_{\tau_1}(\widehat{\bd{\Sigma}})-\bd{\Sigma}\|_1\ge C_1c_pu_1^{1-q}\right] \\
		& =O(p^{-c_1}), \nonumber
	\end{align}
	where $\tau_1=M_1u_1$ with any constant $M_1\ge c_3$ and some constant $C_1>0$ dependent on $M_1$. Thus, we obtain \eqref{prob 2-norm}.
	
	By condition \eqref{(iii)} of the generalized thresholding function and \eqref{max dif sigma}, we have
	\begin{align}
		\label{S_t(Sigma)-Sigma matrix max norm}
		&P\left[
		|S_{\tau_1}(\widehat{\bd{\Sigma}})-\bd{\Sigma}|_{\infty}\ge 2{\tau_1}
		\right] \nonumber \\
		&\qquad =P\left[\max_{1\le i,j\le p}|s_{\tau_1}(\hat{\sigma}_{ij})-\sigma_{ij}|\ge 2{\tau_1}\right]\nonumber\\
		&\qquad\le P\left[\max_{1\le i,j\le p}|s_{\tau_1}(\hat{\sigma}_{ij})-\hat{\sigma}_{ij}|+\max_{1\le i,j\le p}|\hat{\sigma}_{ij}-\sigma_{ij}|\ge 2{\tau_1}\right]\nonumber\\
		&\qquad\le P\left[{\tau_1}+\max_{1\le i,j\le p}|\hat{\sigma}_{ij}-\sigma_{ij}|\ge 2{\tau_1}\right]
		=O(p^{-c_1}).
	\end{align}
	Thus, $|S_{\tau_1}(\widehat{\bd{\Sigma}})-\bd{\Sigma}|_{\infty}=O_P(u_1)$.
	By \eqref{S_t(Sigma)-Sigma matrix 2 norm}, \eqref{S_t(Sigma)-Sigma matrix max norm} and the inequality $\|\mb{M}\|_F^2\le p\|\mb{M}\|_{1}|\mb{M}|_{\infty}$ for any $p\times p$ matrix $\mb{M}$, we have
	\begin{align}
		\label{P(F_n)}
		\lefteqn{P\left[\frac{1}{p}\|S_{\tau_1}(\widehat{\bd{\Sigma}})-\bd{\Sigma}\|_F^2\ge 2{\tau_1}C_1c_pu_1^{1-q}\right]}\nonumber\\
		&\quad \le P\left[\|S_{\tau_1}(\widehat{\bd{\Sigma}})-\bd{\Sigma}\|_1|S_{\tau_1}(\widehat{\bd{\Sigma}})-\bd{\Sigma}|_{\infty}\ge 2{\tau_1}C_1c_pu_1^{1-q}\right]\nonumber\\
		&\quad \le
		P\left[\|S_{\tau_1}(\widehat{\bd{\Sigma}})-\bd{\Sigma}\|_1\ge C_1c_pu_1^{1-q}\right]
		+
		P\left[
		|S_{\tau_1}(\widehat{\bd{\Sigma}})-\bd{\Sigma}|_{\infty}\ge 2{\tau_1}
		\right]
		\nonumber\\
		&\quad =O(p^{-c_1}).
	\end{align}
	Hence, we obtain \eqref{prob Frobenius norm}.

	For the convergence in mean square, we additionally assume $p\ge n^c$ for some constant $c>0$. Now
	\begin{align}
		\label{E|T(hat_Sigma-Sigma)|_2}
		&E\|S_{\tau_1}(\widehat{\bd{\Sigma}})-\bd{\Sigma}\|_2^2
		\nonumber \\
		&\qquad=E\left[\|S_{\tau_1}(\widehat{\bd{\Sigma}})-\bd{\Sigma}\|_2^2\mathds{1}\left(\|S_{\tau_1}(\widehat{\bd{\Sigma}})-\bd{\Sigma}\|_2\ge C_1c_pu_1^{1-q}\right)\right]
		\nonumber\\
		&\qquad\qquad+E\left[
		\|S_{\tau_1}(\widehat{\bd{\Sigma}})-\bd{\Sigma}\|_2^2\mathds{1}\left(\|S_{\tau_1}(\widehat{\bd{\Sigma}})-\bd{\Sigma}\|_2< C_1c_pu_1^{1-q}\right)
		\right]\nonumber\\
		&\qquad\le \left( E\|S_{\tau_1}(\widehat{\bd{\Sigma}})-\bd{\Sigma}\|_2^4\right)^{\frac{1}{2}}\left(P\left[
		\|S_{\tau_1}(\widehat{\bd{\Sigma}})-\bd{\Sigma}\|_2\ge C_1c_pu_1^{1-q}
		\right]\right)^{\frac{1}{2}} \nonumber \\
		&\qquad\qquad
		+(C_1c_pu_1^{1-q})^2\nonumber\\
		&\qquad\le\left( E\|S_{\tau_1}(\widehat{\bd{\Sigma}})-\bd{\Sigma}\|_F^4\right)^{\frac{1}{2}}O(p^{-\frac{c_1}{2}})
		+(C_1c_pu_1^{1-q})^2.
	\end{align}
	We want to show $E\|S_{\tau_1}(\widehat{\bd{\Sigma}})-\bd{\Sigma}\|_F^4=O(p^{c_5})$ with a constant $c_5>0$, and then choose a sufficiently large $c_1$ to obtain desired result.
	By condition \eqref{(iii)} of the generalized thresholding function, we have
	$
	\|S_{\tau_1}(\widehat{\bd{\Sigma}})-\widehat{\bd{\Sigma}}\|_F\le p{\tau_1},
	$
	then
	\begin{align}
		\lefteqn{E\|S_{\tau_1}(\widehat{\bd{\Sigma}})-\bd{\Sigma}\|_F^4} \nonumber\\
		&\quad \le  E
		\left(
		\|S_{\tau_1}(\widehat{\bd{\Sigma}})-\widehat{\bd{\Sigma}}\|_F+\|\widehat{\bd{\Sigma}}-\bd{\Sigma}\|_F
		\right)^4\nonumber \\
		&\quad  =E\Big[
		\|S_{\tau_1}(\widehat{\bd{\Sigma}})-\widehat{\bd{\Sigma}}\|_F^4+\|\widehat{\bd{\Sigma}}-\bd{\Sigma}\|_F^4+4\|S_{\tau_1}(\widehat{\bd{\Sigma}})-\widehat{\bd{\Sigma}}\|_F^3\|\widehat{\bd{\Sigma}}-\bd{\Sigma}\|_F\nonumber\\
		&\quad \quad+4\|S_{\tau_1}(\widehat{\bd{\Sigma}})-\widehat{\bd{\Sigma}}\|_F\|\widehat{\bd{\Sigma}}-\bd{\Sigma}\|_F^3+6\|S_{\tau_1}(\widehat{\bd{\Sigma}})-\widehat{\bd{\Sigma}}\|_F^2\|\widehat{\bd{\Sigma}}-\bd{\Sigma}\|_F^2
		\Big]\nonumber\\
		&\quad \le p^4{\tau_1}^4 + E\|\widehat{\bd{\Sigma}}-\bd{\Sigma}\|_F^4+4p^3{\tau_1}^3E\|\widehat{\bd{\Sigma}}-\bd{\Sigma}\|_F
		+4p{\tau_1}E\|\widehat{\bd{\Sigma}}-\bd{\Sigma}\|_F^3\nonumber\\
		&\quad \quad +6p^2{\tau_1}^2E\|\widehat{\bd{\Sigma}}-\bd{\Sigma}\|_F^2\nonumber\\
		& \quad \le p^4{\tau_1}^4 + E\|\widehat{\bd{\Sigma}}-\bd{\Sigma}\|_F^4+4p^3{\tau_1}^3\left(E\|\widehat{\bd{\Sigma}}-\bd{\Sigma}\|_F^2\right)^{\frac{1}{2}}
		\nonumber\\
		&\quad \quad
		+4p{\tau_1}\left(E\|\widehat{\bd{\Sigma}}-\bd{\Sigma}\|_F^6\right)^{\frac{1}{2}}
		+6p^2{\tau_1}^2E\|\widehat{\bd{\Sigma}}-\bd{\Sigma}\|_F^2.\label{E|T(hat_Sigma-Sigma)|^2}
	\end{align}
	Since $p\ge n^c$,
	it is easy to see that for $d=1,2,3$,
	\be
	\label{poly}
	\|\widehat{\bd{\Sigma}}-\bd{\Sigma}\|_F^{2d}=\left\{\sum_{1\le i,j\le p}\left(
	\frac{1}{n}\sum_{k=1}^nX_{ik}X_{jk}-\bar{X}_i\bar{X}_j-\sigma_{ij}
	\right)^2\right\}^d.
	\ee
	This is a polynomial of variables $X_{ij}$, $1\le i\le p$, $1\le j\le n$, of degree $4d$.  The number of all its terms is bounded by $p^{C_2}$ with a constant $C_2>0$,
	and all its coefficients are bounded by a constant $C_3>0$ that only depends on $v_0$.
	Denote by $P_{\ell}^{(d)}$ the $\ell$-th term of the polynomial in \eqref{poly}.
	By H\"older's inequality, there exist positive constants $c_6$ and $c_7$ for all $\ell$ and $d$ such that
	\be
	\label{Holder ineq}
	E\left(|P_{\ell}^{(d)}|\right)\le C_3\prod_{i,j} \left[E(|X_{ij}|^{c_{ij\ell d}})\right]^{\frac{1}{C_{ij\ell d}}}
	\ee
	with appropriate choices of integer constants $c_{ij\ell d}\in[0, c_6]$ and $C_{ij\ell d}\in[1,c_7]$,
	and $\sum_{i,j}\mathds{1}(c_{ij\ell d}\ne 0)\le 4d$.
	Note that it is assumed $\mu_p=0$ without loss of generality.
	First by Lyapounov's inequality,
	then by Fatou's Lemma, see Theorem~25.11 in \cite{Bill95}, and by \eqref{|bZ|_k} in Proposition~\ref{Erdo13 Lemma}, for all $k\ge 2$, we have
	\begin{align}\label{moment ineq}
		E|X_{ij}| &\le
		(E|X_{ij}|^k)^{1/k} \\
		&\le \underset{m\to\infty} {\lim\inf} \Big(  E\Big| \sum_{s=1}^m h_{i+(j-1)p,s}e_{s} \Big|^k\Big)^{1/k} \nonumber
	\end{align}
	\begin{align}
		&\le c_8 k^{1/2} \eta_k  \bigg(\sum_{s=1}^\infty h_{i+(j-1)p,s}^2\bigg)^{1/2}\nonumber\\
		&=c_8k^{1/2}\eta_k\sigma_{ii}^{1/2} \nonumber \\
		&\le c_8k^{1/2}\eta_kv_0^{1/2}\nonumber
	\end{align}
	with an absolute constant $c_8>0$ and $\eta_k=Kk^{1/2}$ by \ref{C1}.
	Hence, \eqref{Holder ineq} is bounded by a constant uniformly for all $\ell$ and $d$.
	Then by \eqref{poly},
	we have
	\be
	\label{E|hat_Sigma-Sigma|^2}
	E \|\widehat{\bd{\Sigma}}-\bd{\Sigma}\|_F^{2d}=O(p^{C_2}),\ \text{for}\ d=1,2,3.
	\ee
	Then by \eqref{E|T(hat_Sigma-Sigma)|^2},
	we obtain
	$E\|S_{\tau_1}(\widehat{\bd{\Sigma}})-\bd{\Sigma}\|_F^4=O(p^{c_5})$ with some constant $c_5>0$.
	Hence, by \eqref{E|T(hat_Sigma-Sigma)|_2}, we have
	\[
	E\|S_{\tau_1}(\widehat{\bd{\Sigma}})-\bd{\Sigma}\|_2^2\le O(p^{\frac{c_5-c_1}{2}}) +(C_1c_pu_1^{1-q})^2.
	\]

	Since $u_1\ge n^{-1}\ge p^{-1/c}$,
	we can let $c_1$ be sufficiently large such that $p^{\frac{c_5-c_1}{2}}=O\left((c_pu_1^{1-q})^2\right)$, and then we obtain
	the rate of mean-square convergence in the spectral norm.
	Similarly, from \eqref{S_t(Sigma)-Sigma matrix max norm} and \eqref{P(F_n)} we can obtain the rates of mean-square convergence in max norm and Frobenius norm, respectively.
\end{proof}

\begin{proof}[Proof of Theorem~\ref{thm2}]
	For the sparsistency and sign-consistency, the proof follows the similar lines of the proof of Theorem 2 in \cite{Roth09} by replacing their equation (A.4) with \eqref{max dif sigma}. Details are hence omitted.
\end{proof}

\begin{proof}[Proof of Corollary~\ref{cor1}]
	By using part (i) of Lemma~\ref{lem: max diff in sample R},
	the proof follows similar lines of the proof of Theorem~\ref{thm1},
	where
	we simply use
	$\|S_{\tau_2}(\widehat{\mb{R}})-\mb{R}\|_F^4\le 16p^4$ to bound the first term on the right hand side of the counterpart of
	\eqref{E|T(hat_Sigma-Sigma)|_2}.
\end{proof}

\begin{proof}[Proof of Theorem~\ref{minimax}]
	From Theorem~\ref{thm1} and Corollary~\ref{cor1},  when $u_1\le \kappa\sqrt{(\log p)/n}$ for a constant $\kappa$, the generalized thresholding estimators satisfy the upper bounds for their minimax risks. We only need to show that the minimax lower bounds have the same magnitude as the upper bounds.
	
	For the covariance matrix estimation, 
	following \cite{CZ12},
	we only need to show the desirable minimax lower bounds
	for a subset
	of $\mathcal{P}_1$, denoted by $\mathcal{P}_*(c_3,\bd{\mu}_p)$ that will be defined later.
	For any positive $v_0$ and $c_p$, there exists a constant $c_3\in (0,1)$ such that $c_3\le v_0$ and $c_p-c_3^q>0$.
	We consider the set of covariance matrices $\mathcal{F}_*$
	defined in (20) of \cite{CZ12}.
	The set $\mathcal{F}_*$ depends on a parameter $c_{n,p}$ satisfying
	$0\le c_{n,p}\le c_4 n^{(1-q)/2}(\log p)^{-(3-q)/2}$ with some given constant $c_4>0$. If $\mb{\Sigma}\in\mathcal{F}_*$, then $\mb{\Sigma}$ is positive definite, $\sigma_{ii}=1$ for all $i$, and
	\[
	\max_{1\le j\le p}\sum_{i\ne j}|\sigma_{ij}|^q\le c_{n,p}.
	\]
	Define $\mathcal{F}_{*}^{c_3}=\{c_3\mb{\Sigma}:\mb{\Sigma}\in \mathcal{F}_{*} \}$.
	If $\mb{\Sigma}\in \mathcal{F}_*^{c_3}$, then
	$\sigma_{ii}=c_3\le v_0$ for all $i$, and
	\[
	\max_{1\le j\le p}\sum_{i=1}^p|\sigma_{ij}|^q\le c_3^q(c_{n,p}+1).
	\]
	Let $c_p=c_3^q(c_{n,p}+1)$ and $c_4=c_3^{-q}c_2$, then
	$0<c_3^{-q}(c_p-c_3^{q})=c_{n,p}\le c_3^{-q}c_p\le c_3^{-q}c_2n^{(1-q)/2}(\log p)^{-(3-q)/2}
	=c_4n^{(1-q)/2}(\log p)^{-(3-q)/2}
	$, thus
	$\mathcal{F}_*^{c_3}\subset \mathcal{U}(q,c_p,v_0)$.
	Let $\mathcal{P}_*(c_3,\bd{\mu}_p)$ denote the set of distributions of
	$\bd{X}_{pn}=\text{vec}(\mb{X}_{p\times n})$ satisfying that the columns of $\mb{X}_{p\times n}$, i.e., $\bd{X}_1,\dots,\bd{X}_n$ are i.i.d. each following a $p$-variate Gaussian distribution with mean $\bd{\mu}_p$ and covariance matrix $\mb{\Sigma}\in \mathcal{F}_{*}^{c_3}$.
	Then $\bd{X}_{pn}$ can be generated in the form of \eqref{linear form} with
	a linear combination of standard Gaussian variables for which $K=K_G$. That $\bd{X}_{pn}$ can have $g_F=g_2=1$ together with $\sqrt{(\log p)/n}=o(1)$ leads to $u_1=\sqrt{(\log p)/n}$ for sufficiently large $n$.
	Hence,
	$
	\mathcal{P}_*(c_3,\bd{\mu}_p)\subset \mathcal{P}_1
	$
	when $n$ is sufficiently large.
	From the proofs of Theorems 2 and 4 in \cite{CZ12}, we have
	\[
	\inf_{\widetilde{\mb{\Sigma}}} \sup_{\mathfrak{D}\in\mathcal{P}_*(1,\bd{0})}E_{\bd{X}_{pn}|\mathfrak{D}}\left(\| \widetilde{\mb{\Sigma}}-\mb{\Sigma} \|_2^2\right)
	\ge c_5c_{n,p}^2\left(\frac{\log p}{n}\right)^{1-q}
	\]
	and
	\[
	\inf_{\widetilde{\mb{\Sigma}}} \sup_{\mathfrak{D}\in\mathcal{P}_*(1,\bd{0})}\frac{1}{p}E_{\bd{X}_{pn}|\mathfrak{D}}\left(\| \widetilde{\mb{\Sigma}}-\mb{\Sigma} \|_F^2\right)
	\ge c_5c_{n,p}\left(\frac{\log p}{n}\right)^{1-q/2}
	\]
	with some constant $c_5>0$.
	By the argument of changing variables,
	\begin{align*}
		&\inf_{\widetilde{\mb{\Sigma}}}
		\sup_{\mathfrak{D}\in \mathcal{P}_1}
		E_{\bd{X}_{pn}|\mathfrak{D}}\left(\| \widetilde{\mb{\Sigma}}-\mb{\Sigma} \|_2^2\right)\\
		&\qquad \ge
		\inf_{\widetilde{\mb{\Sigma}}}
		\sup_{\mathfrak{D}\in\mathcal{P}_*(c_3,\bd{\mu}_p)}
		E_{\bd{X}_{pn}|\mathfrak{D}}\left(\| \widetilde{\mb{\Sigma}}-\mb{\Sigma} \|_2^2\right)\\
		&\qquad=
		c_3^2\inf_{\widetilde{\mb{\Sigma}}}
		\sup_{\mathfrak{D}\in\mathcal{P}_*(1,\bd{0})}
		E_{\bd{X}_{pn}|\mathfrak{D}}\left(\| \widetilde{\mb{\Sigma}}-\mb{\Sigma} \|_2^2\right)\\
		&\qquad \ge
		c_3^2 c_5c_{n,p}^2\left(\frac{\log p}{n}\right)^{1-q}\\
		&\qquad=
		c_3^2 c_5[c_3^{-q}(c_p-c_3^{q})]^2\left(\frac{\log p}{n}\right)^{1-q}
		\asymp c_p^2\left(\frac{\log p}{n}\right)^{1-q},
	\end{align*}
	and similarly,
	\begin{align*}
		\inf_{\widetilde{\mb{\Sigma}}}
		\sup_{\mathfrak{D}\in \mathcal{P}_1}\frac{1}{p}
		E_{\bd{X}_{pn}|\mathfrak{D}}\left(\| \widetilde{\mb{\Sigma}}-\mb{\Sigma} \|_F^2\right)
		&\ge
		c_3^2\inf_{\widetilde{\mb{\Sigma}}}
		\sup_{\mathfrak{D}\in\mathcal{P}_*(1,\bd{0})}\frac{1}{p}
		E_{\bd{X}_{pn}|\mathfrak{D}}\left(\| \widetilde{\mb{\Sigma}}-\mb{\Sigma} \|_F^2\right)\\
		&\asymp
		c_p\left(\frac{\log p}{n}\right)^{1-q/2}.
	\end{align*}
	
	For the correlation matrix estimation, when $c_p>1$, we have
	\begin{align*}
		\inf_{\widetilde{\mb{R}}}\sup_{\mathfrak{D}\in\mathcal{P}_2}E_{\bd{X}_{pn}|\mathfrak{D}}\left(\| \widetilde{\mb{R}}-\mb{R}  \|_2^2 \right)
		&\ge 
		\inf_{\widetilde{\mb{\Sigma}}} \sup_{\mathfrak{D}\in\mathcal{P}_*(1,\bd{0})}E_{\bd{X}_{pn}|\mathfrak{D}}\left(\| \widetilde{\mb{\Sigma}}-\mb{\Sigma} \|_2^2\right)\\
		&\ge c_5c_{n,p}^2\left(\frac{\log p}{n}\right)^{1-q} \\
		&=c_5(c_p-1)^2\left(\frac{\log p}{n}\right)^{1-q} \\
		&\asymp  c_p^2\left(\frac{\log p}{n}\right)^{1-q},
	\end{align*}
	and similarly,
	\begin{align*}
		\inf_{\widetilde{\mb{R}}}\sup_{\mathfrak{D}\in\mathcal{P}_2}
		\frac{1}{p}
		E_{\bd{X}_{pn}|\mathfrak{D}}\left(\| \widetilde{\mb{R}}-\mb{R}  \|_F^2 \right)
		&\ge 
		\inf_{\widetilde{\mb{\Sigma}}} \sup_{\mathfrak{D}\in\mathcal{P}_*(1,\bd{0})}\frac{1}{p}
		E_{\bd{X}_{pn}|\mathfrak{D}}\left(\| \widetilde{\mb{\Sigma}}-\mb{\Sigma} \|_F^2\right)\\
		&\asymp c_p\left(\frac{\log p}{n}\right)^{1-q/2}.
	\end{align*}
	
	The proof is now complete.
\end{proof}

\begin{proof}[Proofs of Theorems~\ref{clime_thm1} and \ref{clime_thm2}]
	From part (i) of Lemma~\ref{lem: max diff in sample sigma}, we have \eqref{max dif sigma}, then the proofs of the convergence rates for CLIME follow similarly to
	the proofs of Theorems 2, 5 and 6 in \cite{CLL11}, 
	where we also obtain $|\widehat{\bd{\Omega}}_{\varepsilon,\lambda_1}-\bd{\Omega}|_{\infty}\le 4M_p\lambda_1$ with probability tending to 1.
	Also, the proofs of sparsistency and sign-consistency for the thresholded CLIME estimator follow the similar lines of the proof of Theorem 2 in \cite{Roth09}. Details are hence omitted.
\end{proof}

\begin{proof}[Proof of Theorem~\ref{spice_thm2}]
	Since 
	$$\min_{i}\sigma_{ii}-\varphi_{\min}(\mb{\Sigma})=\min_{i}\bd{e}_i^\top(\mb{\Sigma}-\varphi_{\min}(\mb{\Sigma})\mb{I}_{p\times p})\bd{e}_i\ge 0$$ and $\max_i\sigma_{ii}\le |\mb{\Sigma}|_{\infty}\le \|\mb{\Sigma} \|_2=\varphi_{\max}(\mb{\Sigma})$,
	we have
	\be
	\label{bound for sigma}
	v_0^{-1}\le \min_i\sigma_{ii}\le \max_i\sigma_{ii}\le v_0.
	\ee
	Thus,
	\begin{eqnarray}
		\label{ineq for W}
		&& v_0^{-1/2}\le \| \mb{W} \|_2,\|  \mb{W}^{-1} \|_2 \le v_0^{1/2},\\
		\label{ineq for K}
		&& \| \mb{K} \|_2\le \|  \mb{W} \|_2 \| \mb{\Omega}  \|_2 \|  \mb{W} \|_2 \le v_0^{2},\\
		&& \|\mb{R}\|_2\le \|  \mb{W}^{-1} \|_2 \| \mb{\Sigma}  \|_2 \|  \mb{W}^{-1} \|_2 \le v_0^{2},
	\end{eqnarray}
	and
	\be
	\label{ineq for R}
	v_0^{-2}\le\| \mb{K} \|_2^{-1}=\varphi_{\min}(\mb{R})\le \varphi_{\max}(\mb{R})=\|\mb{R}\|_2\le v_0^{2}.
	\ee
	Under $\max_i\sigma_{ii}\le v_0$ and $u_1=o(1)$, 
	we can obtain
	part (i) of Lemma~\ref{lem: max diff in sample sigma},
	part (i) of Lemma~\ref{lem: max diff in sample R},
	and \eqref{P for max difference in pseudo-rhos} with $u^*=C_1u_1$ for some constant $C_1>0$.
	In other words,
	for any constant $C'>0$, there exists a constant $C_1>0$ 
	depending only on $v_0,K$ and $C'$
	such that with
	probability $1-O(p^{-C'})$,
	\begin{eqnarray}
		\label{diff in sigma}
		&& \max_{1\le i,j\le p}|\hat{\sigma}_{ij}-\sigma_{ij}|\le C_1{u_1},\\
		\label{diff in rho}
		&& \max_{1\le i,j\le p}|\hat{\rho}_{ij}-\rho_{ij}|\le C_1{u_1},
	\end{eqnarray}
	and
	\be
	\label{max diff in pseu rho}
	\max_{1\le i,j\le p}\left|\frac{\hat{\sigma}_{ij}-\sigma_{ij}}{\sqrt{\sigma_{ii}\sigma_{jj}}} \right|\le C_1{u_1}.
	\ee
	From \eqref{diff in sigma} and \eqref{bound for sigma},  we obtain 
	$\max_{i} \hat{\sigma}_{ii}^{-1/2}\le 2v_0^{1/2}$ with probability $1-O(p^{-C'})$.
	Letting $i=j$ in \eqref{max diff in pseu rho}, with probability $1-O(p^{-C'})$, we have 
	\begin{align}
		\label{infty diff in inv W}
		o(1)=2C_1{u_1} v_0^{1/2}
		&\ge C_1{u_1}\max_{1\le i\le p}\hat{\sigma}_{ii}^{-\frac{1}{2}}
		\ge\max_{1\le i\le p}\left|\frac{\hat{\sigma}_{ii}-\sigma_{ii}}{\sigma_{ii}} \right| \max_{1\le i\le p}\hat{\sigma}_{ii}^{-\frac{1}{2}} \nonumber\\
		&\ge \max_{1\le i\le p}\left| \sqrt{\frac{\hat{\sigma}_{ii}}{\sigma_{ii}}}-1  \right| \left| \sqrt{\frac{\hat{\sigma}_{ii}}{\sigma_{ii}}}+1  \right|\hat{\sigma}_{ii}^{-\frac{1}{2}}
		\ge \max_{1\le i\le p}\left| \sqrt{\frac{\hat{\sigma}_{ii}}{\sigma_{ii}}}-1  \right|\hat{\sigma}_{ii}^{-\frac{1}{2}} \nonumber\\
		&= \max_{1\le i\le p}| \sigma_{ii}^{-\frac{1}{2}} -\hat{\sigma}_{ii}^{-\frac{1}{2}}    |
		=\| \widehat{\mb{W}}^{-1} -\mb{W}^{-1} \|_2,
	\end{align}
	and then by \eqref{ineq for W},
	\be
	\label{2-norm for est W}
	\| \widehat{\mb{W}}^{-1} \|_2\le \| \widehat{\mb{W}}^{-1} -\mb{W}^{-1} \|_2+ \|\mb{W}^{-1} \|_2= o(1)+v_0^{1/2}.
	\ee
	
	Now recall the assumption that ${u_1}=o(1/\sqrt{1+s_p})$.
	Following similar lines of the proof of Theorem 1 in \cite{Roth08}
	by replacing their line 10 on page 500 by $r_n={u_1}\sqrt{s_p}\to 0$,
	replacing their  line 5 on page 501 by \eqref{diff in rho},
	replacing their inequality (14) by $\text{II}=0$, replacing their equation (15) by
	$\lambda_2=C_1{u_1}/\varepsilon$ with a sufficiently small constant $\varepsilon>0$,
	and replacing  the last line on their page 501 by
	$|\Delta_S^-|_1\le\sqrt{s_p}\|\Delta^-\|_F$, as well as using \eqref{ineq for R} to establish the counterpart of their inequality (18) for $\mb{K}$,
	we obtain
	\be\label{diff spice corr}
	\|\widehat{\mb{K}}_{\lambda_2}-\mb{K}\|_F=O_P(u_1\sqrt{s_p}).	
	\ee


	From the proof of Theorem 2 in \cite{Roth08}, we have
	\begin{align}
		\label{omega diff}
		\lefteqn{\| \widehat{\mb{\Omega}}_{\lambda_2}-\mb{\Omega} \|_2}\nonumber\\
		&\quad\le
		\| \widehat{\mb{K}}_{\lambda_2}- \mb{K}   \|_2
		(\| \widehat{\mb{W}}^{-1} -\mb{W}^{-1} \|_2^2 +\| \widehat{\mb{W}}^{-1} \|_2  \|\mb{W}^{-1} \|_2)
		\nonumber\\
		&\qquad+
		\| \widehat{\mb{W}}^{-1} -\mb{W}^{-1} \|_2
		( \| \widehat{\mb{K}}_{\lambda_2} \|_2  \| \mb{W}^{-1} \|_2   +\|\mb{K}\|_2 \| \widehat{\mb{W}}^{-1} \|_2)
		\nonumber\\
		&\quad\le
		\| \widehat{\mb{K}}_{\lambda_2}- \mb{K}   \|_F
		(\| \widehat{\mb{W}}^{-1} -\mb{W}^{-1} \|_2^2 +\| \widehat{\mb{W}}^{-1} \|_2  \|\mb{W}^{-1} \|_2)
		\nonumber\\
		&\qquad +
		\| \widehat{\mb{W}}^{-1} -\mb{W}^{-1} \|_2
		\left[( \| \widehat{\mb{K}}_{\lambda_2}-\mb{K} \|_F + \|\mb{K}\|_2)\| \mb{W}^{-1} \|_2   +\|\mb{K}\|_2 \| \widehat{\mb{W}}^{-1} \|_2\right].
	\end{align}
	Plugging \eqref{diff spice corr}, \eqref{infty diff in inv W}, \eqref{2-norm for est W}, \eqref{ineq for W} and \eqref{ineq for K} into \eqref{omega diff} yields $\| \widehat{\mb{\Omega}}_{\lambda_2}-\mb{\Omega}  \|_2=O_P(u_1\sqrt{1+s_p})$.
	
	We can obtain $\frac{1}{\sqrt{p}}\| \widehat{\mb{\Omega}}_{\lambda_2}-\mb{\Omega}  \|_F=O_P\left(u_1\sqrt{1+s_p/p}\right)$ similarly from
	\begin{align*}
		\label{omega diff F}
		\lefteqn{\| \widehat{\mb{\Omega}}_{\lambda_2}-\mb{\Omega} \|_F}\nonumber\\
		&=\|\widehat{\mb{W}}^{-1} \widehat{\mb{K}}_{\lambda_2} \widehat{\mb{W}}^{-1}-\mb{W}^{-1}\mb{K}\mb{W}^{-1} \|_{F}\nonumber\\
		&=\|(\widehat{\mb{W}}^{-1} -\mb{W}^{-1})(\widehat{\mb{K}}_{\lambda_2}- \mb{K} )(\widehat{\mb{W}}^{-1} -\mb{W}^{-1})
		+\mb{W}^{-1}\widehat{\mb{K}}_{\lambda_2}(\widehat{\mb{W}}^{-1} -\mb{W}^{-1})\nonumber\\
		&\qquad+(\widehat{\mb{W}}^{-1} -\mb{W}^{-1})\mb{K}\widehat{\mb{W}}^{-1}
		+\widehat{\mb{W}}^{-1} (\widehat{\mb{K}}_{\lambda_2}-\mb{K})   \mb{W}^{-1}\|_F\nonumber\\
		&\le
		\| \widehat{\mb{K}}_{\lambda_2}- \mb{K}   \|_F
		(\| \widehat{\mb{W}}^{-1} -\mb{W}^{-1} \|_2^2 +\| \widehat{\mb{W}}^{-1} \|_2  \|\mb{W}^{-1} \|_2)
		\nonumber\\
		&\quad+
		\| \widehat{\mb{W}}^{-1} -\mb{W}^{-1} \|_F
		( \| \widehat{\mb{K}}_{\lambda_2} \|_2  \| \mb{W}^{-1} \|_2   +\|\mb{K}\|_2 \| \widehat{\mb{W}}^{-1} \|_2)\\
		&\le
		\| \widehat{\mb{K}}_{\lambda_2}- \mb{K}   \|_F
		(\| \widehat{\mb{W}}^{-1} -\mb{W}^{-1} \|_2^2 +\| \widehat{\mb{W}}^{-1} \|_2  \|\mb{W}^{-1} \|_2)
		\nonumber\\
		&\quad+
		\sqrt{p}\| \widehat{\mb{W}}^{-1} -\mb{W}^{-1} \|_2
		\left[( \| \widehat{\mb{K}}_{\lambda_2}-\mb{K} \|_F + \|\mb{K}\|_2)\| \mb{W}^{-1} \|_2   +\|\mb{K}\|_2 \| \widehat{\mb{W}}^{-1} \|_2\right],
	\end{align*}
	where $\|\mb{BA}\|_F=\|\mb{AB}\|_F\le \|\mb{A}\|_2\|\mb{B}\|_F$ for symmetric matrices $\mb{A}$ and $\mb{B}$ [see  Lemma 1 in \cite{Lam09}].
	
\end{proof}

\begin{proof}[Proof of Theorem~\ref{spice_thm3}]
	First, we need to show $|\widehat{\mb{K}}_{\lambda_2}-\mb{K}|_{\infty}=o_P(1)$,
	which is similar to the proof of Theorem 1 in \cite{Ravi11}. We follow some of their notation for convenience.
	In the proof,  their $\Theta$ and $\Sigma$ are now replaced by our $\mb{K}$ and $\mb{R}$, respectively.
	But we keep their $W$ that is our $\widehat{\mb{R}}-\mb{R}$, which should not be confused with our $\mb{W}$ in bold.
	From \eqref{diff in rho}, for any constant $\tau>2$ [note that here we use the notation $\tau$ given in \cite{Ravi11} 
	rather than the one defined as the generalized thresholding parameter], there exist constants $M_1$ and $N_1$ such that when
	$M\ge M_1$ and
	$n>N_1$, we have
	\be
	\label{replace lemma 8}
	P(|W|_{\infty}\le Mu_1)\ge P(|W|_{\infty}\le M_1u_1)=1- O(p^{2-\tau}),
	\ee
	thus we can set their $\bar{\delta}_f(n,p^\tau)=M{u_1}$ and their $1/v_*=\infty$.
	Then, $\lambda_2=8M{u_1}/\beta=8\bar{\delta}_f(n,p^\tau)/\beta$.
	From $\lambda_2
	\le[6(1+\beta/8)d\max\{\kappa_{\mb{R}} \kappa_{\mb{\Gamma}}, \kappa_{\mb{R}}^3  \kappa_{\mb{\Gamma}}^2   \}]^{-1}$, we have
	\be
	\label{delta<}
	\bar{\delta}_f(n,p^\tau)
	\le[6(1+8/\beta)d\max\{\kappa_{\mb{R}} \kappa_{\mb{\Gamma}}, \kappa_{\mb{R}}^3  \kappa_{\mb{\Gamma}}^2   \}]^{-1}.
	\ee
	Then following the proof of their Theorem 1 by using \eqref{replace lemma 8} instead of their Lemma~8, and \eqref{delta<} instead of their (15) and (29),
	with probability $1-O(p^{2-\tau})$ we have
	\be
	\label{infty diff in K}
	|\widehat{\mb{K}}_{\lambda_2}-\mb{K}|_{\infty}\le 2(1+8/\beta)\kappa_{\mb{\Gamma}}\bar{\delta}_f(n,p^\tau)
	=2(1+8/\beta)\kappa_{\mb{\Gamma}}M{u_1}=o(1),
	\ee
	and
	all entries of $\widehat{\mb{K}}_{\lambda_2}$ in $S^c$ are zero.
	By $|\mb{BA}|_{\infty}=|\mb{AB}|_{\infty}\le |\mb{A}|_{\infty}\| \mb{B} \|_1$ for symmetric matrices $\mb{A}$ and $\mb{B}$,
	we have
	\begin{align}
		\label{spice infty diff in Omega}
		\lefteqn{|\widehat{\mb{\Omega}}_{\lambda_2}-\mb{\Omega}|_{\infty}}\nonumber\\
		&= |\widehat{\mb{W}}^{-1} \widehat{\mb{K}}_{\lambda_2} \widehat{\mb{W}}^{-1}-\mb{W}^{-1}\mb{K}\mb{W}^{-1} |_{\infty}\nonumber\\
		&=|(\widehat{\mb{W}}^{-1} -\mb{W}^{-1})(\widehat{\mb{K}}_{\lambda_2}- \mb{K} )(\widehat{\mb{W}}^{-1} -\mb{W}^{-1})
		+\mb{W}^{-1}\widehat{\mb{K}}_{\lambda_2}(\widehat{\mb{W}}^{-1} -\mb{W}^{-1})\nonumber\\
		&\qquad+(\widehat{\mb{W}}^{-1} -\mb{W}^{-1})\mb{K}\widehat{\mb{W}}^{-1}
		+\widehat{\mb{W}}^{-1} (\widehat{\mb{K}}_{\lambda_2}-\mb{K})   \mb{W}^{-1}|_{\infty}\nonumber\\
		&\le| \widehat{\mb{K}}_{\lambda_2}- \mb{K}   |_{\infty}  \| \widehat{\mb{W}}^{-1} -\mb{W}^{-1} \|_1^2
		+ | \widehat{\mb{K}}_{\lambda_2} |_{\infty}  \| \mb{W}^{-1} \|_1    \| \widehat{\mb{W}}^{-1} -\mb{W}^{-1} \|_1\nonumber\\
		&\qquad +|\mb{K}|_{\infty}  \| \widehat{\mb{W}}^{-1} -\mb{W}^{-1} \|_1  \| \widehat{\mb{W}}^{-1} \|_1
		+| \widehat{\mb{K}}_{\lambda_2}- \mb{K} |_{\infty}\| \widehat{\mb{W}}^{-1} \|_1  \|\mb{W}^{-1} \|_1 \nonumber\\
		&=| \widehat{\mb{K}}_{\lambda_2}- \mb{K}   |_{\infty}  \| \widehat{\mb{W}}^{-1} -\mb{W}^{-1} \|_2^2
		+ | \widehat{\mb{K}}_{\lambda_2} |_{\infty}  \| \mb{W}^{-1} \|_2    \| \widehat{\mb{W}}^{-1} -\mb{W}^{-1} \|_2\nonumber\\
		&\qquad +\|\mb{K}\|_2  \| \widehat{\mb{W}}^{-1} -\mb{W}^{-1} \|_2  \| \widehat{\mb{W}}^{-1} \|_2
		+| \widehat{\mb{K}}_{\lambda_2}- \mb{K} |_{\infty}\| \widehat{\mb{W}}^{-1} \|_2  \|\mb{W}^{-1} \|_2.
	\end{align}
	By inequalities \eqref{ineq for K} and \eqref{infty diff in K}, with probability $1-O(p^{2-\tau})$ we have
	\be
	\label{infty norm for est K}
	| \widehat{\mb{K}}_{\lambda_2} |_{\infty}
	\le |\mb{K}|_{\infty}+| \widehat{\mb{K}}_{\lambda_2}-\mb{K} |_{\infty}
	\le\|\mb{K}\|_2+ | \widehat{\mb{K}}_{\lambda_2}-\mb{K} |_{\infty}
	\le v_0^{2}+o(1).
	\ee
	Plugging \eqref{infty diff in K}, \eqref{infty diff in inv W}, \eqref{infty norm for est K}, \eqref{ineq for W}, \eqref{ineq for K}, \eqref{2-norm for est W}
	into \eqref{spice infty diff in Omega} and letting $M\ge\max\{M_1,10C_1v_0^{2}\}$ yields that, with probability $1-O(p^{2-\tau})$,
	\begin{align}
		\label{bound for infty diff in Omega}
		\lefteqn{|\widehat{\mb{\Omega}}_{\lambda_2}-\mb{\Omega}|_{\infty}}\nonumber\\
		&\quad \le 2(1+8/\beta)\kappa_{\mb{\Gamma}}M{u_1}o(1)
		+\left(v_0^{2}+o(1)\right)  v_0^{1/2} 2C_1{u_1} v_0^{1/2}\nonumber\\
		&\qquad+ v_0^{2}  2C_1{u_1} v_0^{1/2}  \left( o(1)+v_0^{1/2}   \right)
		+2(1+8/\beta)\kappa_{\mb{\Gamma}}M{u_1}  \left( o(1)+v_0^{1/2}   \right) v_0^{1/2}\nonumber\\
		&\quad \le 5C_1{u_1}v_0^{3}
		+2.5(1+8/\beta)\kappa_{\mb{\Gamma}}M{u_1}  v_0
		\le\left(
		0.5+2.5(1+8/\beta)\kappa_{\mb{\Gamma}}
		\right)M{u_1}  v_0
		=r, 
	\end{align}
	\begin{align*}
		\|\widehat{\mb{\Omega}}_{\lambda_2}-\mb{\Omega}  \|_2
		&\le \min\{
		\|\widehat{\mb{\Omega}}_{\lambda_2}-\mb{\Omega}  \|_1,
		\|\widehat{\mb{\Omega}}_{\lambda_2}-\mb{\Omega}  \|_F
		\} \\
		&\le \min\{d,\sqrt{p+s_p}\}|\widehat{\mb{\Omega}}_{\lambda_2}-\mb{\Omega}|_{\infty} \\
		& \le \min\{d,\sqrt{p+s_p}\} r, \nonumber
	\end{align*}
	and
	\begin{align*}
		p^{-1/2}\|\widehat{\mb{\Omega}}_{\lambda_2}-\mb{\Omega}  \|_F
		&\le \min\{
		\|\widehat{\mb{\Omega}}_{\lambda_2}-\mb{\Omega}  \|_2,
		p^{-1/2}\sqrt{p+s_p} |\widehat{\mb{\Omega}}_{\lambda_2}-\mb{\Omega}|_{\infty}
		\}\nonumber\\
		&\le  r\min\left\{d,\sqrt{1+s_p/p}\right\} \\
		& = r \sqrt{1+s_p/p},
	\end{align*}
	where the last equality follows from $\sqrt{1+s_p/p} \le \sqrt{1+(d-1)p/p} = \sqrt{d} \le d$. For any $(i,j)\in S$, by \eqref{bound for infty diff in Omega} and $|\omega_{ij}|>r$, $\hat{\omega}_{ij}^{(\lambda_2)}$ cannot differ enough from the nonzero $\omega_{ij}$ to change sign.
	Since $\widehat{\mb{\Omega}}_{\lambda_2}$ has the same sparsity as $\widehat{\mb{K}}_{\lambda_2}$ and we have shown that, with probability $1-O(p^{2-\tau})$, all entries of $\widehat{\mb{K}}_{\lambda_2}$ in $S^c$ are zero, then
	$\widehat{\mb{\Omega}}_{\lambda_2}$ also has this sparsistency result.
\end{proof}

\begin{proof}[Proofs of Theorems~\ref{subexp thm} and \ref{poly thm}]
	The proofs are similar to the proofs of the preceding theorems and corollary by using
	the corresponding results given in Lemmas~\ref{lem: max diff in sample sigma}
	and \ref{lem: max diff in sample R} for conditions~\ref{C2} and \ref{C3}, respectively.
	Details are omitted.
\end{proof}

\section{Some Details for Numerical Considerations}

\subsection{Candidate values for tuning parameters}\label{AppB:tuning}
Let $\eta$ be the generic notation for each considered tuning parameter, i.e., $\eta$ represents one of $\tau_1$, $\tau_2$, $\lambda_1$, and $\lambda_2$ corresponding to the generalized thresholding covariance matrix estimation, the generalized thresholding correlation matrix estimation, CLIME, and SPICE, respectively. The ordered candidate values $\eta_1,\dots,\eta_N$ of $\eta$ are chosen from a logarithmic spaced grid. Specifically, $\log\eta_1,\dots,\log\eta_N$ are equally spaced values, where $\eta_1=r\eta_N$ with some constant  $r\in(0,1)$. In our numerical examples, we use $N=50$ and $r=0.01$.

For the generalized thresholding estimation of a correlation matrix, we let $\eta_N$ be the largest absolute off-diagonal value of the sample correlation matrix so that the thresholding estimator using $\eta_N$ is a diagonal matrix.
For the covariance matrix estimation, we first scale the sample covariance matrix to the sample correlation matrix, and then re-scale the generalized thresholding estimator of the correlation matrix back to the original scale using the sample standard deviations.

For CLIME, we use the same $\eta_N$ generated by the R package {\tt flare}
[\cite{Li15}, version 1.5.0; see the function {\tt sugm}]
based on the following formula
\[
\eta_N=\mathds{1}(\eta^*\ne 0)\eta^*+\mathds{1}(\eta^*=0)\eta^{**}
\]
with
\begin{align*}
	\eta^*&=\min\left\{ \max_{1\le i,j\le p} s_{ij}, \ -\min_{1\le i,j\le p} s_{ij}  \right\},\\
	\eta^{**}&=\max\left\{ \max_{1\le i,j\le p} s_{ij}, \ -\min_{1\le i,j\le p} s_{ij}  \right\},
\end{align*} 
where $(s_{ij})_{p\times p} \coloneqq \hat{\mb{\Sigma}}-{\rm diag} \{ \hat{\sigma}_{11},\dots,\hat{\sigma}_{pp}  \}$.
For SPICE, we generate $\eta_N$ using the same approach implemented in the R package {\tt huge} [\cite{Zhao12}, version 1.2.7; see the function {\tt huge.glasso}] for GLasso, i.e., 
$\eta_N$ is the largest absolute off-diagonal value of the sample correlation matrix.
Note that SPICE is a slight modification of GLasso.

\subsection{Data generating method in simulations}
\begin{figure}[tb!]
	\begin{center}
		\includegraphics[width=3.9in]{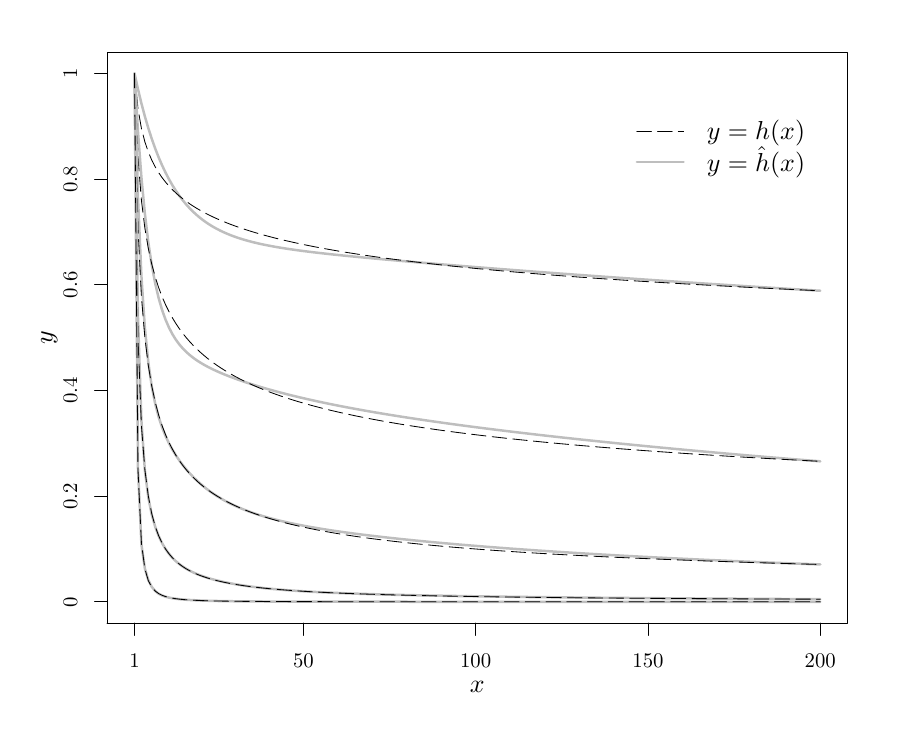}
	\end{center}
	\caption{Approximation of $h(x)=x^{-\alpha}$ for \mbox{$\alpha=0.1,0.25,0.5,1,2$}.
		\label{expfit}}
\end{figure}

It is computationally expensive to simulate data $\bd{X}_{pn}\coloneqq \textup{vec}(\mb{X}_{p\times n})$ directly from a multivariate Gaussian random number generator because of the large dimension of its covariance matrix $\textup{Cov}(\bd{X}_{pn})$.
We use an alternative way
to simulate the data that approximately satisfy
\eqref{simu rate matrix}.
Note that
$h(x)=x^{-\alpha}$ with $x\in[1,n]$ and $\alpha>0$ can be approximated by $\hat{h}(x)=\sum_{i=0}^Na_i\exp(-b_ix)$ with small $N$ and appropriately chosen $\{a_i, b_i\}$ by the method of \cite{Boch07}
[see Figure~\ref{expfit}].
Thus, data $\mb{X}_{p\times n}$ are simulated as follows:
each column of $\mb{X}_{p\times n}$ is generated by
$\bd{X}_t=\sum_{i=0}^Nc_i\bd{Y}_t^{(i)}$ for $t=1,...,n,$
where $c_i=\sqrt{a_i\exp(-b_i)}$,
$\bd{Y}_1^{(i)}$ are i.i.d. $\mathcal{N}(\bd{0},\mb{\Sigma})$ for all $i$,
and for $t=2,...,n$,
$\bd{Y}_{t}^{(i)}=\rho_i \bd{Y}_{t-1}^{(i)}+\bd{e}_{t}^{(i)}$
with $\rho_i=\exp(-b_i)$
and white noise
$(1-\rho_i^2)^{-1/2}\bd{e}_t^{(i)}$ i.i.d. $\mathcal{N}(\bd{0},\mb{\Sigma})$. It is easily seen that
$
\mb{R}^{t,t+j}
=\sum_{i=0}^Nc_i^2  \rho_i^j \mb{R}
=\sum_{i=0}^N a_i \exp\{-b_i(j+1)\}\mb{R}
\approx (j+1)^{-\alpha}\mb{R}.
$

\section{Additional Results of the rfMRI Data Analysis}

\begin{figure}
			\centering
			\begin{tabular}{c c}
				\subfloat[]{
					\includegraphics[scale=0.4]{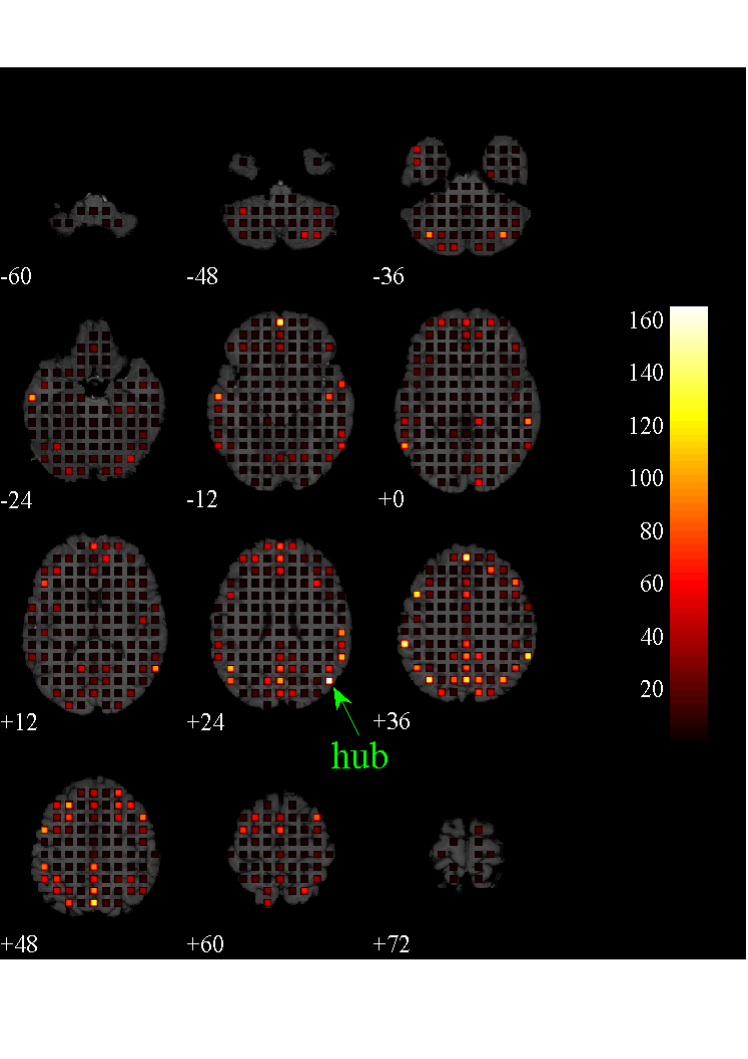}}
				&
				\subfloat[]{
					\includegraphics[scale=0.4]{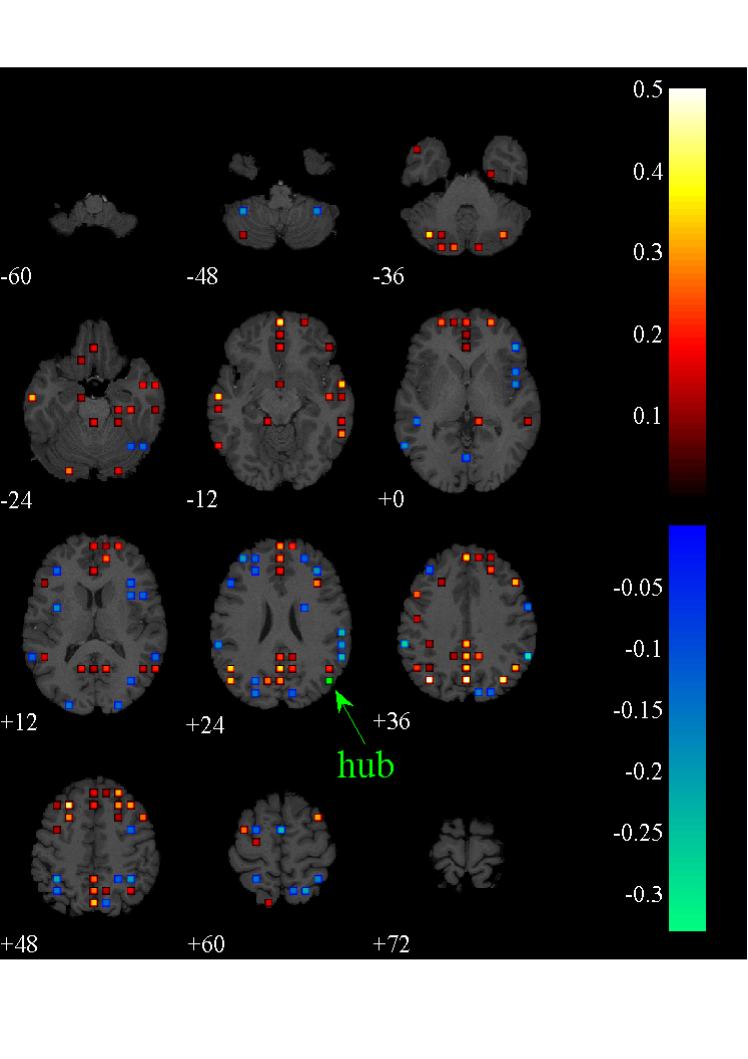}}
				\\
				\subfloat[]{
					\includegraphics[scale=0.4]{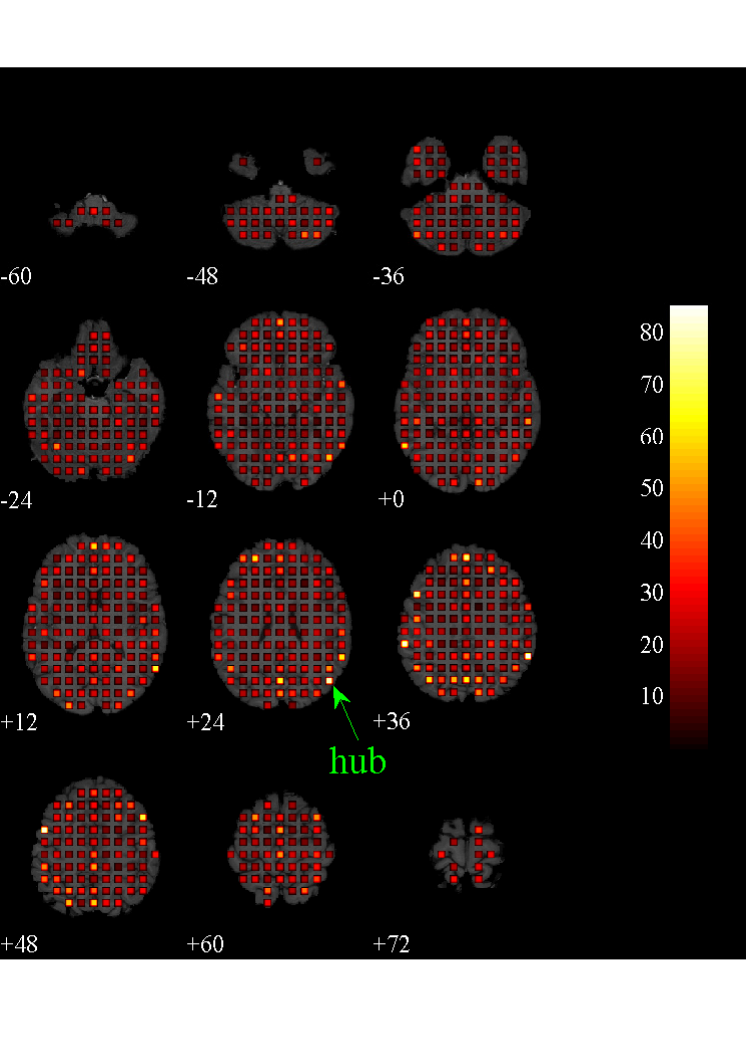}}
				&
				\subfloat[]{
					\includegraphics[scale=0.4]{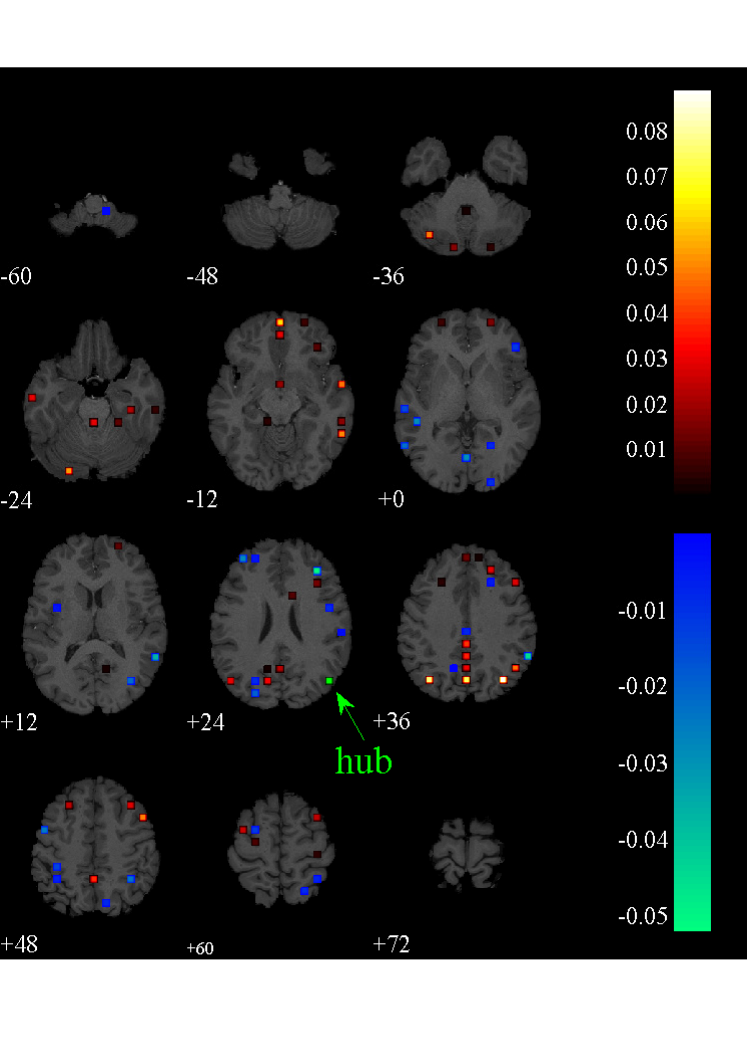}}
			\end{tabular}
			\caption{rfMRI data analysis for brain functional connectivity.
				(a) Node degrees of marginal connectivity found by hard thresholding. (b) Marginally connected nodes and their estimated correlations to the selected hub. (c) Node degrees of direct connectivity found by CLIME.
				(d) Directly connected nodes and their estimated partial correlations to the selected hub. The brain is plotted in the Montreal Neurological Institute 152 space
				with $Z$-coordinates displayed.}
			\label{fcon}
		\end{figure}

We illustrate the node degrees of marginal connectivity and direct connectivity in
Figure~\ref{fcon} (a) and (c), respectively.	
The top 10 hubs for marginal connectivity and the top 10 hubs for direct connectivity are listed in the following two tables.
The coordinates of the center of each hub is given in the Montreal Neurological Institute (MNI) 152 space. The hubs for marginal connectivity with MNI coordinates listed in bold numbers are spatially close to those found in \cite{Buck09} and \cite{Cole10} from studies with multiple subjects.
As an illustration, we plot the marginal connectivity and the direct connectivity of a single hub in
Figure~\ref{fcon} (b) and (d), respectively, which is ranked No. 1 in degree of marginal connectivity and No. 4 in degree of direct connectivity.
This hub has 164 marginally connected nodes and 79 directly connected nodes, where 80\% of the directly connected nodes are also marginally connected. It is located in the right inferior parietal cortex, a part of the so-called default mode network [\cite{Buck08}] that is most active during the resting state.

\begin{table}[h]
	\caption{Top 10 hubs for marginal connectivity found by hard thresholding}
	\centering
	\scalebox{0.9}{
		\begin{tabular}{l l r c c c}
			\hline\hline
			Rank & Location & MNI coordinates & Degree
			& Direct rank & Direct degree\\[0.5ex]
			\hline
			1& Inferior parietal             &\textbf{48, -72, 24}   & 164 & 4 & 79\\
			2& Supramarginal              &-60, -36, 36  & 151 & 3 & 82 \\
			3& Superior frontal            &\textbf{0, 48, 36}       & 150 & 6 & 73 \\
			4&Medial orbitofrontal       &\textbf{0, 60, -12}& 140 &20 & 53\\
			5& Inferior parietal             &\textbf{-36, -72, 36} & 137 &15 &61\\
			6&Supramarginal               &\textbf{60, -48, 36} & 131 &1 &85\\
			7&Precuneus                     &0, -72, 48 & 128 &16 & 58\\
			8& Precuneus                    &0, -72, 36 & 125 &10 & 64\\
			9&Rostral middle frontal    &\textbf{-48, 12, 36} & 121 & 5&74 \\
			10&Inferior parietal            &\textbf{-48, -60, 24} & 109& 37 & 48\\
			\hline
		\end{tabular}}
	\end{table}

	\begin{table}[h]
		\caption{Top 10 hubs for direct connectivity found by CLIME}
		\centering
		\scalebox{0.9}{
			\begin{tabular}{l l r c c c}
				\hline\hline
				Rank & Location & MNI coordinates & Degree
				& Marginal rank & Marginal degree\\[0.5ex]
				\hline
				1& Inferior parietal             & 60, -48, 36          & 85 & 6 & 131\\
				2& Precentral            &-48, 0, 48                    & 82 & 18 & 98 \\
				3& Supramarginal            &-60, -36, 36          & 82 & 2 & 151 \\
				4&Inferior parietal      & 48, -72, 24                & 79 &1 & 164\\
				5& Rostral middle frontal        &-48, 12, 36      & 74 &9 &121\\
				6&Superior frontal               &0, 48, 36            & 73 &3 &150\\
				7&Caudal middle frontal       &48, 12, 48           & 68 &29 & 87\\
				8& Middle temporal               &60, -60, 12          & 66 &19 & 96\\
				9&Precuneus                    &0, -72, 24                 & 65 & 14&101 \\
				10&Precuneus           &0, -72, 36                        & 64& 8 & 125\\
				\hline
			\end{tabular}}
		\end{table}

\bibliographystyle{imsart-number}
\bibliography{reference}

\end{document}